\documentclass[12pt]{amsart}
\usepackage{amsfonts}
\usepackage{fullpage,multirow, amsmath, amsthm,amsfonts,amssymb,stmaryrd, mathrsfs,hhline}
\usepackage{mathdots}
\usepackage{pgf,tikz,tikz-cd}
\usepackage{mathrsfs}
\usepackage{url}
\usepackage{mathtools}

\usetikzlibrary{arrows}
\usetikzlibrary[patterns]
\usetikzlibrary{matrix,decorations.pathreplacing,calc,shapes, decorations.pathmorphing}
\usepackage{hyperref}

\usepackage[all]{xy}

\usepackage[all]{xy}

\newtheorem{theorem}[subsection]{Theorem}
\newtheorem{lemma}[subsection]{Lemma}
\newtheorem{corollary}[subsection]{Corollary}
\newtheorem{conjecture}[subsection]{Conjecture}

\newtheorem{proposition}[subsection]{Proposition}
\theoremstyle{definition}

\newtheorem{definition}[subsection]{Definition}
\newtheorem{definition-proposition}[subsection]{Definition-Proposition}

\newtheorem{example}[subsection]{Example}

\newtheorem{remark}[subsection]{Remark}
\newtheorem{convention}[subsection]{Convention}
\newtheorem{notation}[subsection]{Notation}

\numberwithin{equation}{subsection}

\makeatletter

\newcommand{\Rmnum}[1]{\expandafter\@slowromancap\romannumeral #1@}
\makeatother

\def\calA{\mathcal{A}}

\def\calC{\mathcal{C}}
\def\calD{\mathcal{D}}

\def\calL{\mathcal{L}}
\def\calM{\mathcal{M}}
\def\calO{\mathcal{O}}

\def\calS{\mathcal{S}}

\def\calW{\mathcal{W}}

\def\gothm{\mathfrak{m}}

\def\CC{\mathbb{C}}
\def\FF{\mathbb{F}}
\def\GG{\mathbb{G}}

\def\QQ{\mathbb{Q}}
\def\RR{\mathbb{R}}

\def\ZZ{\mathbb{Z}}

\def\bff{\mathbf{f}}

\def\scrL{\mathscr{L}}

\def\scrS{\mathscr{S}}

\def\a{\alpha}
\def\b{\beta}

\def\d{\delta}
\def\e{\varepsilon}

\def\V{\tilde{V}}

\def\bi{\dbinom}

\DeclareMathOperator{\Gal}{Gal}

\DeclareMathOperator{\NP}{NP}

\DeclareMathOperator{\Spf}{Spf}
\DeclareMathOperator{\Spm}{Spm}
\DeclareMathOperator{\Sym}{Sym}

\DeclareMathOperator{\diag}{diag}

\DeclareMathOperator{\ur}{ur}
\DeclareMathOperator{\Berk}{Berk}
\DeclareMathOperator{\Vtx}{Vtx}

\newcommand{\quash}[1]{}  

\newcommand{\new}{\mathrm{new}}
\newcommand{\old}{\mathrm{old}}
\newcommand{\Iw}{\mathrm{Iw}}
\newcommand{\unr}{\mathrm{unr}}

\newcommand{\rig}{\mathrm{rig}}

\newcommand{\wt}{\mathrm{wt}}

\DeclareMathOperator{\GL}{GL}

\DeclareMathOperator{\Spc}{Spc}

\DeclareMathOperator{\I}{I}
\DeclareMathOperator{\M}{M}

\DeclareMathOperator{\Dig}{Dig}

\DeclareMathOperator{\CS}{CS}
\DeclareMathOperator{\GZ}{GZ}
\DeclareMathOperator{\gl}{gl}
\DeclareMathOperator{\sgn}{sgn}
\DeclareMathOperator{\BV}{BV}

\begin{document}
	\title{Slope Distributions of $\mathcal{L}$-Invariants for Modular Forms}
	\author{Jiawei An} 
	\address{Jiawei An, Morningside Center of Mathematics,
		Chinese Academy of Sciences, 55 Zhongguancun East Road, Haidian District, Beijing 100190,
		China.}
	\email{jiaweian@amss.ac.cn}
	
	\begin{abstract}
		Fix a prime $p\geq5$ and an integer $N\geq1$ coprime to $p$. In this paper, we use ghost series to study the $p$-adic valuations of $\calL$-invariants associated with modular forms. More precisely, let $\bar r:\Gal(\bar\QQ/\QQ)\rightarrow \GL_2(\bar{\FF}_p)$ be an irreducible representation. Assuming that $\bar r$ is locally reducible and very generic, (1) we determine the slopes of $\calL$-invariants associated with $p$-newforms of weight $k$, level $\Gamma_0(Np)$, and residual representation $\bar r$, with at most $O(\log_pk)$ exceptions; (2) we prove the integrality of these slopes; and (3) we prove an equidistribution property for these slopes as the weight $k$ tends to infinity, which confirms a distribution conjecture proposed by Bergdall--Pollack (\cite{BP dis}) recently (as an analogue of Gouvêa’s conjecture on slope distributions). 
	\end{abstract}
	
	\maketitle
	
	\setcounter{tocdepth}{1}
	\tableofcontents
	\section{Introduction} 
	Let $p$ be a prime number, and $N\geq1$ be an integer coprime to $p$. In this paper, we prove an equidistribution property (Theorem \ref{Main theorem on BP dis}) for the $p$-adic valuations of $\mathcal{L}$-invariants attached to normalized $p$-new eigenforms of level $\Gamma_0(Np)$, under certain local assumptions on the associated Galois representation. This confirms specific instances of the following Conjecture \ref{BP conjecture}, which was proposed by Bergdall--Pollack \cite{BP dis} recently.  As a byproduct of our analysis, we also obtain an integrality result for these slopes (Theorem \ref{cor: integrality inroduction}). 
	
	 Let $f\in S_k(\Gamma_0(Np))$ be a \emph{$p$-newform}, namely a normalized eigenform that is new at $p$ (Similarly, a \emph{$p$-oldform} is a normalized eigenform that is old at $p$). Denote the associated $\calL$-invariant by $\calL_f\in\bar{\QQ}_p$. By the \emph{slope} of $\calL_f$, we mean the $p$-adic valuation $v_p(\calL_f)$. In this paper, the $p$-adic valuation $v_p(-)$ is normalized so that $v_p(p) =1$.
	\subsection{Equidistribution conjectures}
	Our primary objective is to verify an equidistribution property for the slopes of $\calL$-invariants associated with $p$-newforms. We first recall a well-known analogous equidistribution conjecture for the \emph{slopes of modular forms}.
	
For a positive integer $M$,	let $f\in S_k(\Gamma_0(M))$ be a normalized eigenform with $p$-th Fourier coefficient $a_p(f)$. The \emph{($p$-adic) slope} of $f$ is the $p$-adic valuation $v_p(a_p(f))$. For instance, when $M=Np$ and $f$ is a $p$-newform, it is well-known that $a_p(f)\in\left\{\pm p^{\frac{k-2}{2}}\right\}$, and hence the slope of $f$ is $\frac{k-2}{2}$. For the slopes of modular forms of level $N$, Gouv\^ea proposed the following distribution conjecture (\cite{gouvea}).
		\begin{conjecture}[Gouv\^ea]
		\label{Gouvea dis}
	
		For a real number $T$, consider the following set of ``normalized" slopes of modular forms 
		\[
		X_T=\left\{\frac{p+1}{k}\cdot v_p(a_p(f)): f\ \textrm{is a normalized eigenform in}\ S_k(\Gamma_0(N)),k\leq T\right\}.
		\]
		Then $X_T$ is equidistributed on $[0,1]$ as $T$ tends to infinity.
	\end{conjecture}
	
A local version of this conjecture has been proved in \cite{Lghost2} recently; see \S~\ref{subsection:ghost} for further discussions. In this paper, we aim to address the following Bergdall--Pollack distribution conjecture about the slopes of $\calL$-invariants associated with $p$-newforms.
	
	\begin{conjecture}[\cite{BP dis}]
		\label{BP conjecture}
		For a real number $T$, consider the following set of ``normalized" slopes of $\calL$-invariants associated with $p$-newforms 
			\begin{equation}
			Y_{T}:=\left\{\frac{2(p+1)}{(p-1)k}\cdot v_p(\calL^{-1}_f): f\ \textrm{is a normalized}\ p\textrm{-new eigenform in}\ S_k(\Gamma_0(Np)),k\leq T\right\}.
		\end{equation}
	Then $Y_{T}\subseteq(-\infty,+\infty]$ is evenly distributed on $[0,1]$ as $T$ tends to infinity.
	
	\end{conjecture}

	Notably, both Conjectures \ref{Gouvea dis} and \ref{BP conjecture} are supported by a substantial amount of numerical evidence (\cite{gouvea,BP dis}). Beyond the computational data, one may see some conceptual link between the two conjectures by interpreting them in terms of equidistribution properties for the slopes of parameters of \emph{modular} local Galois representations. 
	
The local Galois representation of $\Gal(\bar{\QQ}_p/\QQ_p)$ attached to a $p$-oldform $f\in S_k(\Gamma_0(Np))$ is \emph{crystalline}, and hence it is determined by the weight $k$ and the $p$-th Fourier coefficient $a_p(f)$. Similarly, for a $p$-newform $f\in S_k(\Gamma_0(Np))$, the associated local Galois representation is \emph{semi-stable but non-crystalline}, and it is determined by the weight $k$, the $p$-th Fourier coefficient $a_p(f) \in \left\{\pm p^{\frac{k-2}{2}}\right\}$, and the $\calL$-invariant $\mathcal{L}_f \in \bar{\QQ}_p$. Therefore, for a fixed weight $k\geq2$, both $a_p(f)$ in Conjecture \ref{Gouvea dis} and $\calL_f$ in Conjecture \ref{BP conjecture} can be viewed as parameters of certain local Galois representations arising from modular forms. This observation also clarifies our interest in investigating the slopes of $\mathcal{L}$-invariants associated with $p$-newforms. We refer to \cite[\S~6]{BP dis} for further discussions. See \cite[\S~8.2]{BP dis} for a discussion about the relationship between these two scalars $\frac{p+1}{k}$ and $\frac{2(p+1)}{(p-1)k}$ appearing in the ``normalizations" of two distribution conjectures.

	\subsection{Local assumptions and  main results}
	Fix a prime number $p\geq5$. Let $E$ be a sufficiently large finite extension of $\QQ_p$ with ring of integers $\calO$ and residue field $\FF$. For $\bar\a\in\FF^{\times}$, write $\unr(\bar\a)$ for the unramified character of $\Gal(\bar \QQ_p/\QQ_p)$ sending the geometric Frobenius to $\bar\a$. Let $\omega_1$ be the \emph{first fundamental character} of the inertial subgroup $\I_{\QQ_p}$.

	Let $\bar{r}:\Gal(\bar\QQ/\QQ) \to \GL_2(\FF)$ be an irreducible Galois representation. Let $S_k(\Gamma_0(Np))_{\bar r }$ be the space of modular forms of weight $k$, level $\Gamma_0(Np)$, and localized at the Hecke maximal ideal corresponding to $\bar r$. Put $\bar\rho=\bar r|_{\Gal(\bar\QQ_p/\QQ_p)}$.

	\begin{definition}
		\label{wild generic}
		The representation $\bar r$ is called \emph{locally reducible and very generic}, if
		\begin{equation}
			\bar\rho=\bar r|_{\Gal(\bar\QQ_p/\QQ_p)}\simeq \begin{pmatrix}\unr(\bar\a)\omega_{1}^{a+b+1}&*\\0&\unr(\bar\b)\omega_1^b\end{pmatrix},
		\end{equation}
where $a\in\{2,\dots,p-5\}$, $b\in\{0,\dots,p-2\}$, and $\bar\a,\bar\b\in\FF^{\times}$.
	\end{definition}
	Using the local ghost series (\S\ref{Section2}), we explicitly determine
	all but at most $O(\log_p k)$ slopes of $\calL$-invariants associated with
	$p$-newforms in $S_k(\Gamma_0(Np))_{\bar r}$; see
	Remark~\ref{rmk in the introduction of constans} for an effective bound on
	the number of exceptional slopes.
	As a consequence, we obtain two main results of this paper, concerning
	the equidistribution and integrality of these slopes.

	\begin{theorem}[Equidistribution property]
	\label{Main theorem on BP dis}
		Fix $p\geq11$. Assume that $\bar r$ is locally reducible and very generic with $\det(\bar r|_{\I_{\QQ_p}})=\omega_1^c$, for some $c\in\{0,\dots,p-2\}$. Then Conjecture \ref{BP conjecture} holds for $\bar r$. More precisely, for $k\equiv c+1\pmod{p-1}$, let $\mu_{k,\bar r}$ be the uniform probability measure of the multiset
	\begin{equation}
		Y_k:=\left\{\frac{2(p+1)}{(p-1)k}\cdot v_p(\calL^{-1}_f): f\ \textrm{is a normalized}\ p\textrm{-new eigenform in}\ S_k(\Gamma_0(Np))_{\bar r}\right\}.
	\end{equation}
	Then the measure $\mu_{c+1+n(p-1),\bar r}$ weakly converges to the uniform probability measure on the interval $[0,1]$ as $n$ tends to infinity.
\end{theorem}
	
	\begin{theorem}[Integrality]
		\label{cor: integrality inroduction}		
		Fix $p\geq11$. Assume that $\bar r$ is locally reducible and very generic. 
		Let $k\geq2$ be an integer. With at most $O(\log_p k)$ exceptions, each $p$-newform $f\in S_k(\Gamma_0(Np))_{\bar r}$ satisfies the following properties.
		\begin{enumerate}
			\item If $k$ is even, then $v_p(\calL_f)\in\ZZ$.
			\item If $k$ is odd, then $v_p(\calL_f)\in\frac{1}{2}\ZZ$.
		\end{enumerate}
	\end{theorem}
		When $k$ is odd, we refer to Corollary \ref{cor: integrality} and Proposition \ref{nearStprop} for a precise characterization of the cases for which $v_p(\calL_f)$ is not an integer. 
	\begin{remark}
		\label{rmk in the introduction of constans}
	The implied constant in the $O(\log_p k)$ term appearing in
	Theorem~\ref{cor: integrality inroduction} is effective and depends only on
	$p$, $N$, $\bar r$, and
	$a\in\{2,\dots,p-5\}$.
	More precisely, write
	$
	k=k_0+n(p-1)$,
	where \(k_0\in\{0,\dots,p-2\}\).
	With the convention \(\lfloor \log_p(0)\rfloor=0\), the
	number of exceptional forms in Theorem~\ref{cor: integrality inroduction}
	is bounded by \(m(\bar r)\,C(k,a)\), where
	\[
	C(k,a)
	=
	\frac{
		2(p+3-\min\{a+2,p-1-a\})
		+
		4\left\lfloor
		\log_p\!\left(\frac{k-k_0}{p-1}\right)
		\right\rfloor
	}{p-1}
	\le
	2+
	\frac{
		4\left\lfloor
		\log_p\!\left(\frac{k-k_0}{p-1}\right)
		\right\rfloor
	}{p-1}.
	\]
	Here \(C(k,a)\) is the constant appearing in
	Corollary~\ref{cor of thresholds}(3), and \(m(\bar r)\) is characterized
	by the property that
	\[
	\dim S_k(\Gamma_0(Np))_{\bar r}
	-
	\frac{2k}{p-1}m(\bar r)
	\]
	remains bounded as \(k\to\infty\).
	For explicit values of \(m(\bar r)\), see
	Remark~\ref{rmk global multiplicity} and the references therein.	
	Moreover, in Example~\ref{ex4}, this bound is not sharp for
	\(k=6,12,\) and \(18\), in which cases there are no exceptions,
	whereas it is sharp for \(k=24\); see
	Remark~\ref{rmk of log bound meaning} and
	Theorem~\ref{thm:global-thresholds-and-L-invariants}.
	\end{remark}
	\begin{remark}
Theorem \ref{cor: integrality inroduction} can be interpreted as a result on the integrality of the slopes of parameters (i.e. $\calL$-invariants) of certain \emph{modular} semi-stable deformations of $\bar\rho=\bar r|_{\Gal(\bar\QQ_p/\QQ_p)}$. An analogous statement for (not necessarily modular) crystabelline deformations of  $\bar\rho$ is the Breuil--Buzzard--Emerton conjecture. For further discussions of this conjecture and its proof (under certain assumptions on $\bar\rho$), see \cite[\S~1.8]{Lghost2}. 

		In a sense, Theorem~\ref{cor: integrality inroduction} is optimal, since there exist examples of $p$-newforms (of even weight $k$) whose associated \emph{local} Galois representations are reducible, for which $v_p(\mathcal{L}_f)$ is \emph{not} an integer (\cite[\S~9.3]{BP dis}). See \cite[Question~9.1]{BP dis} for a revised question concerning the integrality of the slopes of $\mathcal{L}$-invariants associated with semi-stable deformations of a reducible residual local Galois representation $\bar \rho$.
	\end{remark}

	\subsection{Main input: The ghost conjecture}
	\label{subsection:ghost}
	The proof of Theorems \ref{Main theorem on BP dis} and \ref{cor: integrality inroduction} adheres to the same method used to reduce the Gouvêa distribution Conjecture \ref{Gouvea dis} to the ghost conjecture in \cite{BP2}, namely, analyzing the \emph{ghost series}. The ghost series was introduced 
	and studied by Bergdall--Pollack in a series of papers \cite{BP1,BP2,BP3}. Roughly speaking, it serves as a model for the genuine characteristic power series of the $U_p$-action, and this combinatorial power series unified key conjectures about slopes of modular forms into the \emph{ghost conjecture}, which states that all $U_p$-slopes are given by this ghost series. In particular, Conjecture \ref{Gouvea dis} holds true assuming the validity of the ghost conjecture (\cite{BP2}).
	
	Our approach in this paper relies essentially on the \emph{local ghost theorem} (Theorem \ref{local ghost theorem}), as established in \cite{Lghost2}, which confirms the $\bar r$-\emph{ghost conjecture} for $\bar r$'s that are locally reducible and very generic. In a sense, the techniques used in this paper may be interpreted as deriving implications from the $\bar{r}$-ghost conjecture to the $\bar{r}$-Bergdall--Pollack distribution conjecture.
	
	\begin{remark}
		For the remaining $\bar{r}$'s and for the case in a global context, both the ($\bar r$-) ghost conjecture and the ($\bar r$-) Bergdall--Pollack distribution conjecture remain open questions. 
	\end{remark}

	\subsection{Overview of the proof} 
	\label{overview of the proof}
	In the remaining part of this introduction, we sketch the proof of Theorems \ref{Main theorem on BP dis} and \ref{cor: integrality inroduction} in two key steps (\S~\ref{step1} and \S~\ref{step2}). Fix a locally reducible and very generic $\bar r$ with $\det(\bar r|_{\I_{\QQ_p}})=\omega_1^c$, for some $c\in\{0,\dots,p-2\}$. We begin by recalling several notations and definitions related to \emph{overconvergent modular forms}.
	
Let $\Delta$ be the torsion subgroup of $\ZZ_p^{\times}$ and $\omega:\Delta\simeq\FF_p^{\times}\rightarrow \ZZ_p^{\times}$ be the Teichm\"uller character. Identify $\calO[\![1+p\ZZ_p]\!]\cong\calO[\![w]\!]$ via sending $[\a]$ for $\a\in1+p\ZZ_p$ to $(1+w)^{\log(\a)/p}$, where $\log(-)$ is the formal $p$-adic logarithm. 
Consider the \emph{weight space} $\calW$, whose $\CC_p$-points parameterize the continuous characters from $\ZZ_p^{\times}$ to $\CC_p^{\times}$. It decomposes into a disjoint union of $(p-1)$ connected components $\calW=\bigsqcup_{i=0}^{p-2}\calW_i$, and each $\calW_i:=(\Spf\calO[\![w]\!])^{\rig}$ is the \emph{weight disc} that parameterizes continuous characters $\chi$ of $\ZZ_p^{\times}$ such that $\chi|_{\Delta}=\omega^i$.
	If $k\geq 2$ is an integer such that $k\equiv i\pmod {p-1}$, then $k$ corresponds to the point $
	 w_{k}=\exp(p(k-2))-1
	$ (see \S~\ref{sunsection weight space}) in $\calW_i$.
	For each weight $w\in\calW_i$, let  $S_{w}^{\dagger}(\Gamma_0(Np))_{\bar r}$ be the space of \emph{overconvergent modular forms} of weight $w$, level $\Gamma_0(Np)$, and localized at the Hecke maximal ideal corresponding to $\bar r$ (when $w=w_k$, we retain the classical notation $S_{k}^{\dagger}(\Gamma_0(Np))_{\bar r}$).
	Since $\det(\bar r|_{\I_{\QQ_p}})=\omega_1^c$, this space is non-zero only if $i=\{c+1\}$, where $\{c+1\}$ is the unique integer in $\left\{0,\dots,p-2\right\}$ such that $\{c+1\}\equiv c+1 \pmod{p-1}$. 
	
	Define the \emph{slopes of overconvergent modular forms} (of weight $w$) to be the $p$-adic valuations of the eigenvalues for the $U_p$-action on $S_{w}^{\dagger}(\Gamma_0(Np))_{\bar r}$. This notion coincides with the usual slope of a (classical) normalized eigenform when $w=w_k$, via the inclusion $S_{k}(\Gamma_0(Np))_{\bar r}\subseteq S_{k}^{\dagger}(\Gamma_0(Np))_{\bar r}$. The slopes in $S_{w}^{\dagger}(\Gamma_0(Np))_{\bar r}$ are given by the slopes of the \emph{Newton polygon} $\NP(w)$, namely, the Newton polygon of the \emph{characteristic power series} $C_{\bar r}(w,t)$ of the $U_p$-action on $S_{w}^{\dagger}(\Gamma_0(Np))_{\bar r}$.
	
	Let $k\geq 2$ be an integer such that $k\equiv c+1 \pmod{p-1}$, hence $w_k\in\calW_{\{c+1\}}$.
Denote the subspaces of $S_k(\Gamma_0(Np))_{\bar r}$, that consist of $p$-oldforms and $p$-newforms, by $S_k(\Gamma_0(Np))_{\bar r}^{\old}$ and $S_k(\Gamma_0(Np))_{\bar r}^{\new}$, respectively. Moreover, let $d_k$, $2d_{k,0}$, and $d_{k,1}$ denote the dimensions of the space $S_k(\Gamma_0(Np))_{\bar r}$, $S_k(\Gamma_0(Np))_{\bar r}^{\old}$, and $S_k(\Gamma_0(Np))_{\bar r}^{\new}$, respectively. Then $2d_{k,0}+d_{k,1}=d_k$.
\subsection*{Main objects in the proof}
 Before turning to the detailed explanations of the proofs, we briefly
summarize the central objects involved and their logical relationships in
the following diagram \eqref{diag1}. Let \(m(\bar r)\) be the integer defined in
Remark~\ref{rmk in the introduction of constans}.
The main input of this paper is the local ghost theorem,
Theorem~\ref{local ghost theorem}, which allows us to transfer explicit
slope information obtained from the \emph{local ghost series} (Definition~\ref{def: local ghost series}) to
\(C_{\bar r}(w,t)\).
\begin{equation}
	\label{diag1}
	\begin{tikzcd}
		& \boxed{\text{slopes of \(A_1\)}}
		\arrow[r, leftrightarrow, "{\text{(1)}}", "{\text{GS}}"']
		& \boxed{\text{slopes of \(\mathcal{L}\)-invariants}} \\
		\boxed{\text{\(k\)-derivative polygon \(\underline{\Delta}_k\)}} 
		\arrow[ur,dashed,"{\text{(2)}}"] 
		\arrow[dr,dashed,swap,"{\text{(*)}}"]
		& & \\
		& \boxed{\text{global \(k\)-thresholds}}
	\end{tikzcd}
\end{equation}
\begin{itemize}
	\item The polygon \(\underline{\Delta}_k\) (Definition~\ref{defofDelta}), which has length \(d_{k,1}/m(\bar r)\), is explicitly determined by the local ghost series attached to \(\bar r\). Its equidistribution and integrality follow directly from the study of local ghost series; see Proposition~\ref{nearStprop} and Corollary~\ref{cor:equidistribution of Delta slope}.
	
	\item In (1), using the Greenberg--Stevens method, one may transfer the slopes of the \(\mathcal{L}\)-invariants associated to $p$-newforms in \(S_k(\Gamma_0(Np))_{\bar r}^{\new}\) to the slopes of a \(d_{k,1} \times d_{k,1}\) matrix \(A_1\), which, roughly speaking, is the first coefficient in the Taylor expansion of \(U_p^2\) near weight $k$. See \S~\ref{step1} for details.
	
	\item The main technical part is (2), where we prove that \emph{almost all} slopes of \(A_1\) are given by the \(k\)-derivative polygon \(\underline{\Delta}_k\). This result follows from the study of the variation of slopes of \(p\)-newforms; see Definition~\ref{def:introduction-(k,w_*)-newslope}. In particular, we use both an analytic trick and a linear algebra trick to deduce (2) from these variations. See \S~\ref{step2} for details.
	
	\item The implication (*) is a byproduct of the study of the variation of slopes of \(p\)-newforms; see \S~\ref{subsection in the introduction: global $k$-thresholds} for details.
\end{itemize}
The main results, Theorems~\ref{Main theorem on BP dis} and~\ref{cor: integrality inroduction}, follow from combining (2) with the integrality and equidistribution results for the slopes of $\underline{\Delta}_k$.

We conclude this summary by explaining the precise role of the local
ghost series in our argument.
The local ghost series enters through the study of the
$k$-derivative polygon $\underline{\Delta}_k$ and of the variation of slopes of
$p$-newforms. These provide the explicit slope data used in the proof.
By contrast, the Greenberg--Stevens argument in (1), as well as the
analytic and linear-algebraic arguments used in (2), are independent of
the local ghost series. 

\subsection{Greenberg--Stevens methods in higher rank}
\label{step1}
In \S~\ref{section of Muitiple GS formula}, we establish an analogous Greenberg--Stevens formula that links $\calL$-invariants to certain derivatives of the eigenvalues of the $U_p^ 2$-action. 

Let $\calC_{N,\bar r}$ be the $\bar r$-\emph{Coleman--Mazur eigencurve} of tame level $N$, and $\wt:\calC_{N,\bar r}\rightarrow \calW$ be the \emph{weight map}. Let $\calD':=\Spm \calA'\subseteq\calW_{\{c+1\}}$ be a sufficiently small affinoid closed disc centered at weight $w_k$. Then we may construct a finite free $\calA'$-module $\calS'$ of rank $d_{k,1}$, equipped with an $\calA'$-linear action of $U_p$, such that for any weight $w \in \calD'$, the slopes of the $U_p$-action on $\calS'_w$ (specialized at weight $w$) are all equal to $\frac{k-2}{2}$. Consequently, the $U_p$-action on $\calS'$ can be expressed as a $d_{k,1}\times d_{k,1}$ \emph{diagonal} matrix, whose entries are rigid analytic functions in $\calA'$. Denote the Taylor expansion of the operator $U_p^2$ at $w=w_k$ by
\begin{equation}
	M(w)=U_p^2(w)=p^{k-2}I_{d}+A_1(w-w_k)+A_2(w-w_k)^2+\dots,
\end{equation}
where each matrix $A_i\in\M_{d_{k,1}\times d_{k,1}}(\CC_p)$. Define the \emph{slopes of $A_1$} to be the multiset of the $p$-adic valuations of the eigenvalues of $A_1$. Using the Greenberg--Stevens method, we prove the following result.
\begin{theorem}[Theorem \ref{GS formula}]
	\label{thm step2}
	The $\calL$-invariants associated with $p$-newforms in $S_k(\Gamma_0(Np))_{\bar r}$ are precisely the eigenvalues of the following matrix
	\[
	-p^{2-k}A_1\cdot \left(p\cdot\exp(p(k-2))I_{d_{k,1}}\right).
	\]
	In particular, the multiset of the slopes of $A_1$ coincides with the following multiset:
	\begin{equation}
		\left\{v_p(\calL_f)+k-3:f \ \textrm{is a}\ p\textrm{-newform in}\  S_k(\Gamma_0(Np))_{\bar r}\right\}.
	\end{equation}
\end{theorem}
\begin{remark}
	Theorem~\ref{thm step2} states that the $\mathcal{L}$-invariants under
	consideration are precisely the eigenvalues of the matrix
	$
	A_1 A_0^{-1}\cdot \bigl(p\exp(p(k-2))I_{d_{k,1}}\bigr).
	$
	The matrix may be viewed as defining a linear operator on an abstract
	$d_{k,1}$-dimensional vector space. This operator arises from infinitesimal variation in the weight direction;
	from this construction, however, its intrinsic meaning, and in particular
	its relation to the weight $k$, seems less transparent.
	By contrast, Teitelbaum constructed a related but more concrete operator
	in his study of Teitelbaum’s $\mathcal{L}$-invariants. In particular, the
	underlying vector space in his construction is realized
	explicitly in terms of group cohomology with coefficients in the linear dual of
	$\operatorname{Sym}^{k-2}\QQ_p^2$; see \cite[\S~2.4]{Graf2019} and the
	references therein.
\end{remark}
\subsection{Slopes of the matrix $A_1$ }
\label{step2}
This step is the core of the paper. We now give a more detailed explanation of step (2) in diagram~\eqref{diag1}.

\subsection*{The $(k,w)$-newslopes}
In \S~\ref{dis of k-thresholds}, we use the local ghost series to study the following variation (to other weights $w$'s) of the slopes of $p$-newforms (of weight $k$).
\begin{definition}
	\label{def:introduction-(k,w_*)-newslope}
	Let $w\in\calW_{\{c+1\}}$ be a weight. The \emph{(global) $(k,w)$-newslopes} are defined to be the $\left(d_{k,0}+1\right)$-th to the $\left(d_{k,0}+d_{k,1}\right)$-th slopes of the Newton polygon $\NP(w)$. For $i\in\left\{1,\dots,d_{k,1}\right\}$, the \emph{$i$-th $(k,w)$-newslope} is the $\left(d_{k,0}+i\right)$-th slope of $\NP(w)$.  
\end{definition}
Our main result concerning these slopes is Theorem~\ref{propofslopes}, where we determine all the \((k,w)\)-newslopes for every \(w\) lying in a certain affinoid open disc \(D(k) \subseteq \mathcal{W}_{c+1}\), centered at \(k\) and determined by the local ghost series; see Definition~\ref{good region}. This result provides the source from which we may extract the slopes of \(A_1\).

In \S~\ref{subsection newslope comp}, we also give the following interpretation of these slopes for $w\in D(k)$. Let \(\operatorname{Spm}\mathcal{A}=\mathcal{D}\subseteq D(k)\) be an affinoid closed disc centered at \(k\). Then we construct a finite flat \(\mathcal{A}\)-module \(\mathcal{S}\) of rank \(d_{k,1}\), equipped with an \(\mathcal{A}\)-linear action of \(U_p\), such that, for any weight \(w \in \mathcal{D}\), the eigenvalues of the \(U_p\)-action on the specialization \(\mathcal{S}_w\) coincide with the \(U_p\)-eigenvalues appearing in \(S_w^\dagger(\Gamma_0(Np))_{\bar r}\) that correspond to the \((k,w)\)-newslopes. Consequently, the characteristic polynomial of the \(U_p^2\)-action on \(\mathcal{S}\),
\begin{equation}
	\label{4eq}
	\det(1-U_p^2(w)\cdot t)=1+\sum_{i=1}^{d_{k,1}}(-1)^ic_{i}^{\new}(w)t^i,
\end{equation}
is well defined for every \(w\in D(k)\).
Note that the open disc \(D(k)\) is strictly larger than the disc \(\mathcal{D}'\) appearing in \S~\ref{step1}. This larger domain is essential for the analytic trick introduced later; see Remark~\ref{rmk: analytic trick}.
\subsection*{From \((k,w)\)-newslopes to slopes of \(A_1\)}
In \S~\ref{subsection two step approach}, we present a method for extracting slopes of $A_1$ from $(k,w)$-newslopes for various weights $w\in D(k)$ (especially for those \emph{not} very close to $w_k$). This method relies on the linear equations and an analytic technique developed in \S~\ref{section : linear equations}. 

For \(n\) variables \(\underline{\xi}=(\xi_1,\dots,\xi_n)\), let
\(\mathscr{S}_i(\underline{\xi})\) denote the \(i\)-th elementary symmetric polynomial. For a diagonal matrix
$
A=\operatorname{diag}(a_1,\dots,a_n)$,
we write
$
\mathscr{S}_i(\underline{A})
:=
\mathscr{S}_i(a_1,\dots,a_n)$.
Then
$
\det(1-tA)
=
1+\sum_{i=1}^{n}(-1)^i\mathscr{S}_i(\underline{A})t^i$.
It follows that, in order to determine the slopes of \(A_1\), it suffices to compute the valuations
$
\left\{
v_p\bigl(\mathscr{S}_i(\underline{A_1})\bigr)
\right\}_{1\le i\le d_{k,1}}$.

For \(i\in\{1,\dots,d_{k,1}\}\), let
$
c_i^{\new}(w)
=
\sum_{\ell\ge 0} c_i^{\new,(\ell)}(w-w_k)^\ell$
be the Taylor expansion of the \(i\)-th coefficient in \eqref{4eq},
which converges for \(w\in D(k)\).
 If \(w\in\mathcal D'\), then, by the discussion in \S\ref{step1}, the operator \(U_p^2(w)\) is represented by the matrix
$
M(w)
=
p^{k-2}I_{d_{k,1}}
+
A_1(w-w_k)
+\cdots$.
Consequently,
$
c_i^{\new}(w)
=
\mathscr{S}_i\bigl(\underline{M(w)}\bigr)$.
Using this identity, in
\S\ref{subsection:construction of linear equations}, we express each \(\mathscr{S}_i(\underline{A_1})\) as an explicit \(\mathbb Z\)-linear combination of the coefficients
$
\{\,c_j^{\new,(i)} : 1\le j\le d_{k,1}\,\}$;
see Lemma~\ref{lem:ZZ-linaer comb lemma} and
Proposition~\ref{prop: ZZ-linear comb prop}.

Next, we apply an analytic trick, explained in
\S~\ref{subsection: analytic trick}, to estimate the valuations
\[
\left\{
v_p\bigl(c_i^{\new,(\ell)}\bigr):
i,\ell=1,\dots,d_{k,1}
\right\}.
\]
Let \(r_0\) be the valuation radius of \(D(k)\). For \(r>r_0\), define
the \(r\)-Gauss norm of \(c_i^{\new}\) by
$
\nu_r(c_i^{\new})
:=
\min_{v_p(w-w_k)=r}
\left\{
v_p(c_i^{\new}(w))
\right\}$.
Then the dual graph
$
\left\{
(r,\nu_r(c_i^{\new}))\in\mathbb{R}^2:
r>r_0
\right\}
$
determines the segments of the Newton polygon \(\operatorname{NP}(c_i^{\new})\)
with slopes \(<-r_0\); see Lemma~\ref{lem:dual}. In our setting, the
\((k,w)\)-newslopes for various weights \(w\in D(k)\) provide certain
lower bounds for this dual graph. From these bounds, we obtain the desired
estimates for the valuations \(v_p(c_i^{\new,(\ell)})\); see
Proposition~\ref{main estimate for scrM}.
Notably, the smaller the value of \(r\), the more segments of
\(\operatorname{NP}(c_i^{\new})\) we can recover. It is therefore crucial
to consider \((k,w)\)-newslopes for weights \(w\) that are not very close
to \(w_k\), that is, for weights lying in the larger disc \(D(k)\) rather
than in \(\mathcal{D}'\).

\begin{remark}
	\label{rmk: analytic trick}
	Roughly speaking, we use the valuations \(v_p(c_i^{\new}(w))\), for
	weights \(w\) not very close to \(w_k\), to estimate, and in certain
	cases even determine, the valuations \(v_p(c_i^{\new,(\ell)})\). The
	motivation for this method is that we do not know the explicit values
	of \(c_i^{\new}(w)\), but only their valuations.
	A priori, each coefficient \(c_i^{\new,(\ell)}\) is determined by the
	values of \(c_i^{\new}(w)\) for weights \(w\) sufficiently close to
	\(w_k\). However, for such weights, the valuation
	\(v_p(c_i^{\new}(w))\) is constant. The variations of
	\(v_p(c_i^{\new}(w))\) caused by the valuations
	\(v_p(c_i^{\new,(\ell)})\) appear only when \(w\) is not very close to
	\(w_k\). It is precisely these variations that allow us to estimate
	the valuations \(v_p(c_i^{\new,(\ell)})\).
	
	This approach is similar to the heuristic for
	\(\mathcal{L}\)-invariants as gradients described in
	\cite[\S~8.4]{BP dis}, an idea that dates back to a private
	communication between Buzzard and Pollack in 2007. It is also closely
	related to the method of using slopes of \(\mathcal{L}\)-invariants
	to bound the constant slope radius of Coleman families in \cite{const}.
\end{remark}

Combining these estimates (Proposition \ref{main estimate for scrM}) with the $\ZZ$-linear relations (Proposition \ref{Cor for Simply F_i(u)}) we deduce the following theorem.

\begin{theorem}[Corollary \ref{cor main1}]
	\label{thm step3}
Write \(k=k_0+n(p-1)\), where \(k_0\in\{0,\dots,p-2\}\). Recall that
\(C(k,a)\) and $m(\bar r)$ are constants defined in
Remark~\ref{rmk in the introduction of constans}. Let \(S(k,\bar r)\) be
the multiset of slopes of \(\underline{\Delta}_k\), with each slope
repeated \(m(\bar r)\) times. Arrange both the slopes in \(S(k,\bar r)\) and the slopes of \(A_1\) in
ascending order. 
 Then the multiset of slopes of \(A_1\) is equal to the
multiset
$
\{k-2-s : s\in S(k,\bar r)\,\}
$
except possibly for the last $\lceil m(\bar r)\cdot C(k,a)\rceil$ elements. Moreover, each exceptional
slope of \(A_1\) is greater than or equal to
$
k-4-\log_p\left\lfloor \frac{k-k_0}{p-1}\right\rfloor$.

\end{theorem}

Combining Theorems~\ref{thm step2} and~\ref{thm step3} with the equidistribution and integrality results for the slopes of \(\underline{\Delta}_k\) (see Proposition~\ref{nearStprop} and Corollary~\ref{cor:equidistribution of Delta slope}), we obtain the main results, Theorems~\ref{Main theorem on BP dis} and~\ref{cor: integrality inroduction}.

	\subsection{A byproduct: the global k-thresholds}
	\label{subsection in the introduction: global $k$-thresholds}
In this final subsection, we discuss a byproduct of our analysis of \((k,w)\)-newslopes. This is also related to upper bounds for constant slope families studied in \cite{const}; see Remark~\ref{rmk:Bergdall CS and my CS} for further discussion.
	
By the \(p\)-adic continuity of slopes with respect to weights
 (\cite[Theorem D]{Coleman} or Theorem~\ref{propofslopes}), the global
\((k,w)\)-newslopes remain all equal to \(\frac{k-2}{2}\) whenever \(w\) is
sufficiently close to \(w_k\). The following definition records the
thresholds for this constancy.
	\begin{definition}
		\label{global k-thresholds introduction}
		The \emph{(global) $k$-thresholds} are the $d_{k,1}$ valuations 
		$\CS_{1,\bar r}(k),\dots,\CS_{d_{k,1},\bar r}(k)$, where each valuation $\CS_{i,\bar r}(k)$ is given by
\begin{equation}
	\begin{aligned}
		\CS_{i,\bar r}(k) := \inf \Bigl\{ s : \ &\textrm{for any } w \in \calW_{\{c+1\}} \textrm{ such that } v_p(w-w_k) = s, \\
		&\textrm{the } i\textrm{-th } (k,w)\textrm{-newslope is } \frac{k-2}{2} \Bigr\}\geq0.
	\end{aligned}
\end{equation}
	\end{definition}
An immediate consequence of our main result on \((k,w)\)-newslopes for
\(w\in D(k)\), namely Theorem~\ref{propofslopes}, is a parallel result to
Theorem~\ref{thm step3} for (global) \(k\)-thresholds; see
Corollary~\ref{cor of thresholds}. This result shows that the slopes of
\(\underline{\Delta}_k\) compute almost all global \(k\)-thresholds.
Consequently, the integrality and equidistribution of global
\(k\)-thresholds follow directly from the corresponding results for the
slopes of \(\underline{\Delta}_k\); see
Corollary~\ref{cor:integrality of thresholds} and
Theorem~\ref{distribution thm for constant}.
Since the slopes of \(\underline{\Delta}_k\) compute almost all slopes of the
\(\mathcal{L}\)-invariants and almost all global \(k\)-thresholds, the
following theorem gives a comparison between these two multisets.
	\begin{theorem}
		\label{thm:global-thresholds-and-L-invariants}
		Recall the constant $C(k,a)$ and $m(\bar r)$ are defined in Remark~\ref{rmk in the introduction of constans}.
		Fix \(p\geq 11\). Assume that \(\bar r\) is locally reducible and very generic.
		Then, up to at most $m(\bar r)\cdot \lceil C(k,a)\rceil$ exceptions, the following two multisets of
		valuations coincide:
		\begin{enumerate}
			\item the normalized slopes of \(\mathcal{L}\)-invariants,
			\[
			\left\{
			-v_p(\mathcal{L}_f)+1 :
			f \text{ is a normalized \(p\)-new eigenform in }
			S_k(\Gamma_0(Np))_{\bar r}
			\right\};
			\]
			\item the global \(k\)-thresholds,
			\[
			\left\{
			\CS_{i,\bar r}(k) : i=1,\dots,d_{k,1}
			\right\}.
			\]
		\end{enumerate}
	\end{theorem}
	
\begin{remark}
Let \(f\in S_k(\Gamma_0(Np))\) be a $p$-newform.	A related result is proved in \cite[Theorem~4.3]{const}, where the
author shows that $-v_p(\calL_f)+v_p(2)$ is a lower bound for the
constant slope radius $\CS^k(f)$. He also gives examples in which this
bound is sharp; see \cite[\S5]{const}. 	

Our global $k$-thresholds and Bergdall's constant slope radii appear to
measure related phenomena. However, the two notions are defined in rather
different ways: the former is expressed in terms of the variation of the segments of
\(\NP(C_{\bar r}(w,t))\) corresponding to $p$-newforms, whereas the latter
is attached to individual $p$-newforms $f\in S_k(\Gamma_0(Np))$; see Remark~\ref{rmk:Bergdall CS and my CS} for more discussions.
For this reason,
the precise relation between Theorem~\ref{thm:global-thresholds-and-L-invariants}
and \cite[Theorem~4.3]{const} remains unclear to us.
\end{remark}
	
	\subsection*{Notations and normalizations}
	Throughout this paper, we fix a prime $p\geq5$.
	Let $E$ be a finite extension of $\QQ_p(\sqrt p)$, as the coefficient field. Let $\calO$, $\FF$, and $\varpi$ be its ring of integers, residue field, and a uniformizer, respectively. Let $\I_{\QQ_p}\subset \Gal(\overline{\QQ}_p/\QQ_p)$ be the inertia subgroup, and  $\omega_1: \I_{\QQ_p} \twoheadrightarrow\Gal(\QQ_p(\mu_p)/\QQ_p) \cong \FF_p^\times$ \emph{the first fundamental character}. For $R$ a $p$-adic ring and $\alpha \in R^\times$, let $\unr(\alpha): \Gal(\bar\QQ_p/\QQ_p) \to R^\times $ be the unramified representation that sends the geometric Frobenius to $\alpha$.

	We will encounter both the $p$-adic logarithm $\log(x) = (x-1) - \frac{(x-1)^2}2 + \cdots $ for $x$ a $p$-adic number or a formal element, and the natural logarithm $\log_p(-)=\frac{\ln(-)}{\ln(p)}$ in the real analysis.
	
	For each $m \in \ZZ$, we write $\{m\}$ for the unique integer
	satisfying the conditions 
	$0\leq \{m\}\leq p-2$ and $m\equiv \{m\}
	\pmod{p-1}.$
	
	For a power series $f(X)=1+\sum_{n\geq1}a_nX^n\in\CC_p[\![X]\!]$, we use $\NP(f)$ to denote its \emph{Newton polygon}, namely the lower convex hull of the points $\left(n,v_p(\a_n)\right)$ for all $n$; the slopes of the segments of $\NP(f)$ are referred to as slopes of $f(X)$.
	
For a formal $\calO$-scheme $\Spf(R)$ formally of finite type, let $\Spf(R)^\rig$ be its associated rigid analytic space over $E$. 
	
	\subsection*{Acknowledgements}
	This work grew out of my thesis, and I am deeply grateful to my advisor, Liang Xiao, for introducing me to the slope theory and for his encouragement and support. I also owe a special debt of gratitude to Robert Pollack and John Bergdall, who generously shared their ideas and insights, as well as their hospitality during my visits to Boston University and the University of Arkansas. I also thank Bin Zhao, Peter Gräf, and Rufei Ren for helpful discussions. I thank the anonymous referee for numerous detailed comments and suggestions, which substantially improved the exposition and organization of the paper.

	%
	%
	%
	%
	%
	%
	%
	\section{Recollection of the local ghost theorem}
	\label{Section2}
	In this section, we review the main results of \cite{Lghost} and \cite{Lghost2}.
	We begin by recalling the definition of the \emph{local ghost series} associated
	with a generic \emph{Serre weight}
	\(\sigma=\Sym^a\FF^2\otimes\det^b\), where
	\(a\in\{1,\dots,p-4\}\) and \(b\in\{0,\dots,p-2\}\). After this, we state the local
	ghost theorem, Theorem~\ref{local ghost theorem}, for irreducible Galois
	representations
	$
	\bar r:\Gal(\bar\QQ/\QQ)\rightarrow \GL_2(\FF)
	$
	satisfying certain local conditions. In
	\S~\ref{subsection: vetices} and \S~\ref{subsection: Delta slopes}, we also
	recall several results from \cite{Lghost} concerning the vertices and slopes of
	the \emph{Newton polygon} of the local ghost series. In particular, the equidistribution of the \emph{\(k\)-derivative slopes}
	(Definition~\ref{defhat}) is an immediate consequence of the estimates in
	\S~\ref{subsection: Delta slopes}; see
	Corollary~\ref{cor:equidistribution of Delta slope}. This equidistribution
	result is essential for proving the equidistribution of the slopes of
	\(\calL\)-invariants in the final section.
	
	\subsection{Weight discs}
	\label{sunsection weight space}
	Identify $\calO[\![1+p\ZZ_p]\!]\cong\calO[\![w]\!]$ via sending $[\a]$ for $\a\in1+p\ZZ_p$ to $(1+w)^{\log(\a)/p}$, where $\log(-)$ is the formal $p$-adic logarithm. In particular, for each $k\in\ZZ$, we have
	\begin{equation}
		w_{k}:=\exp(p(k-2))-1,
	\end{equation}
and
	\begin{equation}
		\label{valuation}
		v_{p}\big(w_{k}-w_{k'}\big)=v_{p}\big(\exp(p(k'-2))\cdot(\exp(p(k-k'))-1)\big)=1+v_{p}\big(k-k'\big).
	\end{equation}
	Let $\Delta$ be the torsion subgroup of $\ZZ_p^\times$. Recall for each character $\varepsilon: \Delta\rightarrow \calO^{\times}$, we have the corresponding \emph{weight disc} $\calW^{(\e)}=(\Spf\calO[\![w]\!])^{\rig}$, parameterizing continuous characters $\chi$ of $\ZZ_p^{\times}$ such that $\chi|_{\Delta}=\e$. If we write $\e=\omega^{s_{\e}}$ for some $s_{\e}\in\{0,\dots,p-2\}$, where $\omega$ is the Teichm\"{u}ller character, then $\calW^{(\e)}$ coincides with $\calW_{s_{\e}}$ defined in \S~\ref{overview of the proof}.
	
	\subsection{Serre weights and local Galois representations}
For $a \in \{0, \dots, p-1\}$ and $b \in \{0, \dots, p-2\}$, we use $\sigma_{a,b}$ to denote the \emph{right} $\FF$-representation $\Sym^a\FF^{\oplus 2} \otimes \det^b$ of $\GL_2(\FF_p)$; these exhaust all irreducible right $\FF$-representations of $\GL_2(\FF_p)$. We call them the \emph{Serre weights}. A Serre weight $\sigma_{a,b}$ is called \emph{generic} (resp. \emph{very generic}) if $1\leq a\leq p-4$ (resp. $2\leq a\leq p-5$).
	
	\begin{definition}
		Let $\bar r: \Gal(\bar\QQ/\QQ) \to \GL_2(\FF)$ be an irreducible representation. We say $\bar r$ is \emph{locally reducible and very generic}, if
		\begin{equation}
			\bar\rho=\bar r|_{\Gal(\bar\QQ_p/\QQ_p)}\simeq \begin{pmatrix}\unr(\bar\a)\omega_{1}^{a+b+1}&*\\0&\unr(\bar\b)\omega_1^b\end{pmatrix},
		\end{equation}
		where $a\in\{2,\dots,p-5\}$, $b\in\{0,\dots,p-2\}$, and $\bar\a,\bar\b\in\FF^{\times}$.
	\end{definition}
	
	\subsection{Local ghost series for generic Serre weights} 
	The key idea of defining \emph{ghost series} arises from the dimension formulas related to \emph{oldforms} and \emph{newforms}, which satisfy specific asymptotic properties. Here we only present these dimension formulas (Definition \ref{dimformula}) from \cite[\S~4]{Lghost} to define the local ghost series.
	Fix a character $\e$ of $\Delta$. Let $\sigma=\sigma_{a,0}$ be a \emph{generic} Serre weight.
	 	\begin{notation}
	 	\label{notation1}
	 	For $\varepsilon: \Delta\rightarrow \calO^{\times}$ and $\sigma=\sigma_{a,0}$, we associate the following numbers.
	 	\begin{enumerate}
	 		\item Recall $\e=\omega^{s_{\e}}$, for some $s_{\e}\in\{0,\dots,p-2\}$, where $\omega$ is the Teichm\"{u}ller character.
	 		\item 
	 		Set $k_{\e,a}:=2+\{a+2s_{\e}\}\in\{2,\dots,p\}$. 
	 		\item 
	 			Set
	 		$\d_{\e,a}:=\left\lfloor\frac{s_{\e}+\{a+s_{\e}\}}{p-1}\right\rfloor$ and we introduce two integers $t_{1}^{(\e,a)},t_{2}^{(\e,a)}\in\ZZ$,
	 		\begin{itemize}
	 			\item [$\bullet$] if $a+s_{\e}<p-1$, then $t_{1}^{(\e,a)}=s_{\e}+\d_{\e,a}$ and $t_{2}^{(\e,a)}=a+s_{\e}+\d_{\e,a}+2$;
	 			\item[$\bullet$] if  $a+s_{\e}\geq p-1$, then $t_{1}^{(\e,a)}=\{a+s_{\e}\}+\d_{\e,a}+1$ and $t_{2}^{(\e,a)}=a+s_{\e}+\d_{\e,a}+1$.
	 		\end{itemize}
	 	\end{enumerate}
	 \end{notation}

	\begin{definition}
		\label{dimformula}
		For an integer $k \geq 2$ satisfying $k \equiv k_{\e,a} \mod (p-1)$, we rewrite it as $k=k_{\e,a}+k_{\bullet}(p-1)$. Consider the following data regarding dimensions.
		\begin{enumerate}
			\item 
			Define $d_{k}^{\Iw}(\e,a):=2k_{\bullet}+2-2\d_{\e,a}$, then we have
			$
			d_{k}^{\Iw}(\e,a)=\frac{2k}{p-1}+O(1).		   
			$
			\item
			Define $d_{k}^{\ur}(\e,a):=2\left\lfloor\frac{k_{\bullet}-t_{1}^{(\e,a)}}{p+1}\right\rfloor+1+\left\{\begin{array}{ll}1&\mbox{if $k_{\bullet}-(p+1)\left\lfloor\frac{k_{\bullet}-t_{1}^{(\e,a)}}{p+1}\right\rfloor\geq t_{2}^{(\e,a)}$,}\\0&\mbox{otherwise.}\end{array}\right.$ In particular, the integer $d_{k}^{\ur}(\e,a)$ is non-decreasing in the class $k\equiv k_{\e,a}\pmod{p-1}$. Moreover, we have
			$
			d_{k+p^{2}-1}^{\ur}(\e,a)-d_{k}^{\ur}(\e,a)=2$, and 
			$d_{k}^{\ur}(\e,a)=\frac{2k}{p^{2}-1}+O(1).    
			$
			\item
			Define $d_{k}^{\new}(\e,a):=d_{k}^{\Iw}(\e,a)-2d_{k}^{\ur}(\e,a)$, then we have $d_{k}^{\new}(\e,a)=\frac{2k}{p+1}+O(1)$.	
		\end{enumerate}
	\end{definition}

	\begin{definition}
		\label{def: local ghost series}
		The \emph{local ghost series} attached to the generic Serre weight $\sigma=\sigma_{a,0}$ over the weight disc $\calW^{(\e)}$ is the formal power series
		\[
		G(w,t)=G_{\sigma}^{(\e)}(w,t)=1+\sum_{n=1}^{\infty}g_{n}^{(\e,a)}(w)t^{n}\in\calO[w][\![t]\!],
		\]
		where the \emph{ghost polynomial} $g_{n}^{(\e,a)}(w)$ is given by
		\[
		g_{n}^{(\e,a)}(w)=\prod_{k\geq2,\ k\equiv k_{\e,a}\pmod {p-1}}(w-w_{k})^{m_{n}^{(\e,a)}(k)}\in\ZZ_p[w],
		\]
		with the \emph{ghost multiplicity} $m_{n}^{(\e,a)}(k)$ given by
		\[
		\label{ghostmult}
		m_{n}^{(\e,a)}(k)=\left\{\begin{array}{ll}\min\{n-d_{k}^{\ur}(\e,a),d_{k}^{\Iw}(\e,a)-d_{k}^{\ur}(\e,a)-n\}&\mbox{if $n\in\left(d_{k}^{\ur}(\e,a),d_{k}^{\Iw}(\e,a)-d_{k}^{\ur}(\e,a)\right)$},\\
			0&\mbox{otherwise}.\end{array}\right.\ \
		\]
		For each $k$ for which $m_{n}^{(\e,a)}(k)\neq0$, we refer to $w_{k}$ as a \emph{ghost zero} of the ghost polynomial $g_{n}^{(\e,a)}(w)$. 
	\end{definition}
	\begin{example}
			\label{Appendix ghost series} 
			Consider $p=7$, $a=2$, and the weight disc $\e=\omega$. In this case, $s_{\e}=1$, $k_{\e,a}=6$. We list the first eight ghost polynomials below, omitting the subscripts $(\epsilon, a)$ for simplicity. 
			\begin{equation}
		\begin{aligned}
			g_{1}(w)=&w-w_{6}.\\
			g_{2}(w)=&(w-w_{12})(w-w_{18})(w-w_{24})(w-w_{30}).\\
			g_{3}(w)=&(w-w_{18})^{2}(w-w_{24})^{2}(w-w_{30})^{2}(w-w_{36})\cdots(w-w_{54}).\\
			g_{4}(w)=&(w-w_{18})(w-w_{24})^{3}(w-w_{30})^{3}(w-w_{36})^{2}\cdots(w-w_{54})^{2}(w-w_{60})\cdots(w-w_{78}).\\
			g_{5}(w)=&(w-w_{24})^{2}(w-w_{30})^{4}(w-w_{36})^{3}\cdots(w-w_{54})^{3}(w-w_{60})^{2}\cdots(w-w_{78})^{2}\\&(w-w_{84})\cdots(w-w_{102}).\\
			g_{6}(w)=&(w-w_{24})(w-w_{30})^{3}(w-w_{36})^{4}\cdots(w-w_{54})^{4}(w-w_{60})^{3}\cdots(w-w_{78})^{3}\\&(w-w_{84})^{2}\cdots(w-w_{102})^{2}(w-w_{108})\cdots(w-w_{126}).\\
			g_{7}(w)=&(w-w_{30})^{2}(w-w_{36})^{3}(w-w_{42})^{5}(w-w_{48})^{5}(w-w_{54})^{5}(w-w_{60})^{4}\cdots(w-w_{78})^{4}\\&(w-w_{84})^{3}\cdots(w-w_{102})^{3}(w-w_{108})^{2}\cdots(w-w_{126})^{2}(w-w_{132})\cdots(w-w_{150}).\\
			g_{8}(w)=&(w-w_{30})(w-w_{36})^{2}(w-w_{42})^{4}(w-w_{48})^{6}(w-w_{54})^{6}(w-w_{60})^{5}\cdots(w-w_{78})^{5}\\&(w-w_{84})^{4}\cdots(w-w_{102})^{4}(w-w_{108})^{3}\cdots(w-w_{126})^{3}(w-w_{132})^{2}\cdots(w-w_{150})^{2}\\&(w-w_{156})\cdots(w-w_{174}).
		\end{aligned}
		\end{equation}
	\end{example}

	\subsection{Local ghost theorem}
	\label{subsection local ghost theorem}
Let
$
\bar r \colon \Gal(\overline{\QQ}/\QQ) \to \GL_2(\FF)
$
be an irreducible representation that is locally reducible and very generic.
Assume that
\begin{equation}
	\det(\bar r|_{\I_{\QQ_p}}) = \omega_1^c
	\quad \text{for some } c \in \{0,\dots,p-2\}.
\end{equation}
	Recall $S_k(\Gamma_0(Np))_{\bar r} \subseteq S^\dagger_k(\Gamma_0(Np))_{\bar r}$ are the spaces of classical and overconvergent modular forms of weight $k$, level $\Gamma_0(Np)$, localized at the Hecke maximal ideal corresponding to $\bar r$, respectively. The space $S^\dagger_k(\Gamma_0(Np))_{\bar r}$ is non-zero only if $k\equiv c+1 \pmod {p-1}$. Following the discussions of \cite[\S~1.2]{Lghost2}, there exists a formal power series
	$C_{\bar r}(w,t) \in \calO[\![w,t]\!]$ such that, for every
	$k \equiv c+1 \pmod {p-1}$, the specialization
	$C_{\bar r}(w_k,t)\in\calO[\![t]\!]$ is precisely the
	\emph{characteristic power series} of the $U_p$-action on
	$S^\dagger_k(\Gamma_0(Np))_{\bar r}$, namely,
	\[
	C_{\bar r}(w_k,t)
	=
	\det\bigl(\I_\infty - U_p t;\,
	S^\dagger_k(\Gamma_0(Np))_{\bar r}\bigr).
	\]
We call $C_{\bar r}(w,t)$ the \emph{characteristic power series} of the $U_p$ operator. To state the main result of \cite{Lghost2}, we recall the integer
\(m(\bar r)\) introduced in \cite[\S~1.2]{Lghost2}, which is characterized
by the property that
\[
\dim S_k(\Gamma_0(Np))_{\bar r} - \frac{2k}{p-1}m(\bar r)
\]
is bounded as \(k \to \infty\). We refer to \(m(\bar r)\) as the
\emph{global multiplicity} of \(\bar r\).
\begin{remark}
	\label{rmk global multiplicity}
In \cite[\S~8]{Lghost2}, the authors give a conceptual interpretation of the
integer \(m(\bar r)\). Namely, \(m(\bar r)\) is the rank of the
\emph{\(\calO[\![\GL_2(\ZZ_p)]\!]\)-projective arithmetic module of type
	\(\bar r|_{\Gal(\overline{\QQ}_p/\QQ_p)}\)} determined by the spaces
\(S_k(\Gamma_0(Np))\); see \cite[Definition~8.2, Example~8.4]{Lghost2}.
In particular, \(m(\bar r)\) depends on the tame level \(N\).

For explicit values, \cite[\S~6]{BP2} expresses \(m(\bar r)\) in terms of the
dimensions of the spaces \(S_k(\Gamma_0(Np))_{\bar r}\) for certain small
weights \(k\). More precisely, our \(m(\bar r)\) is
equal to \(\frac{1}{2}Q_{d_t}\) in \cite[Proposition~6.9]{BP2}, with the formula
for \(Q_{d_t}\) given in \cite[Corollary~6.17]{BP2}. Note that \(Q_{d_t}\) is even when \(\bar r\) is locally generic; see
\cite[Equation~(18) in \S~6.3 and Corollary~6.17]{BP2}.
\end{remark}

	\begin{theorem}	[local ghost theorem]
		\label{local ghost theorem}
		Assume $p\geq 11$. 	Let $\bar r: \Gal(\bar\QQ/\QQ) \to \GL_2(\FF)$ be an irreducible representation, locally reducible and very generic, such that
		\[
			\bar\rho:=\bar r|_{\Gal(\bar\QQ_p/\QQ_p)}\simeq \begin{pmatrix}\unr(\bar\a)\omega_{1}^{a+b+1}&*\\0&\unr(\bar\b)\omega_1^b\end{pmatrix},
		\]
	for some $\bar\a,\bar\b\in\FF^{\times}$. Denote its global multiplicity by $m(\bar r)$. Consider the following data
		\[
		\sigma=\sigma_{a,0},\ \e=\omega^{b}.
		\]
		Then for any weight $w_*\in\gothm_{\CC_{p}}$, the Newton polygon $\NP\big(C_{\bar r}(w_*,-)\big)$ is the same as the Newton polygon $\NP\big(G^{(\varepsilon)}_{\sigma}(w_*,-)\big)$, stretched in both $x$- and $y$-directions by $m(\bar r)$, except possibly for the slope zero parts.
	\end{theorem}

	\subsection{Settings}\label{subsection of Settings}
	In this subsection, we fix the setup we will work with in the remaining part of this paper.
	\begin{enumerate}
		\item
		Fix an irreducible global Galois representation $\bar r$  with global multiplicity $m(\bar r)$, which is locally reducible and very generic. Denote its local restriction by $\bar\rho:=\bar r|_{\Gal(\bar\QQ_p/\QQ_p)}\simeq \begin{pmatrix}\unr(\bar\a)\omega_{1}^{a+b+1}&*\\0&\unr(\bar\b)\omega_1^b\end{pmatrix}$. Fix the character $\e=\omega^b$ and the very generic Serre weight $\sigma=\sigma_{a,0}$ as in Theorem \ref{local ghost theorem}.
		\item
		For this pair $(\e,\sigma_{a,0})$ in (1), we have the corresponding integer $k_0:=k_{\e,a}=2+\{a+2s_\e\}=2+\{a+2b\}$. We also simplify notations by omitting the subscripts $(\e,a)$ in
		Notation \ref{notation1}, Definition \ref{dimformula}, and Definition \ref{def: local ghost series}. In particular, the local ghost series is written as $G(w,t)=1+\sum_{n=1}^{\infty}g_{n}(w)t^{n}\in\calO[w][\![t]\!]$.
		\item
		Denote by $\calC_{\bar r}$ the $\bar{r}$-component of the \emph{Coleman--Mazur eigencurve} of tame level $N$, over the weight space $\calW$. By \emph{$\bar r$-newforms} of weight $k$, we mean the $p$-newforms (normalized eigenforms which are new at $p$) in $S_{k}(\Gamma_0(Np))_{\bar r}$.
	\end{enumerate}
Throughout, whenever we discuss asymptotic properties or distributions as
\(k\to\infty\), we mean that \(k\to\infty\) with \(k\) restricted to the fixed
congruence class \(k\equiv k_0 \pmod{p-1}\).
	
		\subsection{Characterization for the vertices of $\NP(G(w,t))$}
		\label{subsection: vetices}
		
		In this subsection, we review a characterization (Theorem \ref{iff}) for the vertices of the \emph{ghost polygon} (i.e. the \emph{Newton polygon} of the $G(w,t)$).
	To begin, we revisit a family of polygons indexed by $k\geq2$, which was introduced in \cite{Lghost} to control the behavior of the ghost polygon. Recall the definition of the ghost series $G(w,t)=1+\sum_{n\geq1}g_n(w)t^n$, where $g_n(w)=\prod_{k\geq2,\ k\equiv k_0\pmod{p-1}}(w-w_{k})^{m_{n}(k)}$.
	\begin{definition}
		\label{defofDelta}
		
		Let $g_{n,\hat{k}}(w):=g_{n}(w)/(w-w_{k})^{m_{n}(k)}$. Define $\Delta'_{k,l}$ to be the following valuation
		\begin{equation}
			\label{defhat}
			\Delta'_{k,l}:=v_{p}\big(g_{\frac{1}{2}d_{k}^{\Iw}+l,\hat{k}}(w_{k})\big)-\frac{k-2}{2}l,\ \ \textrm{for all}\ l=-\frac{1}{2}d_{k}^{\new},\dots,\frac{1}{2}d_{k}^{\new}.
		\end{equation}
		 The \emph{$k$-derivative polygon} is defined to be the lower convex hull of the points $$\left\{(l,\Delta'_{k,l}):l=-\frac{1}{2}d_{k}^{\new},\dots,\frac{1}{2}d_{k}^{\new}\right\},$$ 
		 and we denote this polygon by $\underline{\Delta}_{k}$. The slopes of $\underline{\Delta}_k$ are called \emph{$k$-derivative slopes}.
		 Moreover, we denote by $(l,\Delta_{k,l})$ the corresponding points on this convex hull $\underline{\Delta}_{k}$.		
By the symmetry of this polygon (Proposition \ref{gduality}), we only need to consider $\underline{\Delta}_{k}$, for $l=0,1,\ldots,\frac{1}{2}d_{k}^{\new}$. This half is denoted by $\underline{\Delta}^{+}_{k}$.
	
	\end{definition}
	\begin{proposition}[Ghost duality]
		\label{gduality}
		For $l=1,\dots,\frac{1}{2}d_{k}^{\new}$, we have
		\begin{equation}
			v_{p}\big(g_{\frac{1}{2}d_{k}^{\Iw}+l,\hat{k}}(w_{k})\big)-v_{p}\big(g_{\frac{1}{2}d_{k}^{\Iw}-l,\hat{k}}(w_{k})\big)=(k-2)\cdot l.
		\end{equation}
		In other words, we have
		$
		\Delta'_{k,-l}=\Delta'_{k,l}$, and therefore $\Delta_{k,-l}=\Delta_{k,l}$.
		Moreover, $\left(d_k^{\ur},v_p(g_{d_k^{\ur}}(w_k))\right)$ and  $\left(d_k^{\Iw}-d_k^{\ur},v_p(g_{d_k^{\Iw}-d_k^{\ur}}(w_k))\right)$ are vertices of $\NP(G(w_k,-))$, and thus the $(d_{k}^{\ur}+1)$-th to the $(d_{k}^{\Iw}-d_{k}^{\ur})$-th slopes of $\NP(G(w_{k},-))$ are all equal to $\frac{k-2}{2}$. 
	\end{proposition}
	\begin{proof}
		See \cite[Proposition 4.18 (4) and Notation 5.1]{Lghost}.
	\end{proof}
	Let $w_{*}\in\gothm_{\CC_{p}}$ be a weight. To determine the vertices of $\NP(G(w_*,-))$, we recall the \emph{near-Steinberg range} introduced in \cite{Lghost}.
	\begin{definition}
		\label{nS def}
		Let $n$ be a positive integer. We say the pair $(w_{*},n)$ is \emph{near-Steinberg} for $k$, if the following two conditions hold.
		\begin{equation}
			n\in(d_{k}^{\ur},d_{k}^{\Iw}-d_{k}^{\ur}),\ \ \textrm{and}\ \ v_{p}\big(w_{*}-w_{k}\big)\geq \Delta_{k,|n-\frac{1}{2}d_{k}^{\Iw}|+1}-\Delta_{k,|n-\frac{1}{2}d_{k}^{\Iw}|}.
		\end{equation}
	\end{definition}
	
	\begin{theorem}
		\label{iff}
		For an integer $n\geq1$, $(n,v_{p}(g_{n}(w_{*})))$ is a vertex of $\NP(G(w_{*},-))$ if and only if the pair $(w_{*},n)$ is not  near-Steinberg for any $k$ (in the class $k\equiv k_0\pmod{p-1}$). 
	\end{theorem}
	\begin{proof}
		See \cite[Theorem 5.19]{Lghost}.
	\end{proof}
	The following proposition establishes a partial correspondence between the vertices of $\underline{\Delta}_{k}$ and the vertices of $\NP(G(w_{k},-))$. This will be used in subsequent discussions, particularly in the proofs of the integrality results in Corollary \ref{cor:integrality of thresholds} and Corollary \ref{cor: integrality}.
	\begin{proposition}
		\label{nearStprop}
		Fix $k\equiv k_0\pmod{p-1}$. The following statements are equivalent for any $l\in\{0,\dots,\frac{1}{2}d_{k}^{\new}-1\}$.
		\begin{enumerate}
			\item 
			The point $(l,\Delta'_{k,l})$ is not a vertex of $\underline{\Delta}_{k}$,
			\item 
			$(w_k,\frac{1}{2}d_{k}^{\Iw}-l)$ is near-Steinberg for some $k'<k$, and
			\item 
			$(w_k,\frac{1}{2}d_{k}^{\Iw}+l)$ is near-Steinberg for some $k''>k$.
		\end{enumerate}
		Moreover, the slopes of $\underline{\Delta}_{k}$ with multiplicity one belong to $\ZZ$. Other slopes all have even multiplicity and the slopes belong to $\frac{a}{2}+\ZZ=\frac{k}{2}+\ZZ$.
	\end{proposition}
	\begin{proof}
		See \cite[Proposition 5.26]{Lghost}. 
	\end{proof}

	\subsection{Estimates and equidistribution of $k$-derivative slopes}
	\label{subsection: Delta slopes}
In this subsection, we prove a technical lemma giving estimates for
\(\Delta'_{k,l}\) and \(\Delta_{k,l}\). These estimates immediately imply the
equidistribution of the \(k\)-derivative slopes; see
Corollary~\ref{cor:equidistribution of Delta slope}.
	\begin{lemma}
		\label{lemmaDelta}
		Assume that $p\geq5$. Write $k=k_0+k_{\bullet}(p-1)$. Fix some $l\in\{1,2,\dots,\frac{1}{2}d_{k}^{\new}\}$.
		\begin{enumerate}
			\item 
		We have
			\begin{equation}
				\label{1233}
				\frac{\min\{a+2,p-1-a\}}{2}+\frac{p-1}{2}(l-1)\leq\Delta'_{k,l}-\Delta'_{k,l-1}\leq\frac{3}{2}+\frac{p-1}{2}l+\b_{k,l},
			\end{equation}where $\b_{k,l}=O(\log_pk)$ as $k$ tends to infinity.
			\item
			Assume moreover that $p\geq 7$, then we have 
			\begin{equation}
				\label{12333}
				\Delta'_{k,l}-\Delta_{k,l}\leq 3\big(\log_p l\big)^{2}.		
			\end{equation}
		\item Assume that \(p\geq 7\), and write
		$
		s_l=\frac{2(p+1)}{p-1}\left(\Delta_{k,l}-\Delta_{k,l-1}\right)$.
		Then as \(k\to\infty\), we have
		\[
		\left|
		\frac{s_l}{k}
		-
		\frac{2l}{d_k^{\new}}
		\right|=
		\frac{O\left((\log_p k)^2\right)}{k}.
		\]
		In particular, this difference tends to zero uniformly in \(l\) as \(k\to\infty\).
		\end{enumerate}
	\end{lemma} 
	
	\begin{proof}
	We first prove (1) and (2). The estimates \eqref{1233} and \eqref{12333}
	are contained in \cite[Lemmas 5.2, 5.5, and 5.8]{Lghost}. It remains to prove
	the estimate for \(\beta_{k,l}\). By \cite[Lemma 5.5]{Lghost2}, we have
	$
	\beta_{k,l}=\beta+\lfloor \log_p l\rfloor$,
	where \(\beta\) is the largest \(p\)-adic valuation of an integer in the closed
	interval
	$
	\left[
	\eta-\frac{p+1}{2}(l-1),\,
	\eta+\theta+\frac{p+1}{2}(l-1)
	\right]$.
	Here \(\eta=O(k)\) and \(\eta+\theta=O(k)\) as  \(k\to\infty\), using the notation from the proof
	of \cite[Lemma 5.2]{Lghost2}. Since the endpoints of this interval are
	\(O(k)\), we obtain $\b=O(\log_pk)$  as \(k\to\infty\). Moreover, since
	$
	l\leq \frac{1}{2}d_k^{\new}=O(k)
	$
	as \(k\to\infty\) by Definition~\ref{dimformula}, we have
	$
	\beta_{k,l}=O(\log_p k)$ as  \(k\to\infty\).
	
		We now prove part (3). First note that
		\begin{equation}
			\label{2eq}
			\begin{aligned}
				\left|
				\frac{s_l}{k}
				-
				\frac{2l}{d_k^{\new}}
				\right|
				&\leq
				\left|
				\frac{s_l}{k}
				-
				\frac{(p+1)l}{k}
				\right|
				+
				\left|
				\frac{(p+1)l}{k}
				-
				\frac{2l}{d_k^{\new}}
				\right|  \\
				&=
				\left|
				\frac{s_l}{k}
				-
				\frac{(p+1)l}{k}
				\right|
				+
				\frac{l(p+1)}{k d_k^{\new}}
				\left|
				d_k^{\new}-\frac{2k}{p+1}
				\right|.
			\end{aligned}
		\end{equation}
		Since \(l<d_k^{\new}\) and $
		d_k^{\new}=\frac{2k}{p+1}+O(1)$ (see Definition~\ref{dimformula}), the second term in the last line of
		\eqref{2eq} is \(O(1/k)\). Therefore, it remains to show that
		$
		\left|
		\frac{s_l}{k}
		-
		\frac{(p+1)l}{k}
		\right|$ is
		$\frac{O\left((\log_p k)^2\right)}{k}$
		as \(k\to\infty\).
		We use the inequalities
		$
		\Delta_{k,l}-\Delta_{k,l-1}
		\leq
		\Delta'_{k,l}-\Delta_{k,l-1}
		=
		\left(\Delta'_{k,l}-\Delta'_{k,l-1}\right)
		+
		\left(\Delta'_{k,l-1}-\Delta_{k,l-1}\right)
		$
		and
		$
		\Delta_{k,l}-\Delta_{k,l-1}
		\geq
		\Delta_{k,l}-\Delta'_{k,l-1}
		=
		\left(\Delta'_{k,l}-\Delta'_{k,l-1}\right)
		+
		\left(\Delta_{k,l}-\Delta'_{k,l}\right).
		$
Then, by estimates in (1) and (2), we deduce that
		\[
		\left|
		\Delta_{k,l}-\Delta_{k,l-1}
		-
		\frac{p-1}{2}l
		\right|
		=
		O\left((\log_p k)^2\right).
		\]
		Hence
		$
		\left|
		\frac{s_l}{k}
		-
		\frac{(p+1)l}{k}
		\right|
		=
		\frac{O\left((\log_p k)^2\right)}{k}
		$
		as \(k\to\infty\).
		This completes the proof of (3).
	\end{proof}
	\begin{remark}
		\label{rmk about $l<2p$}
		For $l<2p$, more accurate estimates for $k$-derivative slopes can be found in \cite[Lemma 5.8]{Lghost}, namely, when $l<2p$ but $l\neq p$, $\Delta'_{k,l}=\Delta_{k,l}$; and when $l=p$, $ \Delta'_{k,l}-\Delta_{k,l}\leq1$. In particular, by \eqref{1233}, the first  \(k\)-derivative slope is at least \(	\frac{\min\{a+2,p-1-a\}}{2}\geq2\). These sharper estimates are useful for computing the $k$-derivative slopes for  small weights; see Example \ref{example1}.
	\end{remark}
	
\begin{corollary}
	\label{cor:equidistribution of Delta slope}
	Assume $p\geq7$. Let $\mu_{k}$ be the uniform probability measure of the following multiset
	\begin{equation}
		\begin{aligned}
		\left\{z_{k,i}=\frac{2(p+1)}{(p-1)k}\left(\Delta_{k,i}-\Delta_{k,i-1}\right):i=1,\dots, \frac{1}{2}d_k^{\new}\right\}\subseteq(-\infty,+\infty).
		\end{aligned}
	\end{equation}
Then $\mu_{k}$ weakly converges to the uniform probability measure on the interval $[0,1]$ as $k\to\infty$.
\end{corollary}	
\begin{proof}
For simplicity, write \(n_k=\frac{1}{2}d_k^{\new}\). For \(x\in\RR\), let
\(\delta_x\) denote the Dirac measure at \(x\). Consider
\[
\mu_k=\frac{1}{n_k}\sum_{i=1}^{n_k}\delta_{z_{k,i}},
\qquad
\tilde{\mu}_k=\frac{1}{n_k}\sum_{i=1}^{n_k}\delta_{i/n_k}.
\]
Note that \(\tilde{\mu}_k\) weakly converges to the uniform probability measure on
\([0,1]\) as \(k\to\infty\). Therefore, it suffices to show that \(\mu_k\) and
\(\tilde{\mu}_k\) have the same weak limit as \(k\to\infty\). This follows
directly from the fact that
$
\left|z_{k,i}-\frac{i}{n_k}\right|
$
tends to \(0\) uniformly in \(i\) as \(k\to\infty\), which is precisely
Lemma~\ref{lemmaDelta}(3). This completes the proof.
\end{proof}

	\section{Variation of $(k,w_*)$-newslopes and $k$-thresholds}\label{dis of k-thresholds}
	In this section, we use the local ghost series $G(w,t)$ to study the
	\emph{(global) $(k,w_*)$-newslopes} (Definition~\ref{k-new slopes}) and the
	\emph{(global) $k$-thresholds}  (Definition~\ref{constant-slope region thresholds}).
	The main result is the determination of the (global)
	$(k,w_*)$-newslopes as the weight $w_*$ varies in a specific open affinoid disc
	$D(k)$ centered at $k$; see Theorem~\ref{propofslopes}.
This result reveals the key role of the $k$-derivative slopes of
$\underline{\Delta}_k$ in our study. In \S~\ref{subsection k-thresholds}, we show that \emph{almost all}
(global) $k$-thresholds are given by these slopes. Consequently, the integrality
and equidistribution of the $k$-derivative slopes, established in
Proposition~\ref{nearStprop} and
Corollary~\ref{cor:equidistribution of Delta slope}, yield the corresponding
results for global $k$-thresholds; see
Corollary~\ref{cor:integrality of thresholds} and
Theorem~\ref{distribution thm for constant}.
In the final section, we apply the same strategy to slopes of $\calL$-invariants
at weight $k$, relating them to the $k$-derivative slopes of $\underline{\Delta}_k$ and
deducing their integrality and equidistribution.
	  
	Keep the notations in \S~\ref{subsection of Settings}. In particular, by Theorem \ref{local ghost theorem}, for any weight $w_*\in\gothm_{\CC_p}$, the slopes of $\NP(C_{\bar r}(w_*,-))$ coincide with that of $\NP(G(w_*,-))$ (except for the slope zero parts), up to the global multiplicity $m(\bar r)$. 
	In this section, we always assume that $k\geq 2$ is an integer, such that
	$k\equiv k_0\pmod{p-1}$.
	
	\subsection{The $(k,w_*)$-newslopes}
	\label{subsection _def}
	Recall the slopes of $p$-newforms in $S_k(\Gamma_0(Np))_{\bar r}$ are all equal to $\frac{k-2}{2}$. By Theorem \ref{local ghost theorem} and Proposition \ref{gduality}, these are the slopes of the straight segment of $\NP(C_{\bar r}(w_k,-))$ over the interval $\left(m(\bar r)\cdot d_k^{\ur},m(\bar r)\cdot(d_k^\Iw-d_k^{\ur})\right)$. Our primary interest is studying the \emph{variations} of these slopes (of $p$-newforms), namely the segments of $\NP(C_{\bar r}(w_*,-))$ over this interval for various weights $w_*$. 
	\begin{definition}
		\label{k-new slopes}
		Fix a weight $w_*\in\gothm_{\CC_{p}}$ and an integer $k\geq2$. The \emph{$(k,w_*)$-newslopes} (resp. \emph{global $(k,w_*)$-newslopes}) are the slopes of the polygon $\NP(G(w_*,-))$ (resp. $\NP(C_{\bar r }(w_*,-))$) over the interval $\left(d_k^{\ur},(d_k^\Iw-d_k^{\ur})\right)$ (resp. $\left(m(\bar r)\cdot d_k^{\ur},m(\bar r)\cdot(d_k^\Iw-d_k^{\ur})\right)$).
		
	 For $n\in\left\{1,\dots,d_k^\new\right\}$ (resp. $n\in\left\{1,\dots,m(\bar r)d_k^\new\right\}$), the \emph{$n$-th $(k,w_*)$-newslope} (resp. \emph{$n$-th global $(k,w_*)$-newslope}) is defined to be the $(d_k^{\ur}+n)$-th (resp. $(m(\bar r)\cdot d_k^{\ur}+n)$-th) slope of $\NP(G(w_*,-))$ (resp. $\NP (C_{\bar r}(w_*,-))$).  
	\end{definition}

	\subsection{The $k$-dominant disc $D(k)$}
In general, the $(k,w_*)$-newslopes are related to  $\underline{\Delta}_{k'}$ and $v_p(w_{k'}-w_*)$ for various integers $k'\equiv k_0\pmod{p-1}$ by Theorem \ref{local ghost theorem}.
In this subsection, we introduce the \emph{$k$-dominant disc}, such that for any weight $w_*$ in it, the $(k,w_*)$-newslopes depend essentially on the $k$-derivative slopes of $\underline{\Delta}_k$ and $v_p(w_*-w_k)$.
	
	\begin{definition}
		\label{good region}
		Fix an integer $k\geq2$. Consider the following valuation
		\begin{equation}
			\M(k):=\max\left\{
			v_p(w_{k'}-w_k):
			w_{k'} \neq w_k \text{ is a ghost zero of } g_n(w), n=1,\dots,d_k^{\Iw}-d_k^{\ur}
			\right\}.
		\end{equation}
The \emph{$k$-dominant disc} is defined to be the \emph{open} disc $D(k)\subseteq\gothm_{\CC_p}$, that is centered at $w_k$ with radius $p^{-\M(k)}$.
		Let $\GZ(k)$ be the set of the ghost zeros of  $g_1(w),g_2(w),\dots,g_{d_k^{\Iw}-d_k^{\ur}}(w)$. Clearly $w_*\in D(k)$ if and only if $w_k$ is the unique weight  closest to $w_*$ in $\GZ(k)$. This explains the terminology ``$k$-dominant" disc.
	\end{definition}
	\begin{example}
		\label{ex0}
		Let \(p=7\), \(a=2\), and \(\varepsilon=\omega\). In this case, we have
		\(s=s_{\varepsilon,a}=1\) and \(k_0=k_{\varepsilon,a}=6\). We refer to
		Example~\ref{Appendix ghost series} for the first eight terms of the ghost
		polynomials. For the weights \(k=6,12,18,24\), it is straightforward to check
		that
		\[
		\begin{array}{|c|c|c|c|c|}
			\hline
			k & 6 & 12 & 18 & 24 \\
			\hline
			\M(k) & 1 & 2 & 2 & 2 \\
			\hline
		\end{array}.
		\]
		For general values of \(\M(k)\), we will provide two upper bounds: one in terms of
		\(\log_p\bigl((k-k_0)/(p-1)\bigr)\), see Lemma~\ref{lemmam(k)}, and another in
		terms of the maximal \(k\)-derivative slope of \(\underline{\Delta}_k\), see
		Lemma~\ref{lem:M(k) and s_N}.
	\end{example}
	
	To provide the motivation of this definition, we perform some heuristic computations as follows. Consider $1\leq n\leq d_k^{\Iw}-d_k^{\ur}$. Since 	$v_{p}\big(g_{n}(w_{*})\big)=v_{p}\big(g_{n,\hat{k}}(w_{*})\big)+m_{n}(k)\cdot v_{p}\big(w_{*}-w_{k}\big)$, we have 
	\begin{equation}
		\label{16}
	v_{p}\big(g_{n}(w_{*})\big)=	v_{p}\big(g_{n,\hat{k}}(w_{k})\big)+m_{n}(k)\cdot v_{p}\big(w_{*}-w_{k}\big)+\sum_{k\neq k'}m_{n}(k')\cdot\big(v_{p}\big(w_{*}-w_{k'}\big)-v_{p}\big(w_{k}-w_{k'}\big)\big).
	\end{equation}
	If $w_*\in D(k)$, then the last term of \eqref{16} is equal to zero. Therefore we obtain a much simpler formula for such $w_*$:
	\begin{equation}
		\label{10}
		v_{p}\big(g_{n}(w_{*})\big)=v_{p}\big(g_{n,\hat{k}}(w_{k})\big)+m_{n}(k)\cdot v_{p}\big(w_{*}-w_{k}\big).
	\end{equation}
	From these computations, we deduce the following lemma. This lemma together with \eqref{10} suggests that, when $w_*$ belongs to the $k$-dominant disc $D(k)$, we may use the $k$-derivative slopes and their relationship to $v_p(w_*-w_k)$ to analyze the $(k,w_*)$-newslopes.
	\begin{lemma}
		\label{heuristic computation1}
		Let $w_*\in D(k)$ be a weight.
		\begin{enumerate}
			\item 
			For $n= d_{k}^{\ur}$ or $n=d_{k}^{\Iw}-d_{k}^{\ur}$, we have
			$v_{p}(g_{n}(w_{*}))=v_{p}(g_{n,\hat{k}}(w_{k}))$, which depends only on the integer $k$. 
			\item
			For $d_{k}^{\ur}\leq n<n'\leq\frac{1}{2}d_{k}^{\Iw}$, or for $\frac{1}{2}d_k^{\Iw}\leq n<n'\leq d_{k}^{\Iw}-d_{k}^{\ur}$, we have
			\begin{equation}
				\label{11}
				\begin{aligned}
					&v_{p}\big(g_{n'}(w_{*})\big)-v_{p}\big(g_{n}(w_{*})\big)\\&=
				v_{p}\big(g_{n',\hat{k}}(w_{k})\big)-v_{p}\big(g_{n,\hat{k}}(w_{k})\big)+\left\{\begin{array}{ll}(n'-n)\cdot v_p(w_*-w_k)&\mbox{if $d_{k}^{\ur}\leq n<n'\leq\frac{1}{2}d_{k}^{\Iw},$}\\-(n'-n)\cdot v_p(w_*-w_k)&\mbox{if $\frac{1}{2}d_k^{\Iw}\leq n<n'\leq d_{k}^{\Iw}-d_{k}^{\ur}.$}\end{array}\right.\\
			\end{aligned}
			\end{equation}
Moreover, these differences can be rewritten as 
			\begin{equation}
				\label{9}
					\begin{aligned}
					v_{p}\big(g_{n'}(w_{*})\big)-v_{p}\big(g_{n}(w_{*})\big)&=
				\Delta'_{k,n'-\frac{1}{2}d_k^{\Iw}}-\Delta'_{k,n-\frac{1}{2}d_k^{\Iw}}+(n'-n)\cdot\left(\frac{k-2}{2}\right)\\&+\left\{\begin{array}{ll}(n'-n)\cdot v_p(w_*-w_k)&\mbox{if $d_{k}^{\ur}\leq n<n'\leq\frac{1}{2}d_{k}^{\Iw},$}\\-(n'-n)\cdot v_p(w_*-w_k)&\mbox{if $\frac{1}{2}d_k^{\Iw}\leq n<n'\leq d_{k}^{\Iw}-d_{k}^{\ur}.$}\end{array}\right.\\
			\end{aligned}
			\end{equation}
		\end{enumerate}
	\end{lemma}
	\begin{proof}
(1) follows directly from the definition of the ghost multiplicity
	(Definition~\ref{ghostmult}), since \(m_n(k)=0\) in this case. For
	(2), Definition~\ref{ghostmult} gives
	$$
	m_n(k)=
	\begin{cases}
		n-d_k^{\ur}, & \textrm{if }d_k^{\ur}\le n\le \frac12 d_k^{\Iw},\\
		d_k^{\Iw}-d_k^{\ur}-n, &\textrm{if } \frac12 d_k^{\Iw}\le n\le d_k^{\Iw}-d_k^{\ur}.
	\end{cases}
	$$
	This proves \eqref{11}, and \eqref{9} follows from \eqref{11} and
	Definition~\ref{defofDelta}.
	\end{proof}
	\begin{remark}
		When $w_*$ does not belong to the $k$-dominant disc $D(k)$, controlling the behavior of $(k,w_*)$-newslopes becomes more complex. In this case, $w_k$ is not the unique weight closest to $w_*$ in the set of ghost zeros $\GZ(k)$, and hence this complexity is closely related to the distribution of ghost zeros. See Example \ref{ex beyond log bound} for more discussions.
	\end{remark}
	The following technical lemma proves that $\M(k)$ has an upper bound (i.e. a lower bound for the ``size" of the $k$-dominant disc $D(k)$) that increases logarithmically in $k$.
	\begin{lemma}
		\label{lemmam(k)}
	Write $k=k_0+k_{\bullet}(p-1)$ for some $k_{\bullet}\geq 0$. We use the convention that $\left\lfloor \log_p(0)\right\rfloor=0$. Then we have
		\begin{equation}
			\M(k)\leq  \left\lfloor\log_p(k_{\bullet})\right\rfloor+2.	\end{equation}
	\end{lemma}
	\begin{proof}
		Let  $k''=k_0+k''_{\bullet}(p-1)$ be the maximal integer such that $w_{k''}$ occurs as a ghost zero in $g_{d_{k}^{\Iw}}(w)$. Then 
		$
		\M(k)\leq\max\left\{v_p(w_k-w_{k'}):k'\neq k, k_0\leq k'\leq k''\  \textrm{and}\ k'\equiv k_0\pmod{p-1}\right\}
		$.
		By \eqref{valuation},  $v_p(w_k-w_{k'})=1+v_p(k-k')=1+v_p(k_{\bullet}-k'_{\bullet})$, and hence
		\begin{align*}
		\M(k)&\leq1+\max\left\{\max\{v_p(k_{\bullet}-k'_{\bullet}): 0\leq k'_{\bullet}< k_{\bullet}\},\max\{v_p(k_{\bullet}-k'_{\bullet}): k_\bullet< k'_{\bullet}\leq k''_{\bullet}\}\right\}\\
		&=1+\max\left\{\left\lfloor\log_p(k_{\bullet})\right\rfloor,\left\lfloor\log_p(k''_{\bullet}-k_\bullet)\right\rfloor\right\},
				\end{align*}
		where the last equality follows from the fact that $\left\lfloor\log_p(n)\right\rfloor$ is exactly the maximal valuation in the set $\left\{v_p(i):i=1,\dots,n\right\}$. It suffices to check $\left\lfloor\log_p(k''_{\bullet}-k_\bullet)\right\rfloor\leq\left\lfloor\log_p(k_{\bullet})\right\rfloor+1$. 
		
		We claim that $d_{k}^{\Iw}=d_{k''}^{\ur}+1$. In fact, by Definition \ref{ghostmult}, $k''$ is the maximal integer in the set
		$\left\{k'|d_{k'}^{\ur}<d_{k}^{\Iw}<d_{k'}^{\Iw}-d_{k'}^{\ur}\right\}$. By Definition \ref{dimformula}, both $d_{k'}^{\ur}$ and $d_{k'}^{\Iw}-d_{k'}^{\ur}$ are non-decreasing in the class $k'\equiv k_0\pmod{p-1}$, and 
		$
		d_{l+p-1}^{\ur}-d_{l}^{\ur}\in\{0,1\}$ for any 
		$l\equiv k_0\pmod{p-1}.$ Therefore we have $d_{k}^{\Iw}=d_{k''}^{\ur}+1$. We rewrite the dimension formulas for $d_k^{\ur}$ from Definition~\ref{dimformula} as follows. For $\ell \geq 0$, we have
		\[
		d_k^{\ur} =
		\begin{cases}
			2\ell & \text{if } \ell(p+1) \leq k_\bullet < \ell(p+1) + t_1, \\ 
			2\ell + 1 & \text{if } \ell(p+1) + t_1 \leq k_\bullet < \ell(p+1) + t_2, \\
			2\ell + 2 & \text{if } \ell(p+1) + t_2 \leq k_\bullet < (\ell+1)(p+1),
		\end{cases}
		\]
where $t_1,t_2$ are defined in Notation~\ref{notation1}.
		By Definition~\ref{dimformula}, the relation $d_k^{\Iw}=d_{k''}^{\ur}+1$ then gives
		\begin{equation}
			\label{equation1}
		(p+1)(k_\bullet-\delta)+t_1
		\leq k_{\bullet}''
		<
		(p+1)(k_\bullet-\delta)+t_2.
		\end{equation}
 Equivalently, this condition may be expressed in the following four cases
		\begin{itemize}
			\item 
			If $s+a<p-1$ and $2s+a<p-1$, then $pk_{\bullet}+s\leq k''_{\bullet}-k_{\bullet}< pk_{\bullet}+s+a+2.$
			\item 
			If $s+a<p-1$ and $2s+a\geq p-1$, then $p(k_{\bullet}-1)+s\leq k''_{\bullet}-k_{\bullet}< p(k_{\bullet}-1)+s+a+2.$
			\item 
		If $s+a\geq p-1$ and $2s+a<2p-2$, then $pk_{\bullet}+s+a-p+2\leq k''_{\bullet}-k_{\bullet}< pk_{\bullet}+s+1.$
			\item 
				If $s+a\geq p-1$ and $2s+a\geq2p-2$, then $p(k_{\bullet}-1)+s+a-p+2\leq k''_{\bullet}-k_{\bullet}< p(k_{\bullet}-1)+s+1.$
		\end{itemize}
		From these bounds, we obtain the inequality $\left\lfloor\log_p(k''_{\bullet}-k_\bullet)\right\rfloor\leq\left\lfloor\log_p(k_{\bullet})\right\rfloor+1$.
	\end{proof}
	\subsection{Calculating $(k,w_*)$-newslopes}
	\label{subsection" variation of k-newslopes}
	In this subsection, we determine the (global) $(k,w_*)$-newslopes for $w_*\in D(k)$ (Theorem \ref{propofslopes}). To do this, we introduce more notations for vertices and slopes of the $k$-derivative polygon $\underline{\Delta}_{k}^{+}$.
	\begin{notation}
		\label{notation of Delta slope}
		Given that we use the following notations for a fixed $k$, we omit all subscripts related to $k$ for simplicity.
		\begin{enumerate}
			\item 
			Let $s_{1}<s_{2}<\dots<s_{N}$ be the distinct slopes of $\underline{\Delta}_{k}^{+}$. Set $s_0=0$.
			\item
			For $i\in\left\{0,\dots,N\right\}$, denote the $i$-th vertex of $\underline{\Delta}_{k}^{+}$ by $\left(n_i,\Delta_{k,n_i}'\right)=\left(n_i,\Delta_{k,n_i}\right)$, where we set $n_0=0$. In particular, $n_N=\frac{1}{2}d_k^{\new}$ and the multiplicity of the slope $s_{i}$ is $(n_i-n_{i-1})$. 
		\end{enumerate}
	\end{notation}
	
		\begin{example}
 		\label{ex4}
		Retain the setting of Example~\ref{example1}: \(p=7\), \(a=2\), and
		\(\e=\omega\), so that \(s_{\e,a}=1\), with \(k_0=k_{\e,a}=6\).
		In this setting, for the weights \(k=6,12,18,24\), we have
		\[
		\begin{array}{|c|c|c|c|c|}
			\hline
			k & 6 & 12 & 18 & 24 \\
			\hline
			\M(k) & 1 & 2 & 2 & 2 \\
			\hline
			s_1 & 2 & 3 & 3 & 2 \\
			\hline
			s_N & 2 & 3 & 6 & 9 \\
			\hline
			N & 1 & 1 & 2 & 3 \\
			\hline
		\end{array}.
		\]
Notably, in these cases the maximal $k$-derivative slope $s_N$ is always greater than $\M(k)$. In fact, this strict inequality always holds for every \(k\) such that
\(d_k^{\new}\neq 0\); see Lemma~\ref{lem:M(k) and s_N}. By contrast, the
smallest $k$-derivative slope $s_1$ may be greater than, equal to, or less than $\M(k)$.
	\end{example}
	\begin{lemma}
		\label{lem:M(k) and s_N}
Let
$
s_N
=
\Delta_{k,\frac{1}{2}d_k^{\new}}
-
\Delta_{k,\frac{1}{2}d_k^{\new}-1}
$
denote the maximal slope of \(\underline{\Delta}_k^+\). Then  we have \(s_N>\M(k)\).
			\end{lemma}
\begin{proof}
In this proof, we freely use the notation introduced in Notation~\ref{notation1} and
Definition~\ref{dimformula}. Since \(p\geq 7\) and \(2\leq a\leq p-5\), we first record the values of \(t_1,t_2\), together with some possibly non-sharp bounds, in the following table.
\[
\begin{array}{|c|c|c|c|}
	\hline
	\text{Condition} & t_1 & t_2 &\delta\\ 
	\hline
	2s+a < p-1
	& s\geq0
	& s+a+2\in[4,p]&0 \\ 
	\hline
	2s+a \geq p-1 \text{ and } s+a < p-1
	& s+1\geq3
	& s+a+3\in[7,p+1]&1 \\ 
	\hline
	2s+a < 2(p-1) \text{ and } s+a \geq p-1
	& a+s+1-(p-1)\geq1
	& s+1\in[5,p]&0 \\ 
	\hline
	2s+a \geq 2(p-1)
	& a+s+2-(p-1)\geq3
	& s+2\in[7,p] &1\\
	\hline
\end{array}
\]
Write \(k=k_0+k_{\bullet}(p-1)\), and let
\(k''=k_0+k_{\bullet}''(p-1)\) be the maximal integer such that \(w_{k''}\)
occurs as a ghost zero in \(g_{d_k^{\Iw}}(w)\). Then, by \eqref{equation1} in the proof of Lemma~\ref{lemmam(k)}, we have
\begin{equation}
	\label{equation3}
	k_{\bullet}'' - k_\bullet <p(k_\bullet-\delta)+t_2-\delta
\end{equation}
Here, in the case \(\delta=1\), the condition \(d_k^\Iw\neq 0\) forces
\(k_\bullet\geq 1\).
Suppose now that \(\M(k)\geq s_N\). To obtain a contradiction, it suffices to show that
\begin{equation}
	\label{equation2}
	p^{\M(k)-1}\geq
		p (k_\bullet-\delta) + t_2-\delta.
\end{equation}
Indeed, since \(t_2\geq 4\), it follows directly from \eqref{equation2} that
\(p^{\M(k)-1}>k_\bullet\). Hence
$
k'=k_0+(k_\bullet+p^{\M(k)-1})(p-1)
$
is the smallest integer such that \(v_p(w_{k'}-w_k)=\M(k)\). By the definition
of \(\M(k)\), we must have \(k'\leq k''\), contradicting \eqref{equation3}.
Therefore \(s_N>\M(k)\). It remains to prove \eqref{equation2} under the
assumption that \(\M(k)\geq s_N\).

Write
$
k_\bullet=\ell_0+(p+1)\ell$ for some $\ell_0\in\{0,\dots,p\}$.
Then the dimension formula for \(d_k^{\ur}\) can be rewritten as
\begin{equation}
	\label{equation dim}
	d_k^{\ur} :=2\ell+n(k)=
	\begin{cases}
		2\ell & \text{if } 0 \leq \ell_0 < t_1, \\ 
		2\ell + 1 & \text{if } t_1\leq \ell_0 < t_2, \\
		2\ell + 2 & \text{if } t_2 \leq\ell_0 < p+1.
	\end{cases}
\end{equation}
We first verify that \eqref{equation2} holds in all cases except when
\(\delta=0\) and \(k_\bullet=\ell_0=n(k)=1\); see also
Remark~\ref{rmk:shaper bound} for an explicit example of this special case. Note that
$
s_N
=\Delta_{k,\frac{1}{2}d_k^\new}
-\Delta_{k,\frac{1}{2}d_k^\new-1}
=\Delta'_{k,\frac{1}{2}d_k^\new}
-\Delta_{k,\frac{1}{2}d_k^\new-1}
\geq
\Delta'_{k,\frac{1}{2}d_k^\new}
-\Delta'_{k,\frac{1}{2}d_k^\new-1}$. Taking \(l=\frac{1}{2}d_k^\new\) in \eqref{1233}, we obtain
\begin{equation}
	\label{equation8}
	\M(k)-1\geq\lceil s_N-1\rceil
	\geq
	\frac{p-1}{2}\left(\frac{1}{2}d_k^\new-1\right)+1
	=\frac{p-1}{2}(k_\bullet-\delta-d_k^{\ur})+1.
\end{equation}
Thus, in order to prove \eqref{equation2}, it suffices to show
\begin{equation}
	\label{equation4}
	\frac{p-1}{2}\left(k_\bullet-\delta-d_k^{\ur}\right)
	\geq
		\log_p\left(k_\bullet-\delta+\dfrac{t_2-\delta}{p}\right).
\end{equation}
Using the inequality \(\log_p(1+x)\leq x\) for \(x\geq0\), together with the
facts that \(k_\bullet-\delta\geq0\) and \(0<(t_2-\delta)/p\leq1\), it is enough to prove
$
	\frac{p-1}{2}\left(k_\bullet-\delta-d_k^{\ur}\right)
	\geq k_\bullet-\delta$.
Equivalently, we need to show
\begin{equation}
	\label{equation6}
		\dfrac{(p-1)^2}{2}\ell+\dfrac{p-1}{2}(\ell_0-\delta-n(k))
		\geq (p+1)\ell+\ell_0-\delta.
\end{equation}
It is straightforward to check that \eqref{equation6} holds when \(\ell\geq1\),
or when \(\ell=0\) and \(\ell_0-\delta>n(k)\). From the bounds for \(t_1,t_2\) and the dimension formula
\eqref{equation dim}, it remains to consider the case \(\ell=0\), with the
following possibilities
\begin{itemize}
	\item 
	\(n(k)=0\) and \(k_\bullet=\ell_0=\delta\). In this case, \eqref{equation6} is an equality.
	
	\item 
	$n(k)=1$, $\delta=0$, and	\(k_\bullet=\ell_0=0\). In this case,
	we have \(d_k^{\new}=0\).
	
	\item
	$n(k)=1$, $\delta=0$, and	\(k_\bullet=\ell_0=1\). In this case, a sharper
	bound for \(s_N\) is needed.
\end{itemize}

Consider the case \(k_\bullet=\ell_0=n(k)=1\) and \(\delta=0\).
In this case, \(d_k^\Iw=4\), \(d_k^{\ur}=1\), \(d_k^{\new}=2\), and
\(0\leq t_1\leq \ell_0=1\). Let \(\Dig(m)\) denote the sum of the digits in the \(p\)-based expansion of a
positive integer \(m\). Then, by the proof of \cite[Lemma 5.2]{Lghost}, we have
\[
\Delta'_{k,1}-\Delta'_{k,0}
=
\frac{\theta}{2}
+
\frac{\Dig(\theta)+\Dig(\eta)-\Dig(\eta+\theta)}{p-1}
=
\frac{\theta}{2}
+
\begin{cases}
	0,&\text{if }t_1=0,\\
	1,&\text{if }t_1=1,
\end{cases}
\]
where
$
\theta=p-1-a$ (resp. $a+2$) if $s+a<p-1$ (resp. $s+a\geq p-1$),
so that \(\theta\geq4\), and \(\eta+\theta=p-1+t_1\in\{p-1,p\}\).
Note that \(t_1=0\) can occur only when \(s+a<p-1\); in this case, we have
\(t_1=s=0\). Therefore
$
\left\lceil\Delta'_{k,1}-\Delta'_{k,0}\right\rceil\geq 3$,
except possibly when \(t_1=0\) and \(a=p-5\). We distinguish two cases.
\begin{itemize}
	\item If
	\(\left\lceil\Delta'_{k,1}-\Delta'_{k,0}\right\rceil\geq 3\), then we may
	replace the bound \(\M(k)-1\geq 1\) in \eqref{equation8} by the sharper bound
$
	\M(k)-1
	\geq
	\left\lceil\Delta'_{k,1}-\Delta'_{k,0}\right\rceil-1
	\geq 2$.
	Repeating the argument in the preceding paragraph, \eqref{equation6} becomes
$
	\frac{p-1}{2}(\ell_0+1-n(k))\geq \ell_0$,
	which is immediate when \(\ell_0=n(k)=1\).
	\item In the exceptional case where \(t_1=0\) and \(a=p-5\), we have
	\(k=k_0+(p-1)=2p-4\) and
	$
	s_N=\Delta'_{k,1}-\Delta'_{k,0}=2$.
	We now show directly that \(\M(k)=1\). Let
	\(\tilde{k}''=k_0+\tilde{k}''_\bullet(p-1)\) be the maximal integer such that
	\(w_{\tilde{k}''}\) occurs in \(g_{d_k^{\Iw}-d_k^{\ur}}\). Then, by the same
	argument as in the proof of Lemma~\ref{lemmam(k)}, we have
	$
	d_k^{\Iw}-d_k^{\ur}=d_{\tilde{k}''}^{\ur}+1$.
Using the dimension formulas for \(d_k^{\mathrm{Iw}}\) and \eqref{equation dim}, one checks directly that this relation implies
$
\tilde{k}''_\bullet-k_\bullet \in [t_2-1,p-1]$.
 Hence \(\M(k)=1\).
\end{itemize}
This completes the proof.
\end{proof}
\begin{remark}
	\label{rmk:shaper bound}
For most cases in the proof of Lemma~\ref{lem:M(k) and s_N}, we show that
\(\M(k)\) is strictly smaller than the lower bound
$
\frac{p-1}{2}\left(\frac{1}{2}d_k^{\new}-1\right)+2
$
for \(s_N\) obtained from \eqref{1233}. This occurs, for instance,
for \(k=6,18,24\) in Example~\ref{ex4}. When \(k=12\) in Example~\ref{ex4},
this lower bound is \(2\), which is equal to \(\M(k)\). However, in this case
\(s_N=3\), and this is precisely why we improve the lower bound to \(3\) in the
proof of Lemma~\ref{lem:M(k) and s_N}.
Finally, if \(p\geq 7\), \(s=0\), and \(a=p-5\), then for
\(k=2p-4\), we have
$
d_k^{\Iw}=4$, $d_k^{\ur}=1$, $\M(k)=1$, and $s_N=2$.
Although \(\M(k)=1\), there is a ghost zero
\(w_{p^2+p-4}\) of \(g_4=g_{d_k^{\Iw}}\) such that
$
v_p(w_{p^2+p-4}-w_k)=2=s_N$.
\end{remark}

	\begin{theorem}
		\label{propofslopes}
		Fix $i\in\left\{1,\dots,N\right\}$. Let $w_*\in D(k)$ be a weight in the $k$-dominant disc (i.e. $v_p(w_*-w_k)>\M(k)$). Then we have the following descriptions of the $(k,w_*)$-newslopes (resp. global $(k,w_*)$-newslopes). Put $\nu=v_p(w_*-w_k)$ (note that $\nu>\M(k)$).
		\begin{enumerate}
			\item 
			If $\nu\geq s_N $ (and $\nu>\M(k)$), then the (global) $(k,w_*)$-newslopes are all equal to $\frac{k-2}{2}$.
			\item
			If $s_{i-1}\leq\nu<s_{i}$ (and $\nu>\M(k)$), then the $(k,w_*)$-newslopes (resp. global $(k,w_*)$-newslopes) consist of the following three symmetric parts.
			\begin{enumerate}
				\item
				 Over the interval $\left[\frac{1}{2}d_{k}^{\Iw}-n_{i-1},\frac{1}{2}d_{k}^{\Iw}+n_{i-1}\right]$ (resp. $m(\bar r)\cdot\left[\frac{1}{2}d_{k}^{\Iw}-n_{i-1},\frac{1}{2}d_{k}^{\Iw}+n_{i-1}\right]$), the slopes are $\frac{k-2}{2}$ with multiplicity $2n_{i-1}$ (resp. $2m(\bar r)n_{i-1}$).
				\item 
				Over the interval $\left[\frac{1}{2}d_{k}^{\Iw}-n_j,\frac{1}{2}d_{k}^{\Iw}-n_{j-1}\right]$ (resp. $m(\bar r)\cdot\left[\frac{1}{2}d_{k}^{\Iw}-n_j,\frac{1}{2}d_{k}^{\Iw}-n_{j-1}\right]$) for each $j\in\left\{i,i+1,\dots,N\right\}$, the slopes are $\frac{k-2}{2}+(\nu-s_{j})$ with multiplicity $(n_j-n_{j-1})$ (resp. $m(\bar r)(n_j-n_{j-1})$).
				\item
				Over the interval $\left[\frac{1}{2}d_{k}^{\Iw}+n_{j-1},\frac{1}{2}d_{k}^{\Iw}+n_j\right]$ (resp. $m(\bar r)\cdot\left[\frac{1}{2}d_{k}^{\Iw}+n_{j-1},\frac{1}{2}d_{k}^{\Iw}+n_j\right]$)  for each $j\in\left\{i,i+1,\dots,N\right\}$, the slopes are $\frac{k-2}{2}-(\nu-s_j)$ with multiplicity $(n_j-n_{j-1})$ (resp. $m(\bar r)(n_j-n_{j-1})$).
			\end{enumerate}
		\end{enumerate}
	\end{theorem}
	\begin{remark}
	By Lemma~\ref{lem:M(k) and s_N}, we have $s_N>\M(k)$. Hence
Theorem~\ref{propofslopes} applies and determines the
$(k,w_*)$-newslopes for every weight $w_*\in D(k)$. In particular, for $w_* \in D(k)$, these slopes
depend only on the valuation $v_p(w_*-w_k)$, rather than on the value of $w_*-w_k$ itself.

		In Theorem~\ref{propofslopes}(2), it may happen that
		$0=s_0\leq \M(k)<\nu<s_1$; for instance, this occurs for the weights
		$k=6,12,18$ in Example~\ref{ex4}. In this case, none of the global
		$(k,w_*)$-newslopes for weights $w_*$ satisfying
		$v_p(w_*-w_k)=\nu$ is equal to $\frac{k-2}{2}$. Intuitively, this reflects the
		fact that such a weight $w_*$ is far from $w_k$.
	\end{remark}
By Theorem \ref{local ghost theorem}, it suffices to prove the statement in the local case. Before giving the proof, we provide an example to illustrate the \emph{symmetry} of $(k,w_*)$-newslopes for $w_{*}\in D(k)$.
	\begin{example}
		\label{example1}	
		Retain the setting of Example \ref{ex0}: \(p=7\), \(a=2\), and
		\(\e=\omega\), so that \(s_{\e,a}=1\), with \(k_0=k_{\e,a}=6\).  Let $k=24$, then it is straightforward to check that $d_{k}^{\Iw}=8,d_{k}^{\new}=6, d_{k}^{\ur}=1$, and $\M(k)=2$.
		By Remark \ref{rmk about $l<2p$}, $\{(l,\Delta'_{24,l}):l=0,1,2,3\}$ is already lower-convex (see the following chart), and hence the $k$-derivative slopes of $\underline{\Delta}_{24}^{+}$ are: $s_{1}=2,s_{2}=6$, and $s_{3}=9$ (each with multiplicity one).
		\begin{center}\renewcommand{\arraystretch}{1.2}
			\begin{tabular}{|c|c|c|c|c|}
				\hline
				$l$& $0$ &$1$ &$2$&$3$ \\
				\hline
				$\Delta'_{24,l}=\Delta_{24,l}$ & $17$ &$19$& $25$ & $34$ \\ 
				\hline
			\end{tabular}
		\end{center}
Let $w_*\in D(k)$ be a weight in the $k$-dominant disc, namely $v_p(w_*-w_k)>\M(k)=2$.
Put $r=v_{p}(w_{*}-w_{24})$, we enumerate the first eight slopes of $\NP(G(w_{*},-))$ for various values of $r>\M(k)=2$. 
		\begin{center}\renewcommand{\arraystretch}{1.2}
			\begin{tabular}{|c|c|c|c|c|c|c|c|c|}
				\hline
				$v_p(w_*-w_k)$&&\multicolumn{6}{|c|}{$(k,w_*)$-newslopes}&
				\\
				\hline
				$r\geq9$&$1$&$11$&$11$&$11$&$11$&$11$&$11$&$22$\\
				\hline
				$6\leq r\leq9$&$1$&$11-(s_3-r)$&$11$&$11$&$11$&$11$&$11+(s_3-r)$&$22$\\
				\hline
				$2<r\leq6$&$1$&$8-(s_2-r)$&$11-(s_2-r)$&$11$&$11$&$11+(s_2-r)$&$14+(s_2-r)$&$22$\\
				\hline
			\end{tabular}
		\end{center}
			Given that certain valuations in this table have specific meanings: $11=\frac{k-2}{2}$, $8=11-(s_3-s_2)$, and $14=11+(s_3-s_2)$, one can verify the statements in Theorem \ref{propofslopes} directly from the table. In particular, the $(k,w_*)$-newslopes exhibit the obvious symmetry, namely the sum of the $i$-th and the $(d_k^{\new}-i)$-th $(k,w_*)$-newslope is equal to $k-2=22$. 
	\end{example}
	\begin{proof}[Proof of (1) of Theorem \ref{propofslopes}]
		\label{proof1}
By \eqref{10} and Proposition \ref{gduality}, we have $$	v_p(g_{d_k^{\Iw}-d_k^{\ur}}(w_*))-v_p(g_{d_k^{\ur}}(w_*))=\frac{k-2}{2}(d_k^{\Iw}-d_k^{\ur}-d_k^{\ur}).$$
It suffices to check that $(n,v_p(g_n(w_*)))$ is a vertex of $\NP(G(w_*,-))$ for $n=d_k^{\ur}$ and $d_k^{\Iw}-d_k^{\ur}$, and is not a vertex for $d_k^{\ur}<n<d_k^{\Iw}-d_k^{\ur}$. These two claims can be proved using Theorem \ref{iff} as follows.

For the first claim, by Theorem \ref{iff}, it suffices to prove $(w_*,n)$ is \emph{not} near-Steinberg for any integer $k'\equiv k_0\pmod{p-1}$ when $n=d_k^{\ur}$ or $n=d_k^{\Iw}-d_k^{\ur}$. We only include the details for $n=d_k^{\ur}$ and the proof for $n=d_k^{\Iw}-d_k^{\ur}$ is similar.  Suppose that $(w_*,d_k^{\ur})$ is near-Steinberg for some integer $k'$, then we would have ($k'$ is necessarily smaller than $k$ by Definition \ref{dimformula})
\begin{equation}
	\label{15}
	d_k^{\ur}\in\left(d_{k'}^{\ur},d_{k'}^{\Iw}-d_{k'}^{\ur}\right),\ \ \textrm{and}\ \ v_{p}\big(w_{*}-w_{k'}\big)\geq \Delta_{k',|d_k^{\ur}-\frac{1}{2}d_{k'}^{\Iw}|+1}-\Delta_{k',|d_k^{\ur}-\frac{1}{2}d_{k'}^{\Iw}|}.
\end{equation}
In this case, the ghost multiplicity $m_n(k')\neq0$, and hence $\M(k)\geq v_p(w_*-w_{k'})$. Since $w_*\in D(k)$, namely $v_p(w_*-w_k)>\M(k)$, this implies $v_p(w_k-w_{k'})=v_p(w_*-w_{k'})$. As a consequence, \eqref{15} would imply that the pair $\left(w_k,d_k^{\ur}\right)$ is near-Steinberg for $k'$. By Theorem \ref{iff}, this contradicts the fact that $\left(d_k^{\ur},v_p\left(g_{d_k^{\ur}}(w_k)\right)\right)$ is a vertex of $\NP(G(w_k,-))$ (Proposition \ref{gduality}). The first claim is proved.

For the second claim, by Theorem \ref{iff}, it suffices to show $(w_*,n)$ is near-Steinberg for $k$ when $d_k^{\ur}<n<d_k^{\Iw}-d_k^{\ur}$. Since $d_k^{\ur}<n<d_k^{\Iw}-d_k^{\ur}$, we have $|n-\frac{1}{2}d_k^{\Iw}|<\frac{1}{2}d_k^{\new}$. It follows that 
		$
		v_p(w_*-w_k)\geq s_N=\Delta_{k,\frac{1}{2}d_k^{\new}}-\Delta_{k,\frac{1}{2}d_k^{\new}-1}\geq\Delta_{k,|n-\frac{1}{2}d_{k}^{\Iw}|+1}-\Delta_{k,|n-\frac{1}{2}d_{k}^{\Iw}|}$. So $(w_*,n)$ is near-Steinberg for $k$ by Definition \ref{nS def}.
	We are done.
	\end{proof}
	\begin{proof}[Proof of (2) of Theorem \ref{propofslopes}]
For $j\in\left\{i,\dots,N\right\}$, let $n'=\frac{1}{2}d_k^{\Iw}-n_{j-1}$ (resp. $\frac{1}{2}d_k^{\Iw}+n_{j}$) and $n=\frac{1}{2}d_k^{\Iw}-n_{j}$ (resp. $\frac{1}{2}d_k^{\Iw}+n_{j-1}$). Clearly $d_k^{\ur}<n<n'\leq\frac{1}{2}d_k^{\Iw}$ (resp. $\frac{1}{2}d_k^{\Iw}\leq n<n'<d_k^{\Iw}-d_k^{\ur}$).  By \eqref{9}, Proposition \ref{gduality}, and Notation \ref{notation of Delta slope} we have
			\begin{align*}
				v_{p}\big(g_{n'}(w_{*})\big)-v_{p}\big(g_{n}(w_{*})\big)&=\Delta_{k,-n_{j-1}}-\Delta_{k,-n_{j}}+\left(\frac{k-2}{2}+v_p(w_*-w_k)\right)(n_j -n_{j-1})
				\\&=\left(-s_j+\frac{k-2}{2}+v_p(w_*-w_k)\right)(n_j -n_{j-1}).
			\end{align*}
(resp.	
$
	v_{p}\big(g_{n'}(w_{*})\big)-v_{p}\big(g_{n}(w_{*})\big)=\left(s_j+\frac{k-2}{2}-v_p(w_*-w_k)\right)(n_j-n_{j-1})).
$
It suffices to prove that $(n,v_p(g_n(w_*)))$ is a vertex of $\NP(G(w_*,-))$ for $n=\frac{1}{2}d_k^{\Iw}\pm n_j$, where $j\in\left\{i-1,\dots,N\right\}$, and is not a vertex for the remaining indices $\frac{1}{2}d_k^{\Iw}-n_{j}<n<\frac{1}{2}d_k^{\Iw}-n_{j-1}$ and  $\frac{1}{2}d_k^{\Iw}+n_{j-1}<n<\frac{1}{2}d_k^{\Iw}+n_{j}$, where $j\in\left\{i-1,\dots,N\right\}$. We prove these two claims using Theorem \ref{iff} and Proposition \ref{nearStprop} as follows.
		
For the first claim, we only include the details for $\frac{1}{2}d_k^{\Iw}- n_j$ and the proof for $\frac{1}{2}d_k^{\Iw}+ n_j$ is similar. 
Let $n=\frac{1}{2}d_k^{\Iw}- n_j$ for some $j\in\left\{i-1,\dots,N-1\right\}$ (the case when $j=N$ is included in the proof of (1)). By Theorem \ref{iff}, it suffices to check that $(w_*,n)$ is \emph{not} near-Steinberg for any integer $\tilde{k}\equiv k_0\pmod{p-1}$. By Definition \ref{nS def}, this holds for $\tilde{k}=k$ since
		$
		v_p(w_*-w_k)<s_i\leq s_{j+1}=\Delta_{k,n_j+1}-\Delta_{k,n_j}.$ 
		Suppose that $(w_*,n)$ is near-Steinberg for some $k'<k$, then following the arguments in the proof of (1), we would have
		$
		(w_k,n)$ is near-Steinberg for $k'$, which contradicts Proposition \ref{nearStprop} as $(n_j,\Delta'_{k,n_j})$ is a vertex of the polygon $\underline{\Delta}_{k}^+$.
		Next suppose that $(w_*,n)$ is near-Steinberg for some $k''>k$, we will prove that $(w_*,\frac{1}{2}d_k^{\Iw}+ n_j)$ is near-Steinberg for $k''$, which again contradicts  Proposition \ref{nearStprop}. By Definition \ref{nS def}, we would have
		\begin{equation}
			\label{8}
			n\in\left(d_{k''}^{\ur},d_{k''}^{\Iw}-d_{k''}^{\ur}\right),\ \ \textrm{and}\ \ v_{p}\left(w_{*}-w_{k''}\right)\geq \Delta_{k'',|n-\frac{1}{2}d_{k''}^{\Iw}|+1}-\Delta_{k'',|n-\frac{1}{2}d_{k''}^{\Iw}|}.
		\end{equation}
	 Let $n'=\frac{1}{2}d_k^{\Iw}+ n_j$. Since $\frac{1}{2}d_k^{\Iw}<\frac{1}{2}d_{k''}^{\Iw}$, it follows that $
		|n-\frac{1}{2}d_{k''}^{\Iw}|>|n'-\frac{1}{2}d_{k''}^{\Iw}|
		$. This implies that  $\Delta_{k'',|n-\frac{1}{2}d_{k''}^{\Iw}|+1}-\Delta_{k'',|n-\frac{1}{2}d_{k''}^{\Iw}|}\geq\Delta_{k'',|n'-\frac{1}{2}d_{k''}^{\Iw}|+1}-\Delta_{k'',|n'-\frac{1}{2}d_{k''}^{\Iw}|}$ and that (together with $n\in\left(d_{k''}^{\ur},d_{k''}^{\Iw}-d_{k''}^{\ur}\right)$ ) $n'\in\left(d_{k''}^{\ur},d_{k''}^{\Iw}-d_{k''}^{\ur}\right)$.
		Consequently, we obtain
		\begin{equation}
			\label{111}
			n'\in(d_{k''}^{\ur},d_{k''}^{\Iw}-d_{k''}^{\ur}),\ \ \textrm{and}\ \ v_{p}\big(w_{*}-w_{k''}\big)\geq \Delta_{k'',|n'-\frac{1}{2}d_{k''}^{\Iw}|+1}-\Delta_{k'',|n'-\frac{1}{2}d_{k''}^{\Iw}|}.
		\end{equation}
		Once again, following the arguments in the proof of (1), we would have $(w_k,n')$ is near-Steinberg for $k''$. The first claim is proved.
		
For the second claim, it suffices to check $(w_*,\frac{1}{2}d_k^{\Iw}\pm l)$ is near-Steinberg for some $k'''$, whenever $(l,\Delta'_{k,l})$ is not a vertex of $\underline{\Delta}_k$. For such $l$'s, $(w_k,\frac{1}{2}d_k^{\Iw}\pm l)$ is near-Steinberg for some $k'''$ by Proposition \ref{nearStprop}. Since $w_*\in D(k)$, this implies that $(w_*,\frac{1}{2}d_k^{\Iw}\pm l)$ is near-Steinberg for $k'''$ (following the arguments in the proof of (1)). We are done.
	\end{proof}
	
	\subsection{The $k$-thresholds and slopes of $\underline{\Delta}_k^+$} 
	\label{subsection k-thresholds}
	By Theorem~\ref{propofslopes}, each $(k,w_*)$-newslope is equal to $\frac{k-2}{2}$ for $w_*$ sufficiently close to $w_k$. Inspired by the study of upper bounds for \emph{constant slope families} in
	\cite{const}, we introduce the following definition to characterize the
	\emph{thresholds} for this constancy. See
	Remark~\ref{rmk:Bergdall CS and my CS} for a further comparison between the two
	objects.
	\begin{definition}
		\label{constant-slope region thresholds}
		Fix an integer $k\geq 2$. For $n\in\left\{1,\dots,d_k^\new\right\}$ (resp. $n\in\left\{1,\dots,m(\bar r)d_k^\new\right\}$), the \emph{$n$-th $k$-threshold} $\CS_n(k)$ (resp. \emph{$n$-th global $k$-threshold} $\CS^{\gl}_n(k)$) is defined to be the following valuation
		\begin{align*}
			\CS^{(\gl)}_n(k):=\inf\big\{s: & \textrm{ for any weight}\ w_*\in\gothm_{\CC_p}\ \textrm{such that}\ v_p(w_*-w_k)=s,\\
			&\textrm{the}\ n\textrm{-th (global)}\ (k,w_*)\textrm{-newslope}\ \textrm{is}\ (k-2)/2 \big\}\geq0.
		\end{align*}
		Clearly $\left\{w_*\in\gothm_{\CC_{p}}:v_p(w_*-w_k)> \CS^{(\gl)}_n(k)\right\}$ is the largest \emph{open} disc centered at $w_k$, so that the $n$-th $(k,w_*)$-newslope is equal to $\frac{k-2}{2}$ for any weight $w_*$ in this disc. 
	\end{definition}
Using Theorem \ref{propofslopes}, we may identify the $k$-derivative slopes of $\underline{\Delta}_k^+$ that are larger than $\M(k)$ with \emph{almost all} of the $k$-thresholds; see the following Corollary~\ref{cor of thresholds} and Remark~\ref{rmk of log bound meaning}.
	\begin{corollary}
		\label{cor of thresholds}
		Recall that $s_N$ denotes the largest slope of
		$\underline{\Delta}_k^+$ and that $s_N>\M(k)$. Let $M$ be the smallest index in the range $i\in\left\{1,\dots,N\right\}$ for which $s_i> \M(k)$. 
		\begin{enumerate}
			\item 
		Fix $j\in\left\{M,\dots,N\right\}$. We have $\CS^{(\gl)}_n(k)=s_j$, for
			$
			n\in\left[\frac{1}{2}d_k^{\new}-n_j,\frac{1}{2}d_k^{\new}-n_{j-1}\right)$ or $n\in\left(\frac{1}{2}d_k^{\new}+n_{j-1},\frac{1}{2}d_k^{\new}+n_{j}\right]
			$ (resp. $
			n\in m(\bar r)\cdot\left[\frac{1}{2}d_k^{\new}-n_j,\frac{1}{2}d_k^{\new}-n_{j-1}\right)$ or $n\in m(\bar r)\cdot\left(\frac{1}{2}d_k^{\new}+n_{j-1},\frac{1}{2}d_k^{\new}+n_{j}\right]$).
			\item
			All the remaining $k$-thresholds (resp. global $k$-thresholds), namely when $n$ lies in the interval $\left[\frac{1}{2}d_k^{\new}-n_{M-1},\frac{1}{2}d_k^{\new}+n_{M-1}\right]$ (resp. $m(\bar r)\cdot\left[\frac{1}{2}d_k^{\new}-n_{M-1},\frac{1}{2}d_k^{\new}+n_{M-1}\right]$), satisfy
			$
			\CS^{(\gl)}_n(k)\leq \M(k)\leq\left\lfloor\log_p\left(\frac{k-k_0}{p-1}\right)\right\rfloor+2.
			$
			\item
			The number of the remaining $k$-thresholds (resp. global $k$-thresholds), namely $2n_{M-1}$ (resp. $2m(\bar r)\cdot n_{M-1}$), has the following upper bound 
			\begin{equation}
				\label{23}
				2n_{M-1}\leq C(k,a) \  \textrm{(resp.}\ 	2m(\bar r)n_{M-1}\leq m(\bar r)C(k,a)\ \textrm{)},
			\end{equation}
			where the constant $C(k,a)$, depending on $k$ and $a\in\{2,\dots,p-5\}$, is given by
			\[
			C(k,a)=\frac{2(p+3-\min\{a+2,p-1-a\})+4\left\lfloor\log_p\left(\frac{k-k_0}{p-1}\right)\right\rfloor}{p-1}\leq2+\frac{4\left\lfloor\log_p\left(\frac{k-k_0}{p-1}\right)\right\rfloor}{p-1}.\]
		\end{enumerate}
	\end{corollary}	
	\begin{proof}
		(1) follows directly from Theorem \ref{propofslopes}. For (2), from Theorem \ref{propofslopes}, we deduce that all the remaining $k$-thresholds satisfy $\CS^{(\gl)}_n(k)\leq \M(k)\leq	\left\lfloor\log_p\left(\frac{k-k_0}{p-1}\right)\right\rfloor+2$, where the last inequality is due to Lemma \ref{lemmam(k)}. For (3), note that 
			\begin{align*}
			\left\lfloor\log_p\left(\frac{k-k_0}{p-1}\right)\right\rfloor+2&\geq\M(k)\geq s_{M-1}=\Delta'_{k,n_{M-1}}-\Delta_{k,n_{M-1}-1}\geq\Delta'_{k,n_{M-1}}-\Delta'_{k,n_{M-1}-1}\\
			&\geq\frac{\min\{a+2,p-1-a\}}{2}+\frac{p-1}{2}(n_{M-1}-1),
			\end{align*}
		where the last inequality is included in \eqref{1233}. Now (3) is clear.
	\end{proof}
	\begin{remark}
		\label{rmk of log bound meaning}
If \(s_1>\M(k)\), then Corollary~\ref{cor of thresholds} determines all
(global) \(k\)-thresholds. For instance, in Example \ref{ex4}, we observe that \(s_1>\M(k)\) for \(k=6,12,18\).
Even when this inequality fails,  the $k$-derivative polygon $\underline{\Delta}_k$ calculates \emph{almost all} of the (global) $k$-thresholds, except for a \emph{logarithmic error term} in terms of $k$.
More precisely, the logarithmic error term has a twofold meaning.
\begin{enumerate}
	\item 
By Definition~\ref{dimformula}, we have \(d_k^{\new}=O(k)\), and hence the
total number of (global) \(k\)-thresholds grows linearly in \(k\). In comparison,
by \eqref{23}, the number of (global) \(k\)-thresholds that are not
\(k\)-derivative slopes grows \emph{at most} logarithmically in \(k\).
\item
	These exceptional (global) $k$-thresholds have an upper bound $\M(k)$, which grows logarithmically in $k$ by Lemma  \ref{lemmam(k)}.
\end{enumerate}
Returning to Example~\ref{ex4}, the weight $k=24$ is the smallest weight in the congruence class $k\equiv 6 \pmod{p-1}$ for which the logarithmic error term appears; see Example~\ref{ex beyond log bound}.
	\end{remark}
The remaining $k$-thresholds, as well as the $(k,w_*)$-newslopes when $w_*$ does not belong to the $k$-dominant disc $D(k)$, depend on the distribution of ghost zeros in $\GZ(k)$ in a subtle manner. Intuitively, whenever $w_*\notin D(k)$, other ghost zeros come into play. We present the following example to illustrate the prediction of the remaining $k$-thresholds.
\begin{example}
	\label{ex beyond log bound}
	We continue with the case in Example \ref{example1}, where $k=24$, $d_k^{\new}=6$, and the $k$-derivative slopes of $\underline{\Delta}_{24}$ are given by $s_3=9$, $s_2=6$, and $s_1=2$ (each with multiplicity one). In this case the integer $M$ in Corollary \ref{cor of thresholds} is equal to $2$ and hence
	\[
	\CS_1(24)=\CS_6(24)=9, \ \textrm{and}\ \CS_2(24)=\CS_5(24)=6.
	\] 
	Let $r=v_{p}(w_{*}-w_{24})$. We list the first eight slopes of $\NP(G(w_{*},-))$ for various values of $r\in\left(0,1\right)\cup\left(1,2\right)$. Clearly these weights $w_*$'s do not belong to the $k$-dominant disc $D(k)$.
	\begin{center}\renewcommand{\arraystretch}{1.2}
		\begin{tabular}{|c|c|c|c|c|c|c|c|c|}
			\hline
			$v_p(w_*-w_k)$&&\multicolumn{6}{|c|}{$(k,w_*)$-newslopes}&
			\\
			\hline
			$1< r< s_1 $ & $1$ &$2+r$& $5+r$ & $7+2r$ &$11$ &$13+r$&$14+2r$&$16+3r$ \\ 
			\hline
			$0<r<1 $ & $r$ &$3r$& $6r$ & $9r$ &$11r$ &$14r$&$16r$&$19r$ \\ 
			\hline
		\end{tabular}
	\end{center}
	From this table, it is evident that 
	$\CS_3(24)=2=s_1,\ \textrm{but}\ \CS_4(24)=1<s_1$. Notably the symmetry of $(k,w_*)$-newslopes mentioned in Example \ref{example1} fails when $0<r< s_1$, namely the sum of the third and the fourth $(k,w_*)$-newslope is not equal to $k-2=22$. In fact, when $0< r< s_1$, there are \emph{at least} two ghost zeros $w_k=w_{24}$ and $w_{66}$ that are closest to $w_*$, among all the ghost zeros in $\GZ(24)$. The appearance of these additional ghost zeros makes it more involved to calculate the $k$-thresholds.  
\end{example}
	\begin{remark}
	\label{rmk:Bergdall CS and my CS}
Let $f\in S_k(\Gamma_0(Np))$ be a $p$-newform.
In \cite{const}, Bergdall studies the \emph{constant slope radius}
$\CS^k(f)$ for constant slope families passing through $f$;
see \cite[Definitions~3.1 and~4.2]{const}. Roughly speaking,
$\CS^k(f)$ measures how far such a family can extend in weight space.
A key feature of the construction in \cite{const} is that a constant
slope family $U\subseteq\calC_N$ is required to be rigid analytically
\emph{isomorphic} to its image in weight space.
By contrast, our definitions of the $(k,w_*)$-newslopes and
$k$-thresholds are purely slope-theoretic. They are defined solely in
terms of the Newton polygon
$\NP(C_{\bar r}(w,t))$; see Definition~\ref{k-new slopes}. In particular,
they do not explicitly involve geometric properties of the eigencurve,
such as the ramification of the weight map.

A notable difference is that the construction in \cite{const} studies
constant slope families passing through \emph{individual} $p$-newforms
$f\in S_k(\Gamma_0(Np))_{\bar r}$, whereas our construction treats the
slopes of these families simultaneously through the variation of the
segments of $\NP(C_{\bar r}(w,t))$ corresponding to $p$-newforms.
Consequently, the precise relation between $\CS^k(f)$ and our
global $k$-thresholds remains unclear to us.
Indeed, even if a fixed constant slope family
$U$ contributes to the $(k,w_*)$-newslopes for every
$w_*\in\wt(U)$, it is not clear that the corresponding slope can be
tracked by a fixed index of the Newton polygon as the weight varies.
More precisely, if $w_*,w_*'\in\wt(U)$ are distinct, the corresponding
roots of $C(w_*,t)$ and $C(w_*',t)$ need not give rise to segments with
the same index on $\NP(C(w_*,t))$ and $\NP(C(w_*',t))$.
\end{remark}

\subsection{Integrality and equidistribution of \(k\)-thresholds}
	\label{section_dis}
In this subsection, we give two consequences of Corollary \ref{cor of thresholds} concerning the integrality and equidistribution of (global) $k$-thresholds. The arguments will be reused in \S~\ref{subsection : proof of main theorem} to prove the integrality and equidistribution of the slopes of $\calL$-invariants.
\begin{corollary}
	\label{cor:integrality of thresholds}
	Recall the constant $C(k,a)$ defined in Corollary~\ref{cor of thresholds}(3). 
Then with at most $\lceil C(k,a)\rceil $ (resp. $\lceil m(\bar r)C(k,a)\rceil $) exceptions, each (global) $k$-threshold $\CS_n^{(\gl)}(k)$ satisfies the following properties.
\begin{enumerate}
	\item If $\CS_n^{(\gl)}(k)$ occurs with multiplicity one in $\underline{\Delta}_k$, then $\CS_n^{(\gl)}(k)\in\ZZ$.
	\item If $\CS_n^{(\gl)}(k)$ occurs with multiplicity greater than one in $\underline{\Delta}_k$, then  $\CS_n^{(\gl)}(k)\in\frac{k}{2}+\ZZ$.
\end{enumerate}	
\end{corollary}
\begin{proof}
	This follows from Corollary \ref{cor of thresholds} and Proposition \ref{nearStprop}.
\end{proof}

	Although we do not determine all of the (global) \(k\)-thresholds, 
	Corollary~\ref{cor of thresholds}, together with the equidistribution result for $k$-derivative slopes in
	Corollary~\ref{cor:equidistribution of Delta slope}, suffices to prove the
	following equidistribution property for the (global) \(k\)-thresholds as
	\(k\to\infty\).
		\begin{theorem}
		\label{distribution thm for constant}
	Let \(\nu_k\) (resp. $\nu_k^{\gl}$) be the uniform probability measure of the following
	multiset
		\begin{equation}
\label{3eq}
		\left\{
		\frac{2(p+1)}{(p-1)k}\CS_n(k)
		:
		n=1,\dots,d_k^{\new}
		\right\}\textrm{ (resp. }	\left\{
		\frac{2(p+1)}{(p-1)k}\CS_n^{\gl}(k)
		:
		n=1,\dots,m(\bar r)\cdot d_k^{\new}
		\right\}\textrm{)}.
		\end{equation}
		Then \(\nu_k\) (resp. $\nu_k^{\gl}$) weakly converges to the uniform probability
		measure on \([0,1]\) as \(k\to\infty\).
	\end{theorem}
	\begin{proof}
	It suffices to deal with the local case. Write the \(k\)-derivative slopes of \(\underline{\Delta}_k\), respectively the
	elements of \eqref{3eq}, in ascending order as
	\[
	s_{k,1}\leq \cdots \leq s_{k,d_k^{\new}},
	\qquad
	x_{k,1}\leq \cdots \leq x_{k,d_k^{\new}}.
	\]
	By Corollary~\ref{cor:equidistribution of Delta slope} and its proof, it
	suffices to show that
	$
	\left|
	\frac{2(p-1)}{p+1}\frac{s_{k,i}}{k}
	-
	x_{k,i}
	\right|$
	tends to zero uniformly in \(i\) as \(k\to\infty\).
	By Corollary~\ref{cor of thresholds}, this difference is equal to zero
	for all but the possible exceptional cases. In those exceptional cases, both
	quantities satisfy
	\[
	0\leq
	\frac{2(p-1)}{p+1}\frac{s_{k,i}}{k},
	\ x_{k,i}
	\leq
	\frac{2(p+1)}{p-1}\frac{\M(k)}{k}
	=
	\frac{O(\log_p k)}{k}.
	\]
	Hence 
	$
	\left|
	\frac{2(p-1)}{p+1}\frac{s_{k,i}}{k}
	-
	x_{k,i}
	\right|=\frac{O(\log_p k)}{k}$,
	which tends to zero uniformly in \(i\) as  \(k\to\infty\). This completes the proof.
	\end{proof}

	\section{Greenberg--Stevens methods in higher ranks}\label{chapter: 3}
	\label{section of Muitiple GS formula}
	In this section, we aim to develop a higher-rank version of the Greenberg--Stevens formula, which connects $\mathcal{L}$-invariants associated with $\bar r$-newforms of weight $k$ to the eigenvalues of the \emph{first-derivative matrix} of the $U_p$-action, which we denote by $A_1$. We also briefly introduce the key ideas for computing the slopes of $A_1$, which will be carried out in \S~\ref{section of slopes for A_1}.
	
	We continue with $\bar r$, $\bar \rho$, $\e$, $\calC_{\bar r}$, and $m(\bar r)$ as specified in \S~\ref{subsection of Settings} . Fix an integer $k$, and keep the notations for the slopes of the derivative polygon $\underline{\Delta}_k$ in Notation \ref{notation of Delta slope}. Let $d=m(\bar r)d_k^\new$ be the dimension of the space $S_k(\Gamma_0(Np))_{\bar r}$.
	
	\subsection{Recollection of the eigencurve} Let  $\Spc_{\bar r}$ be the ($\bar r$-) \emph{spectral curve} in the sense of \cite{buzzard}, namely the hypersurface in $\calW\times\GG_m^{\rig}$ defined by the characteristic power series $C_{\bar r}(w,t)\in\calO[\![w,t]\!]$ of the $U_p$-action (see \S~\ref{subsection local ghost theorem}).
	Recall the \emph{weight map}
	$
	\wt:\calC_{\bar r}\rightarrow \Spc_{\bar r}\rightarrow \calW
	$ and the \emph{slope map} $a_p:\Spc_{\bar r}\rightarrow \GG_m^{\rig}$.
	
	\begin{lemma}
		\label{lem:separated}
		The weight map $\wt:\calC_{\bar r}\rightarrow\calW$ is separated.
	\end{lemma}
	\begin{proof}
		First, the natural map $\calC_{\bar r}\rightarrow\Spc_{\bar r}$ is
		finite by \cite[Lemma~5.3]{buzzard}, and hence separated. Moreover, by the construction of the spectral curve $\Spc_{\bar r}$
		(the space denoted by ``$Z$'' in \cite[p.~78]{buzzard}), the weight map
		$
		\wt:\Spc_{\bar r}\longrightarrow \calW
		$
		is locally finite. More precisely, for every affinoid subdomain appearing
		in the admissible covering ``$\mathcal{C}$" of \cite[\S~4]{buzzard}, the
		restriction of the weight map to this affinoid is finite. It follows that
		$\wt:\Spc_{\bar r}\to\calW$ is separated. Hence the composition
		$
		\wt:\calC_{\bar r}\longrightarrow \calW
		$
		is also separated.
	\end{proof}
	
	Let $\Spm \calA=\calD\subset \calW$ be an affinoid subdomain. We use $\calC_{\calD}$ (resp. $\Spc_{\calD}$) to denote the preimage of the disc $\calD$ inside $\calC_{\bar r}$ (resp. $\Spc_{\bar r}$). Let $\calS_{\calD}$ be the $\calA$-Banach module arising from the \emph{eigenvariety machine} in the sense of \cite{buzzard}, and it corresponds to the (eigensystems of) $\bar r$-overconvergent modular forms of weights in $\calD$ (\cite[Lemma 5.9]{buzzard}). 
	Following Buzzard's construction, $\calS_{\calD}$ also corresponds to a coherent sheaf $\calM_\calD$ over $\calC_{\calD}$. In fact, over each affinoid subdomain of the admissible cover ``$\calC$" in \cite[\S~4]{buzzard} of $\calW$, the corresponding module ``$N$" from \cite[p. 84]{buzzard} is finite.
	
	\subsection{The $(k,w_*)$-newslope component over closed discs in $D(k)$} 
	\label{subsection newslope comp}
	Let $\Spm \calA=\calD\subseteq D(k)$ be a closed disc centered at $w_k$.
	In this subsection, we apply Kiehl's Proper Mapping Theorem
	\cite[\S~6.3, Theorem~9]{bosch-rigid-formal} to construct a finite
	$\calA$-module that is flat of constant rank $d$ and whose fiber at
	$w_*\in\calD$ parameterizes the $(k,w_*)$-newslopes. As a consequence, after possibly shrinking $\calD$, the $U_p$-action on this module can be expressed as a $d\times d$ matrix with entries in $\calA$.
	
	\begin{definition}
		\label{definition: eigencurve cut}
Fix a closed disc $\calD$ centered at $w_k$ and contained in the
$k$-dominant disc $D(k)$. Let $\lambda_1(\calD)$ and $\lambda_2(\calD)$
be the greatest lower bound and the least upper bound, respectively, of
the $(k,w_*)$-newslopes as $w_*$ ranges over $\calD$. We define $U_{\calD}$ to
be the following affinoid domain in the spectral curve $\Spc_{\calD}$
\[
U_{\calD}
:=
\Spc_{\calD}^{[\lambda_1(\calD),\lambda_2(\calD)]}
=
\left\{
(w,a_p)\in\wt^{-1}(\calD):
v_p(a_p)\in[\lambda_1(\calD),\lambda_2(\calD)]
\right\}.
\]
Similarly, for any $\epsilon>0$, define
\[
U^\epsilon_{\calD}
:=
\Spc_{\calD}^{[\lambda_1(\calD)-\epsilon,\lambda_2(\calD)+\epsilon]}
=
\left\{
(w,a_p)\in\wt^{-1}(\calD):
v_p(a_p)\in[\lambda_1(\calD)-\epsilon,\lambda_2(\calD)+\epsilon]
\right\}.
\]
We denote by $\widetilde{U}_{\calD}$ (resp. $\widetilde{U}^{\epsilon}_{\calD}$) the preimage of $U_{\calD}$
(resp. $U^{\epsilon}_{\calD}$) in $\calC_{\calD}$ under the map
$\calC_{\calD}\to\Spc_{\calD}$. Note that $\widetilde{U}_{\calD}$ (resp.  $\widetilde{U}^{\epsilon}_{\calD}$) is affinoid since $\calC_{\calD}\to\Spc_{\calD}$ is finite  by \cite[Lemma 5.3]{buzzard}. 
		\end{definition}
		\begin{remark}
			Let $r(\calD)$ be the valuation radius of the closed disc $\calD$ in Definition~\ref{definition: eigencurve cut}. Then
			$r(\calD)>\M(k)$. By Theorem~\ref{propofslopes}, the greatest lower bound $\lambda_1(\calD)$ and the least upper bound $\lambda_2(\calD)$ in Definition \ref{definition: eigencurve cut}, are given respectively by
			\[
			\lambda_1(\calD)
			=
			\min\left\{
			\frac{k-2}{2}-s_N+r(\calD),\frac{k-2}{2}
			\right\},\textrm{ and }
			\lambda_2(\calD)
			=
			\max\left\{
			\frac{k-2}{2}+s_N-r(\calD),\frac{k-2}{2}
			\right\}.
			\]
			Moreover, both bounds are attained for any $w_*\in\calD$ satisfying
			$v_p(w_*-w_k)=r(\calD)$.
		\end{remark}
In what follows, we verify that, for $\calD$ as in
Definition~\ref{definition: eigencurve cut}, the weight map
$
\wt:\widetilde{U}_{\calD}\longrightarrow \calD
$
is \emph{proper} in the sense of
\cite[\S~6.3, Definition~8]{bosch-rigid-formal}, namely
\begin{enumerate}
	\item the weight map
	$
	\wt:\widetilde{U}_{\calD}\longrightarrow \calD
	$
	is separated;
	
	\item there exist two finite admissible affinoid coverings
	$(U_i)_{i\in I}$ and $(U'_i)_{i\in I}$ of $\widetilde{U}_{\calD}$
	such that $U_i\subseteq U'_i$ is relatively compact for all
	$i\in I$, in the sense of
	\cite[\S~6.3, Definition~6]{bosch-rigid-formal}.
\end{enumerate}
This will allow us to apply Kiehl's Proper Mapping Theorem
\cite[\S~6.3, Theorem~9]{bosch-rigid-formal} to the weight map in the proof of
Proposition~\ref{prop:newslope component}.
	\begin{lemma}
		\label{lem eigencurve cut}
		Let $1\leq n<m$ be two integers. Suppose that for any weight $w_*\in\calD$, both $n$ and $m$ correspond to vertices of $\NP(C_{\bar r}(w_*,t))$. Then there exist $\epsilon>0$, such that for any $w_*\in\calD$, the difference between the $n$-th (resp. $(m+1)$-th) and the $(n-1)$-th (resp. $m$-th) slope of $\NP(C_{\bar r}(w_*,t))$ is greater than $\epsilon$. 
	\end{lemma}
	\begin{proof}
			 This follows from \cite[Lemma 2.26]{Lghost2} and the proof of
			 \cite[Corollary 2.27(2)]{Lghost2}. It is enough to prove the corresponding
			 statement for $\NP(G(w_*,t))$.
			 
			 Let $\calD^{\Berk}$ denote the Berkovich space associated with $\calD$. In
			 particular, $\calD^{\Berk}$ is compact. It remains to show that the following
			 two functions from $\calD^{\Berk}$ to $\RR$ are continuous and positive:
			 \begin{align*}
			 	z &\longmapsto
			 	\text{the difference between the } n\text{-th and the }(n-1)\text{-st slopes of }
			 	\NP(G(z,t)), \\
			 	z &\longmapsto
			 	\text{the difference between the } (m+1)\text{-st and the }m\text{-th slopes of }
			 	\NP(G(z,t)),
			 \end{align*}
where the Newton polygon $\NP(G(z,t))$ for $z\in\calD^{\Berk}$ is defined in \cite[Notation 2.25]{Lghost2}.
Continuity follows directly from \cite[Lemma 2.26]{Lghost2}. Note that our assumption ensures that
	$
	\calD \subseteq \Vtx_n^{(\e)} \cap \Vtx_m^{(\e)}
	$
	in the notation of \cite[Corollary 2.27]{Lghost2}; hence positivity follows
	from \cite[Corollary 2.27(1)]{Lghost2}.
	\end{proof}
	\begin{example}
		\label{ex5}
		We apply Lemma \ref{lem eigencurve cut} to the following special case. Let $\calD$ be a closed disc centered at $w_k$ inside the $k$-dominant disc $D(k)$ as in Definition~\ref{definition: eigencurve cut}. For any $w_*\in\calD$, we have shown in the proof of Theorem \ref{propofslopes} that  $n=m(\bar r)d_k^{\ur}$ and $m=m(\bar r)\left(d_k^{\Iw}-d_k^{\ur}\right)$ correspond to vertices of $\NP(C_{\bar r}(w_*,t))$. From this we obtain certain $\epsilon>0$ as claimed in Lemma \ref{lem eigencurve cut}. As a consequence, we have the following key properties
		\begin{equation}
			\label{eq U=U'}
			U_{\calD}=U_{\calD}^{\epsilon}
			\qquad\text{and}\qquad
			\widetilde{U}_{\calD}=\widetilde{U}_{\calD}^{\epsilon}.
		\end{equation}
This implies that  $\widetilde{U}_{\calD}\subseteq\widetilde{U}_{\calD}^{\epsilon}$ is \emph{relatively compact} in the sense of \cite[\S~6.3, Definition~6]{bosch-rigid-formal}.
	\end{example}

\begin{remark}
	\label{rmk eigenvurve cut}
Assume that \(d_k^{\ur}>0\) and that \(\M(k)\) is equal to the following
valuation:
\[
\M'(k):=\max\left\{
v_p(w_{k'}-w_k):
w_{k'} \neq w_k \text{ is a ghost zero of } g_n(w),\ 
n=1,\dots,d_k^{\Iw}
\right\}.
\]
Then one can obtain an effective constant \(\epsilon>0\) as follows. Let
\(0\leq n_1<d_k^{\ur}\) and \(d_k^{\Iw}-d_k^{\ur}<n_2\leq d_k^{\Iw}\)
be, respectively, the greatest and the smallest indices corresponding to
vertices of \(\NP(G(w_k,t))\).
Arguing as in the proof of Theorem~\ref{propofslopes}, one checks that
\(n_1\) and \(n_2\) are also vertices of \(\NP(G(w_*,t))\) for
\(w_*\in D(k)\). Hence the \(\bigl(d_k^{\ur}-1\bigr)\)-th slope and the
\(\bigl(d_k^{\Iw}-d_k^{\ur}+1\bigr)\)-th slope of the polygon
\(\NP(G(w_*,t))\) are independent of $w_*\in D(k)$. Let \(\eta_1\) and \(\eta_2\) denote these
two slopes, respectively.
It follows that
$
\lambda_1(\calD)>\eta_1$
and $\lambda_2(\calD)<\eta_2$.
Thus we may take $\epsilon$ to be any positive number in
$
\left(0,\min\{\lambda_1(\calD)-\eta_1,\eta_2-\lambda_2(\calD)\}\right)$.

In Example \ref{example1}, we have $k=24$, $s_N=s_3=9$, $\M(k)=2$,
	$\eta_1=1$, and $\eta_2=22$. Hence
	$\epsilon\in\left(0,\min\left\{r(\calD)+1,10\right\}\right)
	$,
	where $r(\calD)>2$.
\end{remark}

\begin{corollary}
\label{properness}	
Let $\calD$ be a closed disc centered at $w_k$ inside the $k$-dominant disc $D(k)$. Let $\widetilde{U}_{\calD}$ be the affinoid subdomain in Definition \ref{definition: eigencurve cut}. Then the weight map $
\wt:\widetilde{U}_{\calD}\longrightarrow \calD
$
is proper in the sense of
\cite[\S~6.3, Definition~8]{bosch-rigid-formal}.
\end{corollary}
\begin{proof}
This follows from Lemma \ref{lem:separated} and \eqref{eq U=U'} by taking $\widetilde{U}_{\calD}$ and $\widetilde{U}_{\calD}^{\epsilon}$ to be the two admissible coverings of $\widetilde{U}_{\calD}$.
\end{proof}

	\begin{proposition}
			\label{prop:newslope component}
		Fix a closed disc $\calD$ centered at $w_k$ inside the $k$-dominant disc $D(k)$. 
			Let $\calM_{\widetilde{U}_{\calD}}$ be the restriction to the affinoid subdomain $\widetilde{U}_{\calD}$ of the coherent sheaf $\calM_\calD$ over $\calC_\calD$.
			\begin{enumerate}
				\item The sheaf $\wt_*(\calM_{\widetilde{U}_{\calD}})$ over $\calD=\Spm\calA$ corresponds to a finite $\calA$-module $\calS_{\calD}^{\new}$, which is flat of constant rank $d=m(\bar r)d_k^\new$. 
				\item There is a well-defined degree d characteristic polynomial $C^{\new}(w,t)\in\calA[t]$ of the $U_p$-action on $\calS_{\calD}^{\new}$, such that the $U_p$-slopes at each $w_*\in\calD$ are exactly the $(k,w_*)$-newslopes. 
			\end{enumerate}
	\end{proposition}
\begin{proof}
(2) follows directly from (1). Indeed, flatness implies local freeness, and hence $C^{\new}(w,t)$ is well-defined locally on some (Zariski) open covering of $\calD$. On the overlaps of two of these opens, the transition maps are given by invertible matrices, which do not affect the characteristic polynomial. Therefore $C^{\new}(w,t)\in\calA[t]$ is well-defined. By the construction of $\widetilde{U}_{\calD}$, the slopes of $C^{\new}(w,t)$ at $w_*\in\calD$ are precisely the $(k,w_*)$-newslopes.

We will prove (1) as follows. By Kiehl's Proper Mapping Theorem
\cite[\S~6.3, Theorem~9]{bosch-rigid-formal}, the properness in Corollary~\ref{properness} implies that $\wt_*(\calM_{\widetilde{U}_{\calD}})$ is a coherent sheaf over $\calD$. So we obtain a finite $\calA$-module $\calS_{\calD}^{\new}$. Since $\calD$ is irreducible and the rank of $\calS_{\calD}^{\new}$ at $w_*=w_k$ is $d$ (the dimension of the space $S_k(\Gamma_0(Np))_{\bar r}$), it suffices to prove the flatness. The flatness essentially follows from the construction (in the sense of Buzzard) of the eigencurve.

In fact, when $U_{\calD}$ is an element of the admissible covering ``$\calC$" in \cite[\S~4]{buzzard}, then $\calM_{\widetilde{U}_{\calD}}$ corresponds to a projective  (and hence flat) module ``$N$" in the discussion of \cite[p. 84]{buzzard}. So we need to show $U_{\calD}$ is such a special affinoid subdomain, namely the map $U_{\calD}\rightarrow D$ is finite and surjective, and $U_{\calD}$ is disconnected from its complement in $\Spc_{\calD}$. Surjectivity is trivial and finiteness follows from \cite[Corollary 4.3]{buzzard} (note that $U_{\calD}\rightarrow \calD$ has constant degree by the definition of $U_{\calD}$). By \eqref{eq U=U'}, we have $U_{\calD}=U_{\calD}^{\epsilon}$. Hence the slopes appearing in $\Spc_{\calD}\setminus U_{\calD}$ are
uniformly separated from the slopes appearing in $U_{\calD}$. Consequently, $U_{\calD}$ is disconnected from its complement in $\Spc_{\calD}$. We are done.
\end{proof}
As a direct consequence, when $\Spm\calA=\calD\subseteq D(k)$ is sufficiently small, the $\calA$-module $\calS_\calD^\new$ is finite free  of rank $d$. By choosing a basis of  $\calS_\calD^\new$, the $U_p$-action is given by a matrix in $\M_{d\times d}(\calA)$. In the following subsection, we diagonalize this matrix by choosing sufficiently small $\calD$  and selecting an appropriate basis.

\subsection{Slope functions $a_{p,i}$ and the Taylor expansion of $U_p^2$}
Let $f_1,\dots,f_d$ be the $d$ distinct $\bar r$-newforms of weight $k$. Each $f_i$ corresponds to a point $x_{f_i}$ on $\calC_{\bar r}$, and the weight map $\wt:\calC_{\bar r}\rightarrow\calW$ is \emph{étale} at $x_{f_i}$ (\cite[Proposition 2.6]{const}). As a consequence, there exist \emph{\(p\)-adic families of constant slope} $U_i\subseteq \calC_{\bar r}$, in the sense of \cite[Definition 3.1]{const}, over some sufficiently small closed disc $\Spm\calA^{\textrm{\'et}}=\mathcal{D}^{\textrm{\'et}}\subseteq D(k)$ \emph{centered at $w_k$}. Each $U_i$ is an affinoid open neighborhood of $x_{f_i}$ and the weight map $\wt:U_i\rightarrow \mathcal{D}^{\textrm{\'et}}$ is a rigid analytic isomorphism. 
By \cite[Lemma 5.9]{buzzard} and the work of Coleman (\cite{Coleman}), each $U_i$ corresponds to a formal \emph{$q$-expansion} $
\mathbf{f}_i(w)= \sum_{n \geq 1} a_{n,i}(w)q^n$, where all the coefficient functions \(a_{n,i}(w)\)'s are rigid analytic functions in $\calA^{\textrm{\'et}}$ and $\mathbf{f}_i(w_k)=f_i$. In particular, these $q$-expansions are $U_p$-eigenfunctions:
$
U_p\left(\bff_{i}(w_*)\right)=a_{p,i}(w_*)\cdot\bff_{i}(w_*), \textrm{for}\ w_*\in\calD^{\textrm{\'et}}.
$ 
\begin{lemma}
	\label{lemma: basis}
	For a sufficiently small $\calD^{\textrm{\'et}}$, the finite flat $\calA^{\textrm{\'et}}$-module $\calS^{\new}_{\calD^{\textrm{\'et}}}$ constructed in Proposition \ref{prop:newslope component} is free and the $q$-expansions $\{\bff_i(w):i=1,\dots,d\}$ form a basis of $\calS^{\new}_{\calD^{\textrm{\'et}}}$. In particular, the matrix of the $U_p$-action on $\calS^{\new}_{\calD^{\textrm{\'et}}}$ under this basis is the following diagonal matrix $\diag(a_{p,1}(w),\cdots,a_{p,d}(w))$.
\end{lemma}
\begin{proof}
	After possibly shrinking $\calD^{\textrm{\'et}}$, we may assume $\calS^{\new}_{\calD^{\textrm{\'et}}}$ is a free $\calA^{\textrm{\'et}}$-module of rank $d$, and $U_i$'s are \emph{disjoint} as the weight map is \'etale at each $f_i$. As a consequence, the support of $\calS^{\new}_{\calD^{\textrm{\'et}}}$ on $\calC_{\calD^{\textrm{\'et}}}$ (i.e. the ``$U_{\calD}$" in Definition \ref{definition: eigencurve cut} for $\calD=\calD^{\textrm{\'et}}$) is the disjoint union of $U_i$'s. This implies that $\{\bff_i(w_*):i=1,\dots,d\}$ corresponds to $d$ \emph{distinct} eigensystems (of the Hecke algebra $\mathbf{T}$, see \cite[Lemma 5.9]{buzzard}) for any weight $w_*\in\calD^{\textrm{\'et}}$. Therefore $\{\bff_i(w_*):i=1,\dots,d\}$ is a basis of the $d$-dimensional vector space $\calS^{\new}_{\calD^{\textrm{\'et}}}|_{w=w_*}$ for any $w_*\in\calD^{\textrm{\'et}}$. From this we deduce that $\{\bff_i(w):i=1,\dots,d\}$ is $\calA^{\textrm{\'et}}$-linearly independent, and under this basis we obtain the desired matrix of the $U_p$-action on $\calS^\new_{\calD^{\textrm{\'et}}}$.
\end{proof}
 
 	\begin{remark}
 		\label{rmk: slope function}
 	For a general closed disc $\calD\subseteq D(k)$ centered at $w_k$, the characteristic polynomial of the $U_p$ action on $\calS_{\calD}^\new$ always exists, and one can express $U_p$ as a (not necessarily diagonal) matrix whenever $\calS^{\new}_{\calD}$ is finite free. In a sense,  this characteristic polynomial can be regarded as a \emph{weak} form of extending slope functions $a_{p,i}$'s to a larger disc $\mathcal{D}$; namely one can always extend the \emph{symmetric polynomials} of these functions. 
 	
 	In general, slope functions $a_{p,i}$'s \emph{may not} extend to rigid analytic functions defined on $\calD$, due to the ramifications on the eigencurve; see \cite[\S~5]{const} for more discussions about an concrete example constructed in Emerton's thesis.
 \end{remark}
 
In the rest of this section, we fix a sufficiently small closed disc $\calD^{\textrm{\'et}}\subseteq D(k)$ centered at $w_k$ as in Lemma \ref{lemma: basis}.

	\begin{definition}[Taylor expansions of $U^2_p$]
		\label{def: taylor of Up}
	 Let $M(w)=M_{U^2_p}^\new(w)\in\M_{d\times d}(\calA^{\textrm{\'et}})$ be the matrix of $U^2_p$-action on $\calS^{\new}_{\calD^{\textrm{\'et}}}$ under the basis $\{\mathbf{f}_i(w):i=1,\dots,d\}$. We expand $M(w)$ around $w=w_k$ as $
	 	M(w) = \left(m_{i,j}(w)\right)_{d \times d} \in \M_{d \times d}(\mathcal{A^{\textrm{\'et}}})$.
Write $m_{i,j}(w) =\sum_{\ell\geq0} m_{i,j}^{(\ell)}(w-w_k)^\ell\in\calA^{\textrm{\'et}}$, where each $m_{i,j}^{(\ell)} \in \mathbb{C}_p$.
For $\ell\geq0$, define $A_\ell := \left(m_{i,j}^{(\ell)}\right)_{d \times d} \in \M_{d \times d}(\mathbb{C}_p)$. Then we have
	 \begin{equation}
	 	\label{Taylor expansion}
M(w) =\sum_{\ell\geq0}A_\ell(w-w_k)^\ell, \ \textrm{for}\ w\in\calD^{\textrm{\'et}},
	 \end{equation}
where $A_0$ is the scalar matrix $p^{k-2}I_d$ (as the $U_p$-eigenvalues of $f_i$'s are either $p^{\frac{k-2}{2}}$ or $-p^{\frac{k-2}{2}}$), and each matrix $A_\ell$ is diagonal by Lemma \ref{lemma: basis}.
	\end{definition}

	\subsection{Greenberg--Stevens formula in higher ranks}
	In this subsection, we combine the expansion \eqref{Taylor expansion} with the Greenberg--Stevens method to connect the $\calL$-invariants associated with $f_i$ for $i=1,\dots,d$ to the eigenvalues of the matrix $A_1$.
	
	\begin{lemma}[Greenberg--Stevens]
		\label{lem:Greenberg--stevens}
		For $i\in\left\{1,\dots,d\right\}$, we rewrite the Taylor expansion of the slope function $a_{p,i}(w)$ for $w\in\calD^{\textrm{\'et}}$ as
		\begin{equation}
			\label{eq:expansion of a_{p,i}}
			a_{p,i}(w)=\sgn(i) \cdot p^{\frac{k-2}{2}}+c_i(w-w_k)+\textrm{higher terms in} \ (w-w_k),
		\end{equation}
		where each $c_i$ is a $p$-adic number, and $\sgn(i)\in\{\pm1\}$, depending on the signature of the $U_p$-eigenvalue $a_{p,i}(0)$ of $f_i$.
		Then the $\calL$-invariants associated with $\bar r$-newforms $\{f_i\}_{i=1,\dots,d}$ of weight $k$ are given by
	\begin{equation}
			\left\{\calL_{f_i}=-\sgn(i)\cdot2p^{1-\frac{k-2}{2}}\cdot\exp(p(k-2))\cdot c_i:i=1,\dots,d\right\}.
\end{equation}
		In particular, the slopes of these $\calL$-invariants are given by
\begin{equation}
			\left\{v_p(\calL_{f_i})=v_p(c_i)+1-\frac{k-2}{2}:i=1,\dots,d\right\}.
\end{equation}
	\end{lemma}
	\begin{proof}
		This is a reinterpretation of the Greenberg–Stevens formula in our context. In fact,
		by the Greenberg--Stevens method (see \cite[Lemma 4.1]{const}), we obtain:
		\begin{equation}
			\mathcal{L}_{f_i} = -\frac{2a'_{p,i}(k)}{a_{p,i}(k)} =-\text{sgn}(i) \cdot 2p^{-\frac{k-2}{2}} \cdot a'_{p,i}(k),
		\end{equation}
		where $a_{p,i}(k)=a_{p,i}(w_k)=a_{p,i}(w)|_{w=w_k}$, and $a_{p,i}'(k)$ is the first derivative of the $a_{p,i}$ with respect to $k$ (instead of $w_k$). 
		In our context, $w_k = \exp(p(k-2))-1$ is a function of $k$. Hence
		\begin{equation}
			a'_{p,i}(k) = a'_{p,i}(w_k) \cdot \frac{d}{dk}w_{k} = c_i \cdot p\exp(p(k-2)).
		\end{equation}
		Now the lemma is clear.
	\end{proof}
		\begin{theorem}
		\label{GS formula}
		The $\calL$-invariants associated with $\bar r$-newforms $\{f_i\}_{i=1,\dots,d}$ of weight $k$ are precisely the eigenvalues of the matrix
		\begin{equation}
			\label{equations gs 1}
			-A_1\cdot A_0^{-1}\cdot \left(p\cdot\exp(p(k-2))I_d\right).
		\end{equation}
		In particular, we have the following descriptions for the slopes of these $\calL$-invariants.
		\begin{equation}
			\label{equation gs3}
			\left\{\textrm{slopes of} \ A_1\right\}=\left\{v_p(\calL_f)+k-3:f \ \textrm{is a}\ \bar r\textrm{-newform of weight}\ k
			\right\}.
		\end{equation}
	\end{theorem}
	\begin{proof}
		We maintain the expansions of the slope functions in \eqref{eq:expansion of a_{p,i}}.
	By Lemma \ref{lemma: basis} and \eqref{Taylor expansion},  $M(w)=\diag(a_{p,1}(w),\cdots,a_{p,d}(w))^2$, and hence
	\[
	A_1= \diag\left(\sgn(1)\cdot 2p^{\frac{k-2}{2}}c_1, \dots, \sgn(d)\cdot 2p^{\frac{k-2}{2}}c_d\right).
	\]
Now this theorem follows directly from Lemma \ref{lem:Greenberg--stevens} (note that $\exp(p(k-2))$ is a $p$-adic unit).
	\end{proof}

\begin{remark}
By Remark \ref{rmk: slope function}, the $U_p^2$-action on $\calS_{\calD}^{\new}$ can be expressed as a (not necessarily diagonal) matrix as in \eqref{Taylor expansion}, as long as $\calS_{\calD}^{\new}$ is $\calA$-free. The corresponding ``$A_1$" in the Taylor expansion may not be diagonal in general, but the eigenvalues of $A_1$ coincide with that in Theorem \ref{GS formula}. This can be checked directly by using the fact that $A_0$ is a scalar matrix in both cases. The details are omitted, as they will not be used in this paper. Keep the notations in \eqref{eq:expansion of a_{p,i}}. If we denote the eigenvalues of $A_1$ by $\left\{\theta_i=\sgn(i)\cdot 2p^{\frac{k-2}{2}}\cdot c_i:i=1,\dots,d\right\}$, then the $\calL$-invariants are given by $	\left\{\calL_{f_i}=-\theta_i\cdot p^{3-k}\cdot\exp(p(k-2)):i=1,\dots,d\right\}$.
\end{remark}
	
\begin{remark}
Since $A_0$ is a scalar matrix, calculating the slopes of $A_1$ in \S~\ref{section of slopes for A_1} becomes easier. For this reason, we choose to consider the Taylor expansion of $U_p^2$ instead of $U_p$.
	\end{remark}

\subsection{Calculating the slopes of $A_1$}
\label{subsection two step approach}
In this subsection, we outline the approach for computing the slopes of $A_1$ in \S~\ref{section of slopes for A_1}. To illustrate the method, we study a concrete example (Example \ref{example1}) in \S~\ref{section : p-adic analysis and an example}, before addressing the general case.

We start with several notations.
For the remainder of this paper, we fix a sufficiently small closed disc
$
\Spm\calA^{\textnormal{\'et}}=\calD^{\textnormal{\'et}}\subseteq D(k)
$
centered at $w_k$, such that $\calS_{\calD^{\textnormal{\'et}}}^{\new}$
is a finite free $\calA^{\textnormal{\'et}}$-module of rank $d$. We also
choose a basis in which the $U_p^2$-action on
$\calS_{\calD^{\textnormal{\'et}}}^{\new}$ is given by the
\emph{diagonal} matrix
$
M(w)
=
\sum_{\ell\geq 0} A_\ell (w-w_k)^\ell,$
for $w\in\calD^{\textnormal{\'et}}$, as in Definition~\ref{def: taylor of Up}. We denote the \emph{characteristic polynomial of $M(w)$} by
$C^{\new}(w,t)$, and write
\begin{equation}
	\label{eq16}
	C^{\new}(w,t)
	:=
	\det\bigl(1-tM(w)\bigr)
	=
	1+\sum_{i=1}^d(-1)^i c^{\new}_i(w)t^i.
\end{equation}
By Proposition~\ref{prop:newslope component}, this polynomial is
well-defined for every $w_*\in D(k)$. In particular, for each
$i\in\{1,\dots,d\}$, the rigid analytic function $c_i^{\new}(w)$ is
defined on $D(k)$ and we write its Taylor expansion at $w_k$ as
\begin{equation}
	\label{eq17}
	c_i^{\new}(w)
	=
	\sum_{\ell\geq 0} c_i^{\new,(\ell)}(w-w_k)^\ell .
\end{equation}
For $w_*\in D(k)$, we denote the
\emph{Newton polygon of $C^{\new}(w_*,t)$} by $\NP^{\new}(w_*)$.

We thank the referee for suggesting the following convention of expressing the
coefficients of the characteristic polynomial of a diagonal matrix in
terms of \emph{elementary symmetric polynomials}. This makes the notations and the
subsequent discussions in
\S~\ref{subsection:construction of linear equations} much clearer.
\begin{convention}
	\label{convention:symm [oly]}
	For $d$ variables $\underline{\xi}=(\xi_1,\dots,\xi_d)$, let
	$\mathscr{S}_i(\underline{\xi})$ denote the $i$-th elementary
		symmetric polynomial
		\[
	\mathscr{S}_i(\underline{\xi})=	\scrS_i\left(\xi_1,\dots,\xi_d\right)
		:=
		\sum_{1\leq r_1<\cdots<r_i\leq d}
		\xi_{r_1}\cdots \xi_{r_i}.
		\]
If
	$
	A=\operatorname{diag}(a_{1,1},\dots,a_{d,d})
	$
	is a diagonal $d\times d$ matrix, then we write
	$
	\mathscr{S}_i(\underline{A})
	:=
	\mathscr{S}_i(a_{1,1},\dots,a_{d,d})
	$. Under this convention, the characteristic polynomial of the diagonal matrix $A$ is given by
	\[
	\det(1-tA)=1+\sum_{i=1}^d(-1)^i\scrS_i\left(\underline{A}\right)t^i.
	\]
\end{convention}
\subsection*{Strategy of computing the slopes of $A_1$}
Let 	$
\tau_i:=\scrS_i\bigl(\underline{A_1}\bigr)\in\mathbb{C}_p
$. Computing the slopes of $A_1$ amounts to determining the valuations $\left\{v_p(\tau_i):i=1,\dots,d\right\}$. To achieve this, our key strategy consists of the following two techniques. 
\begin{enumerate}
	\item \textrm{Linear equations:} 
Combining the Taylor expansions \eqref{Taylor expansion} and \eqref{eq17}
with the characteristic polynomial \eqref{eq16}, we obtain certain linear
relations between
\(\{\tau_i\}_{1\leq i\leq d}\) and  \(\left\{c^{\new,(\ell)}_{i}\right\}_{1\leq i,\ell\leq d}\); see Proposition
\ref{Cor for Simply F_i(u)}. Consequently, each \(p^{-i(k-2)}\tau_i\)
can be expressed as a \emph{specific} \(\ZZ\)-linear combination of the
elements
\[
\left\{
p^{-(k-2)}c^{\new,(i)}_{1},
p^{-2(k-2)}c^{\new,(i)}_{2},
\ldots,
p^{-i(k-2)}c^{\new,(i)}_{i}
\right\};
\]
see Lemma \ref{lem:ZZ-linaer comb lemma} and Proposition
\ref{prop: ZZ-linear comb prop}.
	\item \textrm{$p$-adic analysis:} This analytic technique will be explained in detail in \S~\ref{subsection: analytic trick}. Roughly speaking, using the \emph{duality} introduced in \S~\ref{subsection:duality}, we convert lower bounds for the \emph{dual graph} \(\NP^*(c_i^\new)\), obtained from the \((k,w_*)\)-newslopes for various \(w_*\in D(k)\), into lower bounds for the Newton polygon \(\NP(c_i^\new)\). This yields certain estimates for
	\[
	v_p\left(c^{\new,(i)}_{1}\right),\ v_p\left(c^{\new,(i)}_{2}\right),\ \dots,\ v_p\left(c^{\new,(i)}_{i}\right).
	\]
\end{enumerate}
Combining (1) and (2), we obtain estimates for the valuations \(v_p(\tau_i)\); see Proposition~\ref{prop: estimates of X_j's in an example} and Theorem~\ref{prop slopes for A_1}. These estimates determine \emph{almost all} of the slopes of \(A_1\).
\begin{remark}
In (1), the linear relations are obtained using only the Taylor expansions of \(M(w)\) and \(c_i^\new(w)\) for \(w\in\calD^{\textrm{\'et}}\). However, as discussed in Remark~\ref{rmk: analytic trick}, it is crucial in (2) to consider the \emph{variation} of the \((k,w_*)\)-newslopes as \(w_*\) ranges over weights \emph{sufficiently far} from \(w_k\).
The farther the weight \(w_*\) is from \(w_k\), the more information we may obtain about \(\NP(c^\new_i)\). It is therefore necessary to work over a larger region, namely the open disc \(D(k)\), rather than over \(\calD^{\textrm{\'et}}\). Fortunately, this is sufficient to determine \emph{almost all} of the slopes of the \(\calL\)-invariants of interest, up to a logarithmic error term; see Corollary~\ref{cor slopes of L-inv} and Remark~\ref{rmk of log bound meaning}.
\end{remark}

	\section{An analytic technique and the linear equations}
	\label{section : linear equations}
	In this section, we develop the two techniques introduced in \S~\ref{subsection two step approach}. First, we recall the notions of \emph{Gauss norms} and the corresponding \emph{dual graph} associated with a rigid analytic function. The relations between the dual graph and the Newton polygon (Lemma \ref{lem:dual} and Corollary \ref{analytic trick}) serve as our analytic trick for obtaining estimates of $v_p\left(c^{\new,(\ell)}_{i}\right)$'s from the $(k,w_*)$-newslopes for $w_*\in D(k)$; see \S~\ref{subsection: analytic trick}. Next, we construct a system of linear equations relating the entries of $\left(c_{i}^{\new,(\ell)}\right)_{1\leq i,\ell\leq d}$ to $\left\{\tau_j=\scrS_j\left(\underline{A_1}\right)\right\}_{1\leq j\leq d}$; see Proposition \ref{Cor for Simply F_i(u)}. As a consequence, each $\tau_i$ can be expressed as a \emph{specific} $\ZZ$-linear combination of $c_{i}^{\new,(\ell)}$'s (see Proposition \ref{prop: ZZ-linear comb prop}), which plays a key role in estimating the valuations $v_p(\tau_i)$'s (Theorem \ref{prop slopes for A_1}) in \S~\ref{section of slopes for A_1}. 
	
	 Keep the settings in \S~\ref{subsection two step approach} throughout this section.
	 
	 \subsection{Gauss norms and the dual graph}
	 \label{subsection:duality}
	 Let $f(w)$ be a power series convergent on the $k$-dominant disc $D(k)$, or equivalently convergent for $v_p(w-w_k)>\M(k)$.
	 \begin{definition}
	 	\label{def:dual graph}
	 	For $r>\M(k)$, the \emph{$r$-Gauss norm} of $f(w)=a_0+a_1(w-w_k)+\cdots$ is defined to be
	 	\begin{equation}
	 		\nu_r(f):=\min_{n\geq0}\left\{v_p(a_n)+nr\right\}=\min_{v_p(w-w_k)=r}\left\{v_p(f(w))\right\}.
	 	\end{equation}
	 	The \emph{dual graph} of $f$ is defined to be $\NP^*(f):=\left\{(r,\nu_r(f)):r>\M(k)\right\}\subseteq\RR^2$.
	 \end{definition}
	 
	 The dual graph $\NP^*(f)$ and the Newton polygon $\NP(f)$ mutually determine each other. The following lemma is well-known, see \cite[\S~1.5]{FF2018} for a proof.
	 \begin{lemma}
	 	\label{lem:dual}
	 	The dual graph $\NP^*(f)$ is concave. Moreover, for $r>\M(k)$, and two integers $m>n$, the following two statements are equivalent:
	 	\begin{enumerate}
	 		\item 
	 		The point $(r,\nu_r(f))$ is a vertex on the dual graph $\NP^*(f)$, and the slopes of the corresponding two adjacent segments are $m>n$.
	 		\item 
	 		The valuation $-r$ is a slope of the Newton polygon $\NP(f)$ with multiplicity $m-n$, and the corresponding two adjacent vertices are $\left(n,v_p(a_n)\right)$ and $\left(m,v_p(a_m)\right)$.
	 	\end{enumerate}
	 	Moreover, if either of the above condition holds, then the integers $m$ and $n$ are the greatest and the smallest indices respectively, for which $v_p(a_m)+mr=v_p(a_n)+nr=\nu_r(f)$.
	 \end{lemma}
	 
	 \begin{corollary}
	 	\label{analytic trick}
	 	Keep the notations in Lemma \ref{lem:dual}. For $r>\M(k)$, assume that $(r,\nu_r(f))$ is a point on a segment $I$ of the dual graph $\NP^*(f)$ of slope $n$. Then
	 	$v_p(a_n)=\nu_r(f)-nr$.
	 \end{corollary}
	 \begin{proof}
	 	Denote the left (or right) endpoint of this segment $I$ by $(r',\nu_{r'}(f))$, then by Lemma \ref{lem:dual}, it follows that $v_p(a_n)=\nu_{r'}(f)-nr'$. Note that $n=\frac{\nu_r(f)-\nu_{r'}(f)}{r-r'}$, the corollary follows.
	 \end{proof}
\subsection{An analytic technique}
	 \label{subsection: analytic trick}
	 Fix $i\in\left\{1,\dots,d\right\}$. 
In what follows, we outline the analytic trick used to obtain estimates for $v_p\left(c_{i}^{\new,(\ell)}\right)$'s from the $(k,w_*)$-newslopes for various $w_*\in D(k)$. To this end, our goal is to estimate the Newton polygon $\NP(c_i^{\new})$. Recall $c_i^{\new}(w)=\sum_{\ell\geq0}c_{i}^{\new,(\ell)}(w-w_k)^\ell$ converges when $v_p(w-w_k)>\M(k)$. We summarize the trick in the following commutative diagram.
\\
\usetikzlibrary{arrows.meta, positioning}
\tikzset{
	box/.style = {draw, rectangle, rounded corners, minimum height=1.2em, inner sep=6pt},
	arrow/.style = {thick, ->, >=Stealth}
}
\begin{tikzpicture}[node distance=1cm and 2.5cm]
	
	\node[box, align=center] (newslopes) {$(k, w_*)$-newslopes \\   for $w_*\in D(k)$};
	\node[box, align=center, right=of newslopes] (lowerNP) {Lower bounds \\ for $\NP^*(c^\new_i)$};
	\node[box, align=center, right=of lowerNP] (lowerNPdual) {Lower bounds \\for $\operatorname{NP}(c^\new_i)$};
	\node[box, align=center, below=of lowerNPdual] (vpEstimates) {Estimates for \\$v_p\left(c_i^{\new,(\ell)}\right)$'s};
	
	\draw[arrow] (newslopes) -- (lowerNP) node[midway, above] {\small $\NP^\new(w_*)$};
	\draw[arrow] (newslopes) -- (lowerNP) node[midway, below] {(1)};
	\draw[arrow] (lowerNP) -- (lowerNPdual) node[midway, above] {\small duality};
	\draw[arrow] (lowerNP) -- (lowerNPdual) node[midway, below] {\small (2)};
	\draw[arrow] (lowerNPdual) -- (vpEstimates);
	\draw[arrow] (newslopes.south) -- (vpEstimates.west) node[midway, sloped, below] {\small analytic trick};
\end{tikzpicture}
\\

In (1), the \((k,w_*)\)-newslopes for \(w_*\in D(k)\) determine the Newton polygon \(\NP^\new(w_*)\), and hence give lower bounds for the valuations \(v_p(c^\new_i(w_*))\) of the coefficients of \(C^\new(w_*,t)\). As $w_*$ varies in $D(k)$, these lower bounds provide a \emph{lower bound} for the dual graph $\NP^*(c_i^{\new})$; see Lemmas \ref{lemma: lower bd in the example} and \ref{y_i function}. 
In (2), Lemma~\ref{lem:dual} and Corollary~\ref{analytic trick} transfer the previous lower bound for the dual graph \(\NP^*(c_i^\new)\) to a lower bound for \(\NP(c_i^\new)\), and transfer any \emph{common segment} between \(\NP^*(c_i^\new)\) and its lower bound to a corresponding point on \(\NP(c_i^\new)\).
	 From these analysis of $\NP(c_i^\new)$, we obtain the desired estimates of the valuations $v_p\left(c_{i}^{\new,(\ell)}\right)$'s; see
	 Proposition \ref{prop example of estimate of F_{i,j}} and Proposition \ref{main estimate for scrM}.
	 
We illustrate this technique with an explicit example in the proof of Proposition~\ref{prop example of estimate of F_{i,j}}. The general case is addressed in Proposition~\ref{main estimate for scrM}.

	\subsection{Constructions of the linear equations}
\label{subsection:construction of linear equations}
Keep the notations in \S~\ref{subsection two step approach}. Consider the two Taylor expansions
$
M(w)=\sum_{\ell\geq0}A_\ell(w-w_k)^\ell,
$
and
$
c_i^\new(w)=\sum_{\ell\geq0}c_i^{\new,(\ell)}(w-w_k)^\ell.
$
By Convention~\ref{convention:symm [oly]}, they satisfy
\begin{equation}
	\label{eq20}
	c_i^\new(w)=\scrS_i\left(\underline{M(w)}\right), \qquad w\in\calD^{\textrm{\'et}}.
\end{equation}
The goal of this subsection is Proposition~\ref{Cor for Simply F_i(u)}, where we use the relation \eqref{eq20} to construct certain linear equations relating \(\{\tau_i=\scrS_i\left(\underline{A_1}\right)\}_{1\leq i\leq d}\) and \(\left\{c^{\new,(\ell)}_{i}\right\}_{1\leq i,\ell\leq d}\). The key idea is to \emph{isolate} the contribution of the scalar matrix \(A_0=p^{k-2}I_d\) on the right-hand side of \eqref{eq20}. Intuitively, this allows us to separate the contribution
of the scalar matrix $A_0$ from the higher-order terms,
thereby making the contribution of $A_1$ visible in the
Taylor expansion of $c_i^{\new}(w)$.
\begin{notation}
	Put $\a=p^{k-2}$. Consider the diagonal matrix $M(w)^{+}=\sum_{\ell\geq1}A_\ell(w-w_k)^\ell$. Then we have $M(w)=M(w)^{+}+\a I_d$. For any \(i\geq0\), write
	\begin{equation}
		\scrS_i\left(\underline{M(w)^+}\right)
		=
		\sum_{\ell\geq0}m_i^{+,(\ell)}(w-w_k)^\ell\in\calA^{\textrm{\'et}}.
	\end{equation}
	By construction, \(m_i^{+,(\ell)}=0\) for \(i>\ell\), and
	\begin{equation}
		\label{eq22}
	m_i^{+,(i)}=\scrS_i\left(\underline{A_1}\right)=\tau_i.
	\end{equation}
\end{notation}
The following lemma, suggested by the referee, describes the relation between \(\scrS_i\left(\underline{M(w)}\right)\) and \(\scrS_i\left(\underline{M(w)^+}\right)\). It simplifies the proof of the subsequent Proposition~\ref{Cor for Simply F_i(u)}.
\begin{lemma}
	\label{lem:referee's suggestion}
	For $i\geq0$, we have
	$
	\scrS_i\left(\underline{M(w)}\right)=\sum_{j=0}^i\binom{d-j}{i-j}\a^{i-j}
	\cdot\scrS_j\left(\underline{M(w)^+}\right)
	$.
\end{lemma}
\begin{proof}
	It suffices to show that $\scrS_i(\underline{\xi}+\underline{\alpha})
	=
	\sum_{j=0}^i
	\binom{d-j}{i-j}
	\scrS_j(\underline{\xi})\alpha^{i-j}$, where $\underline{\xi}=(\xi_1,\dots,\xi_d)$ and $\underline{\a}=(\a,\a,\dots,\a)$. This follows directly from the generating function for elementary
	symmetric polynomials: $\sum_{i=0}^d \scrS_i(\underline{\xi})t^i=\prod_{j=1}^d (1+\xi_jt)$; see \cite[Chapter I, Section 2]{MacdonaldSF}.
	Indeed, we have
\[
\sum_{i=0}^d \scrS_i(\underline{\xi}+\underline{\alpha})t^i
=
\prod_{j=1}^d (1+(\xi_j+\alpha)t)
=
(1+\alpha t)^d
\prod_{j=1}^d\left(1+\frac{\xi_jt}{1+\alpha t}\right)=\sum_{\ell=0}^d \scrS_\ell(\underline{\xi})t^\ell(1+\alpha t)^{d-\ell},
\]
where the last equality follows from the generating function with \(t\) replaced by \(\frac{t}{1+\alpha t}\).
Comparing the coefficient of \(t^i\), we obtain the required equality.
\end{proof}
\begin{proposition}
	\label{Cor for Simply F_i(u)}
Fix $\ell\in\left\{1,\dots,d\right\}$. We have the following $d$ linear equations, indexed by $i\in\left\{1,\dots,d\right\}$, with $\ell$ unknowns: $\left\{\a^{-1}\cdot m_1^{+,(\ell)},\a^{-2}\cdot m_2^{+,(\ell)},\dots,\a^{-\ell}\cdot m_\ell^{+,(\ell)}\right\}$.
	\begin{equation}
		\label{5.20.1}
		\a^{-i}\cdot c_{i}^{\new,(\ell)}=\sum_{j=1}^\ell\binom{d-j}{d-i}\cdot\left(\a^{-j}\cdot m_j^{+,(\ell)}\right), \ \textrm{for}\ i\in\left\{1,\dots,d\right\}.
	\end{equation}
	More precisely, the $d\times \ell$ matrix that corresponds to the $d$ equations in \eqref{5.20.1} is given by the first $\ell$ columns of the following matrix $D_d$.
	\begin{equation}
		\label{5.21.2}
		D_d:=\begin{pmatrix}
			1&0&0&\cdots&0\\
			\bi{d-1}{d-2}&1&0&\cdots&0\\
			\bi{d-1}{d-3}&\bi{d-2}{d-3}&1&\cdots&0\\
			\vdots&\vdots&\ddots&\vdots\\
			\bi{d-1}{0}&\bi{d-2}{0}&\bi{d-3}{0}&\cdots&\bi{0}{0}\\
		\end{pmatrix}_{d\times d}=\left(\binom{d-s}{d-r}\right)_{1\leq r,s\leq d}.
	\end{equation}
\end{proposition}
\begin{proof}
	By Lemma~\ref{lem:referee's suggestion}, for any $i\in\left\{1,\dots,d\right\}$, we have
	\begin{equation}
		\label{eq21}
	c_i^{\new}(w)=\scrS_i\left(\underline{M(w)}\right)=\sum_{j=0}^i\binom{d-j}{i-j}\a^{i-j}\cdot\scrS_j\left(\underline{M(w)^+}\right).
	\end{equation}
Therefore, \(c_i^{\new,(\ell)}\) is equal to the coefficient of \((w-w_k)^\ell\) on the right-hand side of \eqref{eq21}, namely
$
\sum_{j=0}^i \binom{d-j}{d-i}\alpha^{i-j}m_j^{+,(\ell)}
$.
Since \(\ell>0\), we have \(m_0^{+,(\ell)}=0\). Moreover, \(m_j^{+,(\ell)}=0\) for \(j>\ell\), and \(\binom{d-j}{d-i}=0\) for \(j>i\). Hence the above sum is equal to
$
\sum_{j=1}^\ell\binom{d-j}{d-i}\alpha^{i-j}m_j^{+,(\ell)}
$.
This proves the proposition.
\end{proof}

\subsection{The $\ZZ$-linear relations between $c^{\new,(\ell)}_{i}$'s and $\tau_i$'s} In this subsection, we present two consequences of Proposition \ref{Cor for Simply F_i(u)}. 
The first one is the following lemma, which may be viewed as a partial inverse to the system of linear equations constructed in Proposition~\ref{Cor for Simply F_i(u)}. 
Recall that \(\tau_i=\scrS_i(\underline{A_1})=m_i^{+,(i)}\); see \eqref{eq22}.

\begin{lemma}
	\label{lem:ZZ-linaer comb lemma}
	Write $\a=p^{k-2}$. For each $j\in\left\{1,\dots,d\right\}$, $\a^{-j}\tau_j$ is a $\ZZ$-linear combination of $\left\{\a^{-1}c^{\new,(j)}_{1},\a^{-2}c^{\new,(j)}_{2},\dots,\a^{-j}c_j^{\new,(j)}\right\}$.
\end{lemma}
\begin{proof}
By Proposition \ref{Cor for Simply F_i(u)}, for any $j\in\left\{1,\dots,d\right\}$, we obtain the following $j$ equations.
\begin{equation}
	\left(\begin{matrix}
		\a^{-1}c_{1}^{\new,(j)}\\
	\vdots\\
		\a^{-j}c_{j}^{\new,(j)}
	\end{matrix}\right)=\left(\begin{matrix}
1\\
\binom{d-1}{d-2}&1\\
\binom{d-1}{d-3}&\binom{d-2}{d-3}&1\\
\vdots&\vdots&\vdots&\ddots\\
\binom{d-1}{d-j}&\binom{d-2}{d-j}&\binom{d-3}{d-j}&\dots&\binom{d-j}{d-j}
	\end{matrix}\right)_{j\times j}\cdot\left(\begin{matrix}
	\a^{-1}m^{+,(j)}_1\\
	\vdots\\
	\a^{-j}m^{+,(j)}_j
	\end{matrix}\right),
\end{equation}
where the $j \times j$ matrix lies in $\GL_j(\ZZ)$ and its determinant is equal to $1$. Therefore its inverse is also a matrix in $\GL_j(\ZZ)$ and this implies that   $\a^{-j}\tau_j=\a^{-j}m^{+,(j)}_j$ can be expressed as a $\mathbb{Z}$-linear combination of $\left\{\a^{-1}c^{\new,(j)}_{1},\a^{-2}c^{\new,(j)}_{2},\dots,\a^{-j}c_j^{\new,(j)}\right\}$.
\end{proof}
Next, we modify the $\ZZ$-linear combinations obtained in Lemma \ref{lem:ZZ-linaer comb lemma} for the case when $j$ is even. A key result is that the \emph{$\ZZ$-coefficient} of $\a^{-\frac{j}{2}}c_{\frac{j}{2}}^{\new,(j)}$ in the modified combination is $\pm1$. This fact plays a crucial role in proving the key estimates of $v_p(\tau_j)$'s in Theorem \ref{prop slopes for A_1}.
\begin{definition}
		Let $j\in\left\{1,\dots,d\right\}$ be an even integer. Consider \emph{the first $\frac{j}{2}$ and the last $\frac{j}{2}$} equations of the $d$ linear equations constructed in \eqref{5.20.1}, and denote the corresponding $j \times j$ matrix by $D_d(j)\in\M_{j\times j}(\ZZ)$; see the following $j\times j$ matrix in \eqref{eqeq3}. Clearly, $D_d(j)$ is obtained from the $d\times d$ matrix $D_d$ by deleting the last $(d-j)$ columns and the $(d-j)$ ``middle" rows, namely the $(\frac{j}{2} + 1)$-th row  to the $(d - \frac{j}{2})$-th row. We obtain the following $j$ modified linear equations. 
		\begin{equation}
			\label{eqeq3}
		\left(	\begin{matrix}\a^{-1}c_{1}^{\new,(j)}\\
			\vdots\\
			\a^{-\frac{j}{2}}c^{\new,(j)}_{\frac{j}{2}}\\
			\a^{\frac{j}{2}-d-1}c^{\new,(j)}_{d-\frac{j}{2}+1}\\
			\vdots\\
			\a^{-d}c^{\new,(j)}_{d}\end{matrix}\right)=\left(\begin{matrix}
			1\\
			\binom{d-1}{d-2}&1\\
			\vdots&\vdots&\ddots\\
			\binom{d-1}{d-\frac{j}{2}}&\binom{d-2}{d-\frac{j}{2}}&\dots&1&0&\dots&0\\
			\binom{d-1}{\frac{j}{2}-1}&\binom{d-2}{\frac{j}{2}-1}&\dots&\binom{d-\frac{j}{2}}{\frac{j}{2}-1}&\binom{d-\frac{j}{2}-1}{\frac{j}{2}-1}&\dots&\binom{d-j}{\frac{j}{2}-1}
			\\
			\binom{d-1}{\frac{j}{2}-2}&\binom{d-2}{\frac{j}{2}-2}&\dots&\binom{d-\frac{j}{2}}{\frac{j}{2}-2}&\binom{d-\frac{j}{2}-1}{\frac{j}{2}-2}&\dots&\binom{d-j}{\frac{j}{2}-2}
			\\
			\vdots&\vdots&\ddots&\vdots&\vdots&\ddots&\vdots\\
			\binom{d-1}{0}&\binom{d-2}{0}&\dots&\binom{d-\frac{j}{2}}{0}&\binom{d-\frac{j}{2}-1}{0}&\dots&1\\
			\end{matrix}\right)_{j\times j}
			\cdot\left(\begin{matrix}
			\a^{-1}m_1^{+,(j)}\\
			\vdots\\
			\a^{-j}m_j^{+,(j)}
			\end{matrix}\right).
		\end{equation}
		\end{definition}
\begin{proposition}
	\label{prop: ZZ-linear comb prop}
 	Write $\a=p^{k-2}$. For each $j\in\left\{1,\dots,d\right\}$, $\a^{-j}\tau_j$ can be written as a $\ZZ$-linear combination of
 $
 \left\{
 \a^{-1}c_{1}^{\new,(j)},\dots,
 \a^{-\frac{j}{2}}c^{\new,(j)}_{\frac{j}{2}},
 \a^{\frac{j}{2}-d-1}c_{d-\frac{j}{2}+1}^{\new,(j)},
 \dots,
 \a^{-d}c_{d}^{\new,(j)}
 \right\}$,
 in which the coefficient of
 \(\a^{-\frac{j}{2}}c_{\frac{j}{2}}^{\new,(j)}\)
 is \(\pm1\).
\end{proposition}
To prove this, we introduce the following concepts of \emph{binomial Vandermonde determinants}. 

	\begin{definition}
	For $n$ variables $x_1,\dots,x_n$, let $\V(x_1,\dots,x_n)$ be the \emph{standard Vandermonde determinant} 
	$
	\det\left(\left(x_j^i\right)_{0\leq i\leq n-1,1\leq j\leq n}\right).
	$
	Define the \emph{binomial Vandermonde determinant} associated with $x_1,\dots,x_n$, denoted by $\BV(x_1,\dots,x_n)$, to be the following determinant. 
	\begin{equation}
		\BV(x_1,\dots,x_n):=\det\left(\bi{x_j}{i}_{0\leq i\leq n-1,1\leq j\leq n}\right).
	\end{equation}
\end{definition}
\begin{lemma}
	\label{lem:linear alg}
	\begin{enumerate}
		\item 	We have the relationship $
		\BV(x_1,\dots,x_n)=\frac{1}{0!1!\cdots (n-1)!}\V(x_1,\dots,x_n).
		$
		In particular, for a fixed integer $n_0\geq0$, we have
		\begin{equation}
			\label{5.25.1}
			\BV(n_0+n-1,n_0+n-2,\dots,n_0+1,n_0)=1.
		\end{equation}
		\item 
		For $0\leq m\leq d-1$, consider a matrix obtained by deleting the last $m$ columns and any $m$ \emph{consecutive} rows from $D_d$ (the $j\times j$ matrix in \eqref{eqeq3} is an example). Then the determinant of this matrix is equal to $1$. 
	\end{enumerate}
	
\end{lemma}	
\begin{proof}
	We first prove (1).
	Consider $\BV(x_1,\dots,x_n)$. Using the first row, one can cancel out the constant terms of the $l$-th row for $l\in\left\{2,\dots,n\right\}$; and then one can use the second row to cancel out the linear terms of the $l$-th row for $l\in\left\{3,\dots,n\right\}$. In general, for $l\in\left\{2,\dots,n\right\}$, using the first $(l-1)$ rows, we may transfer the $l$-th row into the following form
	$
	\begin{pmatrix}\frac{x_1^{l-1}}{(l-1)!}& \frac{x_2^{l-1}}{(l-1)!}&\cdots&\frac{x_n^{l-1}}{(l-1)!}
	\end{pmatrix}
	$.
	From this we deduce the first assertion of (1). The second assertion follows from the first one and the formula for the standard Vandermonde determinants.
	\[
	\BV(n_0+n-1,n_0+n-2,\dots,n_0+1,n_0)=\frac{1}{0!1!\cdots (n-1)!}\cdot\prod_{n_0\leq i<j\leq n_0+n-1}j-i=1.
	\] 

Now we prove (2). Consider a matrix described in (2), which we simply denote by $D_{d,m}$ (We refer to the $j\times j$ matrix in \eqref{eqeq3} as a typical example). Let $m_1\geq 0$ be the index, such that the $(m_1+1)$-th to the $(m_1+m)$-th row of $D_d$ has been deleted. Clearly, $D_{d,m}$ is necessarily a block-upper triangular-matrix of the following form 
	\[
	\begin{pmatrix}
		A_{m_1\times m_1}&0\\
		B_{m_2\times m_1}&C_{m_2\times m_2}
	\end{pmatrix},
	\]
	where $m_1+m_2=d-m$. It is sufficient to show that $\det(A_{m_1\times m_1})=\det(C_{m_2\times m_2})=1$.
	
	Since $D_d$ itself is upper-triangular, the matrix $A_{m_1\times m_1}$ is also upper-triangular and all its diagonal entries are equal to $1$. Therefore $\det(A_{m_1\times m_1})=1$. For $\det(C_{m_2\times m_2})$, we claim it takes the form of some $\BV(m_2-1+?,m_2-2+?,\dots,?)$, which is equal to $1$ by \eqref{5.25.1}. In fact, by \eqref{5.21.2}, $D_d=\left(\binom{d-j}{d-i}\right)_{1\leq i,j\leq d}$, and hence the matrix $C_{m_2\times m_2}$ takes the desired form when we interpret zeros in $C_{m_2\times m_2}$ as appropriate binomial coefficients $\binom{r}{s}$ for $r<s$. So $\det(C_{m_2\times m_2})=1$ and we are done.
\end{proof}

\begin{proof}[Proof of Proposition \ref{prop: ZZ-linear comb prop}]
By Lemma \ref{lem:linear alg}(2) and the construction of $D_d(j)$, we deduce that $\det(D_d(j))=1$. Hence, applying \(D_d(j)^{-1}\in\GL_j(\ZZ)\) to \eqref{eqeq3}, we obtain a specific \(\ZZ\)-linear combination expressing
$
\a^{-j}\tau_j=\a^{-j}m_j^{+,(j)}
$
in terms of the elements
\[
\left\{
\a^{-1}c_{1}^{\new,(j)},\dots,
\a^{-\frac{j}{2}}c^{\new,(j)}_{\frac{j}{2}},
\a^{\frac{j}{2}-d-1}c_{d-\frac{j}{2}+1}^{\new,(j)},
\dots,
\a^{-d}c_{d}^{\new,(j)}
\right\}.
\]
It remains to check that, the coefficient of
$
\a^{-\frac{j}{2}}c_{\frac{j}{2}}^{\new,(j)}
$ in this \(\ZZ\)-linear combination
is \(\pm1\). By Cramer's rule, this coefficient is equal, up to sign, to the \((j-1)\times(j-1)\) minor of \(D_d(j)\) obtained by deleting the row and column corresponding to the entry with index \(\left(\frac{j}{2},j\right)\). By Lemma~\ref{lem:linear alg}(2), this particular minor is equal to \(1\). This proves the proposition.
	\end{proof}

\section{Further analysis of Example \ref{example1}}
\label{section : p-adic analysis and an example}

In this section, we revisit Example \ref{example1} and apply the analytic technique of \S~\ref{analytic trick} to estimate the valuations $v_p\left(c^{\new,(\ell)}_{i}\right)$ under the assumption that $m(\bar r)=1$; see Proposition \ref{prop example of estimate of F_{i,j}}. Combining these estimates with the linear relations in Lemma \ref{lem:ZZ-linaer comb lemma} and Proposition \ref{prop: ZZ-linear comb prop}, we obtain corresponding estimates for the valuations $v_p\left(\tau_i\right)$; see Proposition \ref{prop: estimates of X_j's in an example}. As a consequence, we determine two vertices of the Newton polygon of $A_1$, and hence determine four slopes of $A_1$. The general case is treated in \S~\ref{section of slopes for A_1}, using essentially the same approach.

Keep the notations in \S~\ref{subsection two step approach} throughout this section.
\subsection{Recollection of Example \ref{example1}}
In this case, $k=24$, $\M(k)=2$, and the $k$-derivative slopes are given by $s_3=9$, $s_2=6$, and $s_1=2$. For simplicity, we assume the global multiplicity $m(\bar r)=1$, and hence $d=m(\bar r)d_k^{\new}=6$. Let $w_*\in D(k)$, namely $v_p(w_*-w_k)>\M(k)=2$.
Set $r=v_{p}(w_{*}-w_{24})$. 
We copy the table of  $(k,w_*)$-newslopes for various values of $r>\M(k)=2$, from Example \ref{example1}. 
\begin{center}\renewcommand{\arraystretch}{1.2}
	\begin{tabular}{|c|c|c|c|c|c|c|}
		\hline
		$v_p(w_*-w_k)$&\multicolumn{6}{|c|}{$(k,w_*)$-newslopes}
		\\
		\hline
		$r\geq9$&$11$&$11$&$11$&$11$&$11$&$11$\\
		\hline
		$6\leq r\leq9$&$11-(s_3-r)$&$11$&$11$&$11$&$11$&$11+(s_3-r)$\\
		\hline
		$2<r\leq6$&$8-(s_2-r)$&$11-(s_2-r)$&$11$&$11$&$11+(s_2-r)$&$14+(s_2-r)$\\
		\hline
	\end{tabular}
\end{center}

\subsection{The lower bounds}

Fix $w_*\in D(k)$. Recall that $\NP^\new(w_*)$ denotes the Newton polygon of the characteristic polynomial
$
C^\new(w_*,t)=1+\sum_{i=1}^d(-1)^i c_i^{\new}(w_*)t^i
$
for the action of $U_p^2$. Since the $(k,w_*)$-newslopes are obtained from the action of $U_p$ by Proposition \ref{prop:newslope component}, it follows that, for each $i\in\{1,\dots,d\}$, the $i$-th slope of $\NP^\new(w_*)$ is one half of the $i$-th $(k,w_*)$-newslope.
Therefore, from the table above, we obtain a table of \emph{lower bounds} for the valuations $\frac{1}{2}v_p\left(c_i^{\new}(w_*)\right)$'s for various values $r>2$. 

\begin{example}
	\label{example 3}
In the following table, each valuation in the last six columns is equal to $x/2$, where $(i,x)$ is the corresponding point on the Newton polygon $\NP^\new(w_*)$. An entry marked with ``$\geq$" indicates that the corresponding point $(i,v_p(c_i^\new(w_*)))$ is \emph{not} a vertex of $\NP^\new(w_*)$, and hence we obtain only a lower bound for $\frac{1}{2}v_p(c_i^\new(w_*))$. The remaining entries (i.e. without ``$\geq$") correspond to vertices of $\NP^\new(w_*)$, and therefore we obtain the exact values of $\frac{1}{2}v_p(c_i^\new(w_*))$'s.
\begin{center}\renewcommand{\arraystretch}{1.2}
	\tiny
	\begin{tabular}{|c|c|c|c|c|c|c|}
		\hline
		$v_p(w_*-w_k)$&$\frac{1}{2}v_p(c_1^\new(w_*))$&$\frac{1}{2}v_p(c^\new_2(w_*))$&$\frac{1}{2}v_p(c^\new_3(w_*))$&$\frac{1}{2}v_p(c^\new_4(w_*))$&$\frac{1}{2}v_p(c^\new_5(w_*))$&$\frac{1}{2}v_p(c^\new_6(w_*))$
		\\
		\hline
		$r\geq9$&$\geq11$&$\geq22$&$\geq33$&$\geq44$&$\geq55$&$66$\\
		\hline
		$6\leq r\leq9$&$11-(s_3-r)$&$\geq22-(s_3-r)$&$\geq33-(s_3-r)$&$\geq44-(s_3-r)$&$55-(s_3-r)$&$66$\\
		\hline
		$2<r\leq6$&$8-(s_2-r)$&$19-2(s_2-r)$&$\geq30-2(s_2-r)$&$41-2(s_2-r)$&$52-(s_2-r)$&$66$\\
		\hline
	\end{tabular}
\end{center}
\end{example}

We summarize the lower bounds for the dual graphs  $\NP^*(c_i^\new)$'s, obtained from the table above, in the following lemma.
	\begin{lemma}
		\label{lemma: lower bd in the example}
Consider the following piecewise linear functions $f_i(r)$'s for $r>2$.
		\begin{equation}
			\begin{aligned}
				&f_1(r)=\left\{\begin{array}{ll}22&\mbox{if $r\geq9,$}\\4+2r&\mbox{if $2<r\leq 9$.}\end{array}\right.\ \ 
	f_5(r)=\left\{\begin{array}{ll}110&\mbox{if $r\geq9,$}\\92+2r&\mbox{if $2<r\leq 9$.}\end{array}\right.\ \
	f_6(r)=132.
	\\&
	f_2(r) = \left\{
				\begin{array}{lll}
					44 & \text{if } r \geq 9, \\
					26 + 2r & \text{if } 6 \leq r \leq 9, \\
					14 + 4r & \text{if } 2 < r \leq 6.
				\end{array}
				\right.\ \ 
				f_3(r) = \left\{
				\begin{array}{lll}
					66 & \text{if } r \geq 9, \\
					48 + 2r & \text{if } 6 \leq r \leq 9, \\
					36 + 4r & \text{if } 2 < r \leq 6.
				\end{array}
				\right.\ \ 
				\\&					
					f_4(r) = \left\{
				\begin{array}{lll}
					88 & \text{if } r \geq 9, \\
					70 + 2r & \text{if } 6 \leq r \leq 9, \\
					58+ 4r & \text{if } 2 < r \leq 6.
				\end{array}
				\right.
			\end{aligned}
		\end{equation}
Let $\boldsymbol{\Gamma}(f_i)$ be the graph of the function $f_i(r)$ in $\RR^2$.
Then for any $i\in\left\{1,\dots,6\right\}$, the dual graph $\NP^*(c_i^\new)$ lies above $\boldsymbol{\Gamma}(f_i)$ in $\RR^2$. Moreover, two graphs coincide when $i=1$ (resp. $i=2$, $i=4$, $i=5$, and $i=6$) and $2<r\leq9$ (resp. $2<r\leq6$, $2<r\leq6$, $2<r\leq9$, and $r>2$). 
	\end{lemma} 
	\begin{proof}
	Note that the valuation $v_p(c_i^\new(w_*))$ under consideration depends only on the valuation $r=v_p(w_*-w_k)$, the lemma follows directly from the table in Example \ref{example 3}.
	\end{proof}
	\subsection{Estimates of $v_p\left(c_{i}^{\new,(\ell)}\right)$'s}
	Recall the Taylor expansion $c_i^\new(w)=\sum_{\ell\geq0}c_i^{\new,(\ell)}(w-w_k)^\ell$, where $w\in D(k)$. The following example treats the case $i=1$ and illustrates both the proof strategy and the form of the estimates that appear in Proposition~\ref{prop example of estimate of F_{i,j}}.
\begin{example}
	\label{subex 1}
In this example, we study the coefficients of $c_1^{\new}$. Recall that the
$k$-derivative slopes are
$
s_3=9, s_2=6, s_1=2$.
Our goal is to prove the following estimates for the coefficients
$
\left\{c_1^{\new,(1)},\dots,c_1^{\new,(6)}\right\}$:
\begin{equation}
	\label{eq24}
	\begin{aligned}
		v_p\left(c_1^{\new,(1)}\right)&\geq 22-s_3, &
		v_p\left(c_1^{\new,(2)}\right)&= 22-2s_3,\\
		v_p\left(c_1^{\new,(3)}\right)&> 22-2s_3-s_2, &
		v_p\left(c_1^{\new,(4)}\right)&> 22-2s_3-2s_2,\\
		v_p\left(c_1^{\new,(5)}\right)&> 22-2s_3-3s_2, &
		v_p\left(c_1^{\new,(6)}\right)&> 22-2s_3-4s_2.
	\end{aligned}
\end{equation}
By Lemma~\ref{lemma: lower bd in the example}, the dual graph $\NP^*(c_1^\new)$ lies above the graph $\boldsymbol{\Gamma}(f_1)$, which is shown in the following picture. The \emph{solid segment} indicates the part where $\boldsymbol{\Gamma}(f_1)$ coincides with $\NP^*(c_1^\new)$, while the dashed segment indicates the part where we only know that $\NP^*(c_1^\new)$ lies above $\boldsymbol{\Gamma}(f_1)$.

\begin{center}
	\begin{tikzpicture}[xscale=0.5,yscale=0.15]
		\draw[->] (1.5,0) -- (12.5,0) node[right] {$r$};
		\draw[->] (1.5,6) -- (1.5,25) node[above] {$f_1(r)$};
		
		\draw (2,0.2) -- (2,-0.2) node[below] {$s_1=2$};
		\draw (6,0.2) -- (6,-0.2) node[below] {$s_2=6$};
		\draw (9,0.2) -- (9,-0.2) node[below] {$s_3=9$};
		
		\draw (1.7,8) -- (1.3,8);
		\draw (1.7,22) -- (1.3,22) node[left] {$22$};
		
		\draw[thick] (2,8) -- (9,22);
		
		\draw[thick,dashed] (9,22) -- (12,22);
		
		\draw[thick,fill=white] (2,8) circle (3pt);
		
		\fill (9,22) circle (2.5pt);
		
		\node[below right] at (2,8) {};
	\end{tikzpicture}
\end{center}
Since $\NP^*(c_1^\new)$ lies above $\boldsymbol{\Gamma}(f_1)$, we have, for every $r>2$,
$
v_p\left(c_1^{\new,(j)}\right)+jr\geq \nu_r(c_1^\new)\geq f_1(r)$.
Hence
\begin{equation}
	\label{eq25}
	v_p\left(c_1^{\new,(j)}\right)\geq f_1(r)-jr.
\end{equation}
In what follows, we prove \eqref{eq24} using \eqref{eq25} and Corollary~\ref{analytic trick}. Consider the solid segment of slope $2$ on $\boldsymbol{\Gamma}(f_1)$.
First, taking $r=s_3$ and $j=1$ in \eqref{eq25}, we obtain
	$
	v_p\left(c_1^{\new,(1)}\right)\geq 22-s_3,
	$
	which proves the first inequality in \eqref{eq24}.
Next, by Corollary~\ref{analytic trick}, the vertex $(s_3,22)$ on this segment gives
	$
	v_p\left(c_1^{\new,(2)}\right)=22-2s_3$,
	which proves the equality in \eqref{eq24}. It remains to deal with the strict inequalities in \eqref{eq24}. Since $2<r\leq s_3$, we have
	$
	f_1(r)=4+2r=22-2s_3+2r$.
	Thus \eqref{eq25} gives
	\[
	v_p\left(c_1^{\new,(j)}\right)
	\geq 22-2s_3-(j-2)r.
	\]
	Taking any $r$ with $2<r<s_2$, we obtain the remaining strict inequalities in \eqref{eq24}.

\end{example}
We need to introduce the following notation in order to state the general case.
\begin{notation}
	\label{notation : in the example}
For $i,j\in\left\{1,\dots,6\right\}$, let $\boldsymbol{k}_i:=i(k-2)=22i$ and $\boldsymbol{l}_j$ be the following values:
	\[\boldsymbol{l}_j:=\left\{\begin{array}{llllll}
			s_3=9 & \text{if } j=1, \\
		2s_3=18 & \text{if } j=2, \\
		 2s_3+s_2=24& \text{if }j=3,\\
		 	2s_3+2s_2=30 & \text{if } j=4, \\
		 2s_3+3s_2=36 & \text{if } j=5, \\
		 2s_3+4s_2=42& \text{if }j=6.
	\end{array}
	\right.
	\]
Consider the following $6\times 6$ table $\boldsymbol{M}=(\boldsymbol{m}_{i,j})_{6\times 6}$, whose entries $\boldsymbol{m}_{i,j}\in\{>,=,\geq\}$. 
		\begin{center}
		\begin{tabular}{|c|c|c|c|c|c|c|}
			\hline
			$\boldsymbol{m}_{i,j}$&$j=1$&$j=2$&$j=3$&$j=4$&$j=5$&$j=6$\\
			\hline
			$i=1$&$\geq$&=&$>$&$>$&$>$&$>$\\
			\hline
			$i=2$&$\geq$&$\geq$&$\geq$&=&$>$&$>$\\
			\hline
			$i=3$&$\geq$&$\geq$&$\geq$&$\geq$&$\geq$&$\geq$\\
			\hline
			$i=4$&$\geq$&$\geq$&$\geq$&$=$&$>$&$>$\\
			\hline
			$i=5$&$\geq$&$=$&$>$&$>$&$>$&$>$\\
			\hline
			$i=6$&$>$&$>$&$>$&$>$&$>$&$>$\\
			\hline
		\end{tabular}
	\end{center}
We use $v_p\left(c_{i}^{\new,(j)}\right)\approx^{\boldsymbol{m}_{i,j}}(\boldsymbol{k}_{i}-\boldsymbol{l}_{j})$ to denote the relation between $v_p\left(c_{i}^{\new,(j)}\right)$ and $\boldsymbol{k}_i-\boldsymbol{l}_j$ is given by the entry $\boldsymbol{m}_{i,j}\in\{>,=,\geq\}$.
\end{notation}
\begin{remark}
	\label{rmk1}
The table $\boldsymbol{M}$ is symmetric in the sense that the first and second rows coincide with the fifth and fourth rows, respectively. Moreover, the pattern of each row of $\boldsymbol{M}$, which we have already seen for the first row in Example~\ref{subex 1}, can be summarized as follows.

As in Example~\ref{subex 1}, we call the common segment of $\NP^*(c_i^\new)$ and $\boldsymbol{\Gamma}(f_i)$ the \emph{solid segment}. If $\boldsymbol{\Gamma}(f_i)$ contains no solid segment, then all entries $\boldsymbol{m}_{i,j}$ are non-strict inequalities. Otherwise, a solid segment of slope $n$ gives rise to exactly $d-n$
strict inequalities $\boldsymbol{m}_{i,j}$ with $j>n$. If $n\neq0$, it
also gives rise to one equality $\boldsymbol{m}_{i,n}$ and $n-1$
non-strict inequalities $\boldsymbol{m}_{i,j}$ with $j<n$.
\end{remark}

	\begin{proposition}
		\label{prop example of estimate of F_{i,j}}
		Keep the notations in Notation \ref{notation : in the example}. For $1\leq i,j\leq 6$, we have
		\[
		v_p\left(c_{i}^{\new,(j)}\right)\approx^{\boldsymbol{m}_{i,j}}(\boldsymbol{k}_i-\boldsymbol{l}_j).
		\]
	\end{proposition}
	\begin{proof}
	Since $\boldsymbol{M}$ is symmetric in the sense that the first (resp. second) row coincides with the fifth (resp. fourth) row, it suffices to deal with the cases for $i=1,2,3,6$.
		By Lemma \ref{lemma: lower bd in the example}, the graph $\boldsymbol{\Gamma}(f_i)$ lies below the dual graph $\NP^*(c_i^\new)$ for $r>2$. Consequently, we have
		\begin{equation}
			\label{equ1}
		v_p\left(c_i^{\new,(j)}\right)+jr\geq \nu_r\left(c^\new_i\right)\geq f_i(r),\  \textrm{for any}\  j\geq0, \ \textrm{and}\ r>2.
		\end{equation}
We verify the proposition case by case as follows. The argument follows the strategy used in Example~\ref{subex 1}. For the strict and non-strict inequalities in $\boldsymbol{M}$, we choose suitable values of $r$ so that the relation between $f_i(r)-jr$ and $\boldsymbol{k}_i-\boldsymbol{l}_j$ is given by $\boldsymbol{m}_{i,j}$. For the equalities in $\boldsymbol{M}$, we apply Corollary~\ref{analytic trick} directly.

Consider $i=6$. For any $2<r<s_2=6$ and any $j\in\left\{1,\dots,6\right\}$, we have $jr<\boldsymbol{l}_j$. Therefore by \eqref{equ1}, we obtain
\[
v_p\left(c_{i}^{\new,(j)}\right)\geq f_i(r)-jr=132-jr>132-\boldsymbol{l}_j=(k-2)i-\boldsymbol{l}_j=\boldsymbol{k}_i-\boldsymbol{l}_j.
\]

Now we focus on the cases for $i=1,2,3$. By Lemma \ref{lemma: lower bd in the example}, $f_i(r)-jr=\boldsymbol{k}_i-\boldsymbol{l}_j$ holds for the triples $(i,j,r)$'s listed in the following table.
\begin{center}
	\begin{tabular}{|c|c|c|c|c|c|c|}
		\hline
		&$j=1$&$j=2$&$j=3$&$j=4$&$j=5$&$j=6$\\
		\hline
		$i=1$&$r=9$&&&&&
		\\
		\hline
		$i=2$&$r=9$&$6\leq r\leq 9$&$r=6$&&&	\\
		\hline
		$i=3$&$r=9$&$6\leq r\leq 9$&$r=6$&$2< r\leq 6$&$r=6$&$r=6$	\\
		\hline
	\end{tabular}
\end{center}
Combining this with \eqref{equ1}, we obtain the desired non-strict inequalities in $\boldsymbol{M}$. For the strict inequalities in $\boldsymbol{M}$, we need to consider $j=2$ (resp. $j=4$) and $i=\frac{j}{2}$.
In this case, consider $r_0=9$ (resp. $r_0=6$), then from the table above, we have $f_i(r_0)-jr_0=\boldsymbol{k}_i-\boldsymbol{l}_j$. Now for $j'>j$ and some $2<r<s_2=6$, by \eqref{equ1}, we have 
\begin{equation}
	\begin{aligned}
		v_p\left(c_{i}^{\new,(j')}\right)&\geq f_i(r)-j'r\\
		&=(f_i(r_0)-jr_0)+(f_i(r)-f_i(r_0)+jr_0-jr)+(j-j')r\\
		&=\boldsymbol{k}_i-\boldsymbol{l}_j+0+(j-j')r>\boldsymbol{k}_i-\boldsymbol{l}_{j'},
	\end{aligned}
\end{equation}
where we use the fact that $f_i(r)$ is a linear function of slope $j$ over the interval $[r,r_0]$, and that $\boldsymbol{l}_{j'}>\boldsymbol{l}_j+(j'-j)r$.
Finally, we handle the equalities in $\boldsymbol{M}$, and it suffices to consider $j\in\{2,4\}$ and $i=\frac{j}{2}$. Consider $2< r\leq9$ (resp. $2< r\leq 6$), then by Lemma \ref{lemma: lower bd in the example}, $(r,\nu_{r}(c_1^\new))=(r,f_1(r))$ (resp. $(r,\nu_{r}(c^\new_2))=(r,f_2(r))$) is a point on some segment (this is the solid segment mentioned in Remark~\ref{rmk1}) of slope $2$ (resp. $4$) on the dual graph $\NP^*(c_1^\new)$ (resp. $\NP^*(c^\new_2)$). In this case, by Corollary \ref{analytic trick}, we have
\begin{equation}
	\label{equ2}
\begin{aligned}
	v_p\left(c^{\new,(j)}_{i}\right)=\nu_{r}(c^\new_i)-jr=f_i(r)-jr=\boldsymbol{k}_i-\boldsymbol{l}_j,
\end{aligned} 
\end{equation} where the last equality is included in the above table.
The proposition is proved.
	\end{proof}
	
	\subsection{Estimates of $v_p(\tau_i)$'s} Recall $\a=p^{k-2}=p^{22}$. Combining Proposition \ref{prop example of estimate of F_{i,j}} and the results in \S~\ref{section : linear equations}, we obtain the following estimates of $v_p(\tau_i)$'s.
	\begin{proposition}
		\label{prop: estimates of X_j's in an example}
		For $1\leq j\leq 6$, we have  $v_p(\tau_j)\geq\boldsymbol{k}_j-\boldsymbol{l}_j= j(k-2)-\boldsymbol{l}_j$. Moreover, for $j=2$ or $j=4$, the above inequality is an equality. As a consequence, the first four slopes of the matrix $A_1$ in this special case are given by 
		$
		\{k-2-s_3,k-2-s_3,k-2-s_2,k-2-s_2\}=\{13,13,16,16\}
		$, and the remaining two slopes are greater than or equal to $k-2-s_1=20$.
	\end{proposition}
	\begin{proof}
	By Lemma~\ref{lem:ZZ-linaer comb lemma}, each $\alpha^{-j}\tau_j$ is a $\ZZ$-linear combination of
	$
	\left\{\alpha^{-s}c_s^{\new,(j)}:s=1,\dots,j\right\}$.
	By Proposition~\ref{prop example of estimate of F_{i,j}}, each term $\alpha^{-s}c_s^{\new,(j)}$ with $1\leq s\leq j$ satisfies
	\[
	v_p\left(\alpha^{-s}c_s^{\new,(j)}\right)
	= -s(k-2)+v_p\left(c_s^{\new,(j)}\right)
	\geq -s(k-2)+\boldsymbol{k}_s-\boldsymbol{l}_j
	= -\boldsymbol{l}_j.
	\]
	It follows that
$
	v_p\left(\alpha^{-j}\tau_j\right)\geq -\boldsymbol{l}_j$, and hence
$
	v_p(\tau_j)\geq j(k-2)-\boldsymbol{l}_j$.
		
For $j=2$ and $j=4$, we prove equality. By Proposition~\ref{prop: ZZ-linear comb prop}, $\a^{-2}\tau_2$ can be written as a $\ZZ$-linear combination of
$
\left\{\a^{-t}c^{\new,(2)}_{t}:t=1,6\right\}
$,
with the coefficient of $\a^{-1}c^{\new,(2)}_{1}$ equal to $\pm 1$. Similarly, $\a^{-4}\tau_4$ can be written as a $\ZZ$-linear combination of
$
\left\{\a^{-t}c^{\new,(4)}_{t}:t=1,2,5,6\right\}
$,
with the coefficient of $\a^{-2}c^{\new,(4)}_{2}$ equal to $\pm 1$.
By Proposition~\ref{prop example of estimate of F_{i,j}}, we have
\[
v_p\left(\a^{-1}c^{\new,(2)}_{1}\right)=-\boldsymbol{l}_2,
\qquad
v_p\left(\a^{-6}c^{\new,(2)}_{6}\right)>-\boldsymbol{l}_2,
\]
and
\[
v_p\left(\a^{-2}c^{\new,(4)}_{2}\right)=-\boldsymbol{l}_4,
\qquad
v_p\left(\a^{-1}c^{\new,(4)}_{1}\right),
v_p\left(\a^{-5}c^{\new,(4)}_{5}\right),
v_p\left(\a^{-6}c^{\new,(4)}_{6}\right)>-\boldsymbol{l}_4.
\]
It follows that
$
v_p\left(\a^{-j}\tau_j\right)=-\boldsymbol{l}_j
$
for $j=2,4$. Equivalently,
$
v_p(\tau_j)=j(k-2)-\boldsymbol{l}_j
$
for $j=2,4$. This gives the desired description of the slopes of $A_1$ in this case.
	\end{proof}
\begin{remark}
By Theorem \ref{GS formula}, the first four slopes of $\calL$-invariants in this case are given by $\{1-s_3,1-s_3,1-s_2,1-s_2\}=\{-8,-8,-5,-5\}$. The remaining slopes are greater than or equal to $1-s_1=-1$.
\end{remark}

		\section{Slopes of the matrix $A_1$}
		\label{section of slopes for A_1}
	In this section, we implement the method outlined in \S~\ref{subsection two step approach} in the general case. As a consequence, we determine \emph{almost all} of the slopes of $A_1$, and provide an explicit bound for the remaining slopes; see Theorem \ref{prop slopes for A_1}. As an application, we deduce the main results concerning the slopes of $\calL$-invariants associated with $\bar r$-newforms of weight $k$; see Corollaries \ref{cor slopes of L-inv}, \ref{dis of L-inv} and \ref{cor: integrality}.
	
Keep the notations for $\bar r$, $m(\bar r)$, and $k_0$ in \S~\ref{subsection of Settings}. Fix an integer $k\geq 2$, such that
		$
		k\equiv k_{0}\pmod{p-1}$. Recall $D(k)$ is the $k$-dominant disc with radius $p^{-\M(k)}$ and $d:=m(\bar r)d_k^\new=\dim S_k(\Gamma_0(Np))^{\new}_{\bar r}$. We now introduce two key notations that will be used to state the main results
		of this section. The first provides a more detailed description of the
		$k$-derivative slopes, while the second generalizes
		Notation~\ref{notation : in the example}.
\begin{notation}[$k$-derivative slopes]
	\label{notation: refined k-derivative slopes}
		For an integer $k$, let $s_1<s_2<\cdots<s_N$ be the distinct $k$-derivative slopes of the $k$-derivative polygon $\underline{\Delta}_k^+$ in Notation \ref{notation of Delta slope}. Recall in Corollary \ref{cor of thresholds}, $M$ is defined to be the smallest index in the range $i\in\left\{1,\dots,N\right\}$ for which $s_i> \M(k)$. In particular, we have
		$$
		s_1<s_2<\cdots<s_{M-1}\leq \M(k)<s_{M}<\cdots<s_N.
		$$
Let $d_i$ be the integer such that the slope $s_i$ occurs in $\underline{\Delta}_{k}^+$ with multiplicity $d_i/m(\bar r)$. Clearly, we have $
		d=m(\bar r)\cdot d_k^{\new}=2\cdot\sum_{i=1}^Nd_i$. Note that $M$, $N$, $d$, and each $s_i$, $d_i$ depend on $k$, we omit the subscript for simplicity.
\end{notation}
\begin{notation}[Definitions of $\boldsymbol{l}_j$'s]
		\label{notation of Li}
 Set $\boldsymbol{l}_0=0$. For $j\in\left\{1,\dots,d\right\}$, let $\boldsymbol{l}_j$ be the sum of the $j$ largest valuations in the multiset of $k$-derivative slopes (where we replace the $k$-derivative slopes $s_1,\dots,s_{M-1}$ by $s_M$)
 $$
 \left\{[s_M;2d_1+\cdots 2d_{M-1}],[s_M;2d_M],\dots,[s_N;2d_N]\right\},
 $$ where $[m;n]$ indicates that the slope $m$ occurs with multiplicity $n$. See Notation \ref{notation : in the example} for an explicit example of these $\boldsymbol{l}_j$'s.
\end{notation}

\subsection{Statements of the main results}
\label{subsection:main result}
In this subsection, we state the key result (Theorem \ref{prop slopes for A_1}) of this section, which is a generalization of Proposition \ref{prop: estimates of X_j's in an example}, providing estimates for $v_p(\tau_j)$'s. After this, we present its consequences concerning the slopes of $\calL$-invariants. The proof of Theorem \ref{prop slopes for A_1} appears in \S~\ref{subsection : proof of main theorem}.
	\begin{theorem}
		\label{prop slopes for A_1}
		Keep the notations in Notations \ref{notation: refined k-derivative slopes} and \ref{notation of Li}. For $j\in\left\{1,\dots,d\right\}$, we have
		\[
		v_p\left(\scrS_j\left(\underline{A_1}\right)\right) = v_p(\tau_j) \geq j(k - 2) - \boldsymbol{l}_j.
		\]
		Furthermore, for \(j \in\{ 2d_N, 2(d_N+d_{N-1}), \dots, 2(d_N+\cdots+d_M)\}\), this inequality is an equality.
	\end{theorem}
	\begin{corollary}[Slopes of $\calL$-invariants]
		\label{cor slopes of L-inv}
		Keep the notations in Notation \ref{notation: refined k-derivative slopes} and \ref{notation of Li}.
		We have the following descriptions of the slopes of $A_1$ and the slopes of $\calL$-invariants associated with $\bar r$-newforms of weight $k$ (in terms of the $k$-derivative slopes of $\underline{\Delta}^+_k$). We rearrange the slopes in ascending order. 
		\label{cor main1}
		\begin{enumerate}
			\item The \emph{first} $2d_{N} + \cdots + 2d_M$ slopes of $A_1$ are given by $$\left\{[k-2-s_N;2d_N],[k-2-s_{N-1};2d_{N-1}],\dots[k-2-s_{M};2d_M]\right\}.$$
			Furthermore, the remaining slopes are greater than or equal to $k-2-\M(k)$.
			\item 
			The \emph{first} $2d_{N} + \cdots + 2d_M$ slopes of the $\calL$-invariants are given by  $$\left\{[1-s_N;2d_N],[1-s_{N-1};2d_{N-1}],\dots,[1-s_{M};2d_M]\right\}.$$ 
Furthermore, the remaining slopes are greater than or equal to $1-\M(k)$.
		\end{enumerate}
		
	\end{corollary}
	\begin{proof}
(2) follows from (1) and Theorem \ref{GS formula}. (1) is a direct consequence of Theorem \ref{prop slopes for A_1}. Indeed, the estimates for $j\in\{1,\dots,d\}$ in
Theorem~\ref{prop slopes for A_1} give a lower bound for the Newton polygon of
$A_1$, and this lower bound is already lower convex. Moreover, the estimate at
each
\[
j\in\left\{2d_N,\,2(d_N+d_{N-1}),\,\dots,\,2(d_N+\cdots+d_M)\right\}
\]
determines a vertex of the Newton polygon of $A_1$. These vertices yield the desired slopes.
	\end{proof}

\begin{remark}
	\label{rmk2}
By Corollary \ref{cor main1}, the $k$-derivative polygon $\underline{\Delta}_k$ calculates (up to adding a constant $-1$) almost all of the slopes of (the inverse of) the $\calL$-invariants associated with $p$-newforms in $S_k(\Gamma_0(Np))_{\bar r}$, except for a \emph{logarithmic error term} in terms of $k$. More precisely, the logarithmic error term has the following twofold meaning analogous to that
discussed in Remark~\ref{rmk of log bound meaning}.
\begin{enumerate}
	\item
The number of $\mathcal{L}$-invariants associated with $p$-newforms in
$S_k(\Gamma_0(Np))_{\bar r}$ whose slopes are not $k$-derivative slopes grows
at most logarithmically in $k$ by Corollary~\ref{cor of thresholds}(3). Note that, in the setting of Corollary~\ref{cor main1}, this number is  $2(d_1+\cdots+d_M)$, while in Corollary~\ref{cor of thresholds}(3) the same
quantity is written as $2m(\bar r)n_{M-1}$.
	
	\item
	The slopes of these exceptional $\mathcal{L}$-invariants have a lower
	bound $1-\M(k)$, which grows logarithmically in $k$ by
	Lemma~\ref{lemmam(k)}.
\end{enumerate}
\end{remark}
Following Remark~\ref{rmk2}, the equidistribution and integrality of the
$k$-derivative slopes yield the following two results on the equidistribution and
integrality of the slopes of the $\mathcal{L}$-invariants.
	\begin{corollary}[Equidistribution property]
		\label{dis of L-inv}
		For an integer $k\equiv k_0\pmod{p-1}$, let $\mu_{k,\bar r}$ be the uniform probability measure of the following multiset
		$$
		Y_k:=\left\{\frac{2(p+1)}{(p-1)k}\cdot v_p(\calL^{-1}_f): f\ \textrm{is a}\ p\textrm{-newform in}\ S_k(\Gamma_0(Np))_{\bar r}\right\}\subseteq (-\infty,+\infty].
		$$
		Then the measure $\mu_{k,\bar r}$ weakly converges to the uniform probability measure on the interval $[0,1]$ as $k$ tends to infinity.
	\end{corollary}
	\begin{proof}
By Remark \ref{rmk2}(2), this corollary follows from identical arguments in the proof of Theorem \ref{distribution thm for constant} (the only differences come from the constant $-1$ we mentioned in Remark~\ref{rmk2}).
	\end{proof}
	\begin{corollary}[Integrality]
		\label{cor: integrality}
		Let $k$ be an integer such that $k\equiv k_0\pmod{p-1}$. 
	Let $C(k,a)$ be the constant defined in Corollary~\ref{cor of thresholds}(3).
Then with at most $\lceil C(k,a)\cdot m(\bar r)\rceil$ exceptions, each $p$-newform $f\in S_k(\Gamma_0(Np))_{\bar r}$ satisfies the following properties.
\begin{enumerate}
	\item If $1-v_p(\calL_f)$ occurs with multiplicity one in $\underline{\Delta}_k$, then $v_p(\calL_f)\in\ZZ$.
	\item If $1-v_p(\calL_f)$ occurs with multiplicity greater than one in $\underline{\Delta}_k$, then  $v_p(\calL_f)\in\frac{k}{2}+\ZZ$.
\end{enumerate}
	\end{corollary}
	\begin{proof}
		This follows from Proposition \ref{nearStprop}, Corollary \ref{cor main1} and Remark~\ref{rmk2}(1).
	\end{proof}
In the remaining part of this section, we carry out the strategy outlined in \S~\ref{subsection two step approach} to prove Theorem \ref{prop slopes for A_1}. We point out that the proof is similar to that of Proposition~\ref{prop: estimates of X_j's in an example}.
\subsection{The lower bounds}
 We first construct lower bounds for the dual graphs $\NP^*(c^\new_i)$ from the $(k,w_*)$-newslopes for various weights $w_*$ such that $v_p(w_*-w_k)>\M(k)$. This is a generalization of Lemma \ref{lemma: lower bd in the example}.

\begin{definition}
	For $r> \M(k)$, let $\left(i,f_i(r)\right)\in\RR^2$ be the $i$-th point on the Newton polygon $\NP^\new(w_*)$ for some (or any) $w_*$ satisfying $v_p(w_*-w_k)=r$. For these weights $w_*$'s, we have $v_p(c^\new_i(w_*))\geq f_i(r)$.
	Denote by $\boldsymbol{\Gamma}(f_i)$ the graph of the function $f_i(r)$ in $\RR^2$ over the interval $\left(\M(k),+\infty\right)$. Clearly, $\boldsymbol{\Gamma}(f_i)$ lies below the dual graph $\NP^*(c^\new_i)$.
\end{definition}
\begin{lemma}
	\label{y_i function}
	We have the following descriptions for the piecewise linear functions $f_i(r)$, for $i\in\left\{1,\dots,d\right\}$, over the interval $r\in\left(\M(k),+\infty\right)$.
	First, these functions are symmetric in the sense					
	$f_{d-i}(r)=f_i(r)+(d-2i)(k-2)$ for $i\in\left\{1,\dots,\frac{d}{2}\right\}$, and $f_d(r)=d(k-2)$.
	Furthermore, for $i\in\left\{1,\dots,\frac{d}{2}\right\}$, we have
	\begin{enumerate}
		\item Consider $1\leq i\leq d_M+\cdots+d_N$. Let $t\in\{M,\dots,N\}$ be such that, with the convention $d_{N+1}=0$,
		\[
		i\in\left\{1+\sum_{m=t+1}^N d_m,\dots,\sum_{m=t}^N d_m\right\}.
		\]
	Then the graph $\boldsymbol{\Gamma}(f_i)$ is determined by the initial value
	$f_i(s_N)=(k-2)i$, together with its vertices and the slopes between consecutive
	vertices, as shown below. The numbers above the arrows indicate the
	corresponding slopes.
		\begin{equation}
			\label{5.13.2}
			\begin{aligned}
				&\left(\infty,f_i(s_N)\right)\xrightarrow{0}(s_N,f_i(s_N)) \xrightarrow{2d_N} (s_{N-1},f_{i}(s_{N-1})) \xrightarrow{2d_N+2d_{N-1}}\cdots
				\\& \xrightarrow{2d_N+\cdots+2d_{t+1}} 
				(s_{t},f_{i}(s_{t})) \xrightarrow{2i} (\M(k),f_i(\M(k))).
			\end{aligned}
		\end{equation}
 		\item For the remaining larger indices
 		$
 		i\in\left\{1+\sum_{m=M}^N d_m,\dots,\frac{d}{2}\right\}$,
 		the graph $\boldsymbol{\Gamma}(f_i)$ is given by the initial value
 		$f_i(s_N)=(k-2)i$ together with the following data
		\begin{equation}
			\label{5.14.3}
			\begin{aligned}
				&\left(\infty,f_i(s_N)\right)\xrightarrow{0}(s_N,f_i(s_N))\xrightarrow{2d_N}\cdots
				\\&\xrightarrow{2d_N+\cdots+2d_{M+1}}(s_{M},f_{i}(s_M))\xrightarrow{2d_N+\cdots+2d_{M}}(\M(k),f_i(\M(k))).
			\end{aligned}
		\end{equation}
	\end{enumerate}
	Moreover, for $i=d_N,d_N+d_{N-1},\dots,d_N+\cdots+d_M$, the graph $\boldsymbol{\Gamma}(f_i)$ coincides with the dual graph $\NP^{*}(c^\new_i)$ over the interval $$(\M(k),s_N],(\M(k),s_{N-1}],\dots,(\M(k),s_{M}],$$ respectively.
\end{lemma}
\begin{proof}
	This is a reinterpretation of Theorem \ref{propofslopes}.
	For a weight $w_*$ for which $v_p(w_*-w_k)=r$, $\frac{1}{2}f_i(r)$ is equal to the sum of the $i$ smallest $(k,w_*)$-newslopes (as we did in Example \ref{example 3}). Since these slopes can be written down explicitly by Theorem \ref{propofslopes}, this lemma follows from direct computations. In particular, for the last assertion, let
	$
	i=d_N+d_{N-1}+\cdots+d_t$
	for some $t\in\{M,\dots,N\}$. Then, for any
	$r\in(\M(k),s_t]$ and any weight $w_*$ satisfying
	$v_p(w_*-w_k)=r$, the point $
	\left(i, v_p(c_i^{\new}(w_*))\right)
	$
	is always a vertex of $\NP^{\new}(w_*)$. As a result, $v_p(c^\new_i(w_*))=f_i(r)$ for any $v_p(w_*-w_k)=r$ and $r\in(\M(k),s_t]$. This implies that the dual graph $\NP^*(c^\new_i)$ coincides with the graph $\boldsymbol{\Gamma}(f_i)$ over the interval $(\M(k),s_t]$.
\end{proof}
		\subsection{Estimates of $v_p\left(c_{i}^{\new,(\ell)}\right)$'s} Keep the notations in Notations \ref{notation: refined k-derivative slopes} and \ref{notation of Li}. For $i\in\left\{1,\dots,d\right\}$, let $\boldsymbol{k}_i=(k-2)i$. The following definition of the $d\times d$ table $\boldsymbol{M}$ is a generalization of that in Notation \ref{notation : in the example}.
	\begin{definition}
		\label{scrM notation}
		Consider the \(d\times d\) table
		\(
		\boldsymbol{M}=(\boldsymbol{m}_{i,j})_{1\le i,j\le d}
		\),
		whose entries belong to \(\{>,\geq,=\}\).
		The entries of \(\boldsymbol{M}\) are defined as follows.
		\begin{enumerate}
			\item
			The entries in the last row are all strict inequalities ``$>$".
			\item
			The entry \(\boldsymbol{m}_{i,j}\) is ``$=$"
			if \(j=2i\) and
			\[
			i\in
			\left\{
			d_N,\dots,d_N+\cdots+d_M
			\right\}
			\cup
			\left\{
			d-d_N,\dots,d-(d_N+\cdots+d_M)
			\right\}.
			\]
			
			\item
			The entry \(\boldsymbol{m}_{i,j}\) is a strict inequality
			``$>$"
			if either
			$
			i\in
			\left\{
			1,\dots,\sum_{l=M}^{N}d_l
			\right\}$ and	$j>2i$, or
			$i\in
			\left\{
			d-\sum_{l=M}^{N}d_l,\dots,d-1
			\right\}$ and 
			$j>2(d-i)$.
			\item
			All remaining entries are non-strict inequalities ``$\geq$".
		\end{enumerate}
	\end{definition}
		
		\begin{remark}
			\label{symmetry1}
		Similar to the special case discussed in Remark~\ref{rmk1}, the table
		$\boldsymbol{M}$ is \emph{symmetric} in the following sense: for every
		$i\in\left\{1,\dots,\frac{d}{2}\right\}$, the $i$-th row agrees with the
		$(d-i)$-th row. The essential difference in the general case is that the multiplicities
		$d_1,\dots,d_N$ may produce additional rows. More precisely, for each
		$t\in\{1,\dots,N\}$ with $d_t\geq 2$, the $d_t-1$ rows between the
		$(d_N+\cdots+d_{t+1})$-th row and the $(d_N+\cdots+d_t)$-th row do not appear
		in the table in Notation~\ref{notation : in the example}.

Another key observation is that $\boldsymbol{M}$ contains exactly $2(N-M+1)$
equalities, and that any column containing an equality ``$=$'', say the $j$-th
column, contains at least $j-1$ strict inequalities ``$>$''. This will help us
determine $N-M+1$ vertices of the Newton polygon of $A_1$ in the proof of
Theorem~\ref{prop slopes for A_1}, as in the proof of
Proposition~\ref{prop: estimates of X_j's in an example}.
			\end{remark}
			
Similar to Proposition \ref{prop example of estimate of F_{i,j}}, we have the following general result that provides crucial estimates of $v_p\left(c_{i}^{\new,(\ell)}\right)$'s in proving Theorem \ref{prop slopes for A_1}.
		\begin{proposition}
			\label{main estimate for scrM}
			For $1\leq i,j\leq d$, we have the following estimates of $c_{i}^{\new,(j)}$, denoted by $$v_p\left(c_{i}^{\new,(j)}\right)\approx^{\boldsymbol{m}_{i,j}}(\boldsymbol{k}_i-\boldsymbol{l}_j),$$ 
which means that the relation between $v_p\left(c_{i}^{\new,(j)}\right)$ and $\boldsymbol{k}_i-\boldsymbol{l}_j$ is given by the entry $\boldsymbol{m}_{i,j}$.
In particular, we have
			\begin{enumerate}
				\item For $i,j\in\left\{1,\dots,d\right\}$, we have $v_p\left(c_{i}^{\new,(j)}\right) \geq i(k-2) - \boldsymbol{l}_j$.
				\item
				Fix $t\in\left\{M,\dots,N\right\}$, and $j = 2\left(d_t+\cdots+d_N\right)$. Then for $i\in\left\{1,\dots,\frac{j}{2}-1\right\}$ or $i\in\left\{d-\frac{j}{2}+1,\dots,d\right\}$, we have $v_p\left(c_{i}^{\new,(j)}\right) >i(k-2) -\boldsymbol{l}_j$.
				\item Fix $t\in\left\{M,\dots,N\right\}$, and $j = 2\left(d_t+\cdots+d_N\right)$. Then for $i=\frac{j}{2}$, we have $v_p\left(c_{i}^{\new,(j)}\right)=v_p\left(c_{\frac{j}{2}}^{\new,(j)}\right) = i(k-2)-\boldsymbol{l}_j$.
			\end{enumerate}
		\end{proposition}
	
		\begin{proof}
			The proof is almost the same as that of Proposition \ref{prop example of estimate of F_{i,j}}. By symmetry discussed in Remark~\ref{symmetry1}, we only need to consider $1\leq i\leq \frac{d}{2}$. By constructions of $f_i(r)$, we have 
			\begin{equation}
				\label{eqeqeq}
				v_p\left(c_{i}^{\new,(j)}\right)\geq\nu_r(c^\new_i)-jr\geq f_i(r)-jr,\  \textrm{for any}\ j\geq0,\  \textrm{and}\  r>\M(k).
			\end{equation}
	Once again, for the inequalities in $\boldsymbol{M}$, we choose appropriate values of $r$ such that the relation between $f_i(r)-jr$ and $\boldsymbol{k}_i-\boldsymbol{l}_j$ is given by $\boldsymbol{m}_{i,j}$; for the
	equalities in $\boldsymbol{M}$, we apply Corollary \ref{analytic trick} directly. In what follows, we prove the proposition row by row, by fixing
	$i\in\left\{1,\dots,\frac{d}{2}\right\}$ and verifying the estimates given by the $i$-th row of
	$\boldsymbol{M}$.
	
	(0) Consider $i=d$. For any $r\in\left(\M(k),s_M\right)$ and any $j\in\left\{1,\dots,d\right\}$, we have $jr<\boldsymbol{l}_j$. Therefore by \eqref{eqeqeq}, we have
	\[
	v_p\left(c_{i}^{\new,(j)}\right)\geq f_i(r)-jr=d(k-2)-jr>i(k-2)-\boldsymbol{l}_j=\boldsymbol{k}_i-\boldsymbol{l}_j.
	\]
			 
(1) Consider $1\leq i\leq d_N+\cdots+d_M$.  Let $t\in\{M,\dots,N\}$ be the index such that
$
i\in\left\{1+\sum_{m=t+1}^N d_m,\dots,\sum_{m=t}^N d_m\right\}$. By Lemma~\ref{y_i function}, one checks directly that, for every
$\ell\in\{t,\dots,N\}$ and every
$
j\in\left\{1+2\sum_{m=\ell+1}^N d_m,\dots,
2\sum_{m=\ell}^N d_m\right\}$,
we have
\begin{equation}
	\label{eqeq22}
	f_i(s_{\ell})-js_{\ell}
	=i(k-2)-\boldsymbol{l}_j
	=\boldsymbol{k}_i-\boldsymbol{l}_j.
\end{equation}
\begin{itemize}
	\item We first treat the non-strict inequalities in the $i$-th row of $\boldsymbol{M}$. For
	$j\leq 2i$, taking $r=s_t$, it follows from \eqref{eqeqeq} and \eqref{eqeq22}
	that
	\begin{equation}
		v_p\left(c_i^{\new,(j)}\right)
		\geq f_i(s_t)-js_t
		=\boldsymbol{k}_i-\boldsymbol{l}_j.
	\end{equation}
	This gives the desired non-strict inequalities $\boldsymbol{m}_{i,j}$ in the $i$-th row of $\boldsymbol{M}$.
	\item We next prove the strict inequalities in the $i$-th row. Consider $j>2i$ and
	let $r\in\left(\M(k),s_M\right)$. By \eqref{eqeqeq}, \eqref{eqeq22}, and the fact that
	$f_i$ is linear of slope $2i$ over the interval $[r,s_t]$
	(Lemma~\ref{y_i function}), we have
	\begin{equation}
		\label{EQ1}
		\begin{aligned}
			v_p\left(c_i^{\new,(j)}\right)
			&\geq f_i(r)-jr  \\
			&=\left(f_i(s_t)-2is_t\right)
			+\left(f_i(r)-f_i(s_t)+2is_t-jr\right) \\
			&=\boldsymbol{k}_i-\boldsymbol{l}_{2i}-(j-2i)r.
		\end{aligned}
	\end{equation}
	Since $r\in(\M(k),s_M)$, we have
	\begin{equation}
		\label{EQ2}
		\boldsymbol{k}_i-\boldsymbol{l}_{2i}-(j-2i)r
		>
		\boldsymbol{k}_i-\boldsymbol{l}_{2i}-(j-2i)s_M
		\geq
		\boldsymbol{k}_i-\boldsymbol{l}_j,
	\end{equation}
	where the last inequality follows from the construction of
	$\boldsymbol{l}_j$. Combining \eqref{EQ1} and \eqref{EQ2}, we obtain the
	desired strict inequalities in the $i$-th row of $\boldsymbol{M}$.
	
	\item 
	It remains to treat the equality in the $i$-th row of $\boldsymbol{M}$. In
	this case, we have $i=\sum_{m=t}^N d_m$. By Lemma~\ref{y_i function}, the point
	$
	\left(s_t,f_i(s_t)\right)
	=
	\left(s_t,\nu_{s_t}(c_i^{\new})\right)
	$
	lies on the dual graph $\NP^*(c_i^{\new})$. Therefore, for $j=2i$, Corollary
	\ref{analytic trick} gives
	\begin{equation}
		v_p\left(c_i^{\new,(j)}\right)
		=
		\nu_{s_t}(c_i^{\new})-js_t
		=
		f_i(s_t)-js_t
		=
		\boldsymbol{k}_i-\boldsymbol{l}_j,
	\end{equation}
	where the last equality follows from \eqref{eqeq22}. Hence we obtain the
	desired equality $\boldsymbol{m}_{i,2i}$ in the $i$-th row of
	$\boldsymbol{M}$.
\end{itemize}

(2) Consider
$
1+d_N+\cdots+d_M \leq i \leq \frac{d}{2}$.
For each $j\in\{1,\dots,d\}$, let $t(j)\in\{1,\dots,N\}$ be the index such that
$
j\in\left\{1+2\sum_{m=t(j)+1}^N d_m,\dots,
2\sum_{m=t(j)}^N d_m\right\}$.
To obtain the analogue of \eqref{eqeq22} in this case, define
\[
r(j)=
\begin{cases}
	s_{t(j)}, & \text{if } t(j)\in\{M,\dots,N\},\\
	s_M, & \text{if } t(j)\in\{1,\dots,M-1\}.
\end{cases}
\]
By Lemma~\ref{y_i function}, one checks directly that, for every
$j\in\{1,\dots,d\}$,
\[
f_i(r(j))-jr(j)
=
i(k-2)-\boldsymbol{l}_j
=
\boldsymbol{k}_i-\boldsymbol{l}_j.
\]
Therefore, by \eqref{eqeqeq}, we obtain the desired non-strict inequalities
$
v_p\left(c_i^{\new,(j)}\right)\geq \boldsymbol{k}_i-\boldsymbol{l}_j
$ for $j\in\{1,\dots,d\}$.

The proof is complete.
		\end{proof}
		\subsection{Proof of Theorem \ref{prop slopes for A_1}}
		\label{subsection : proof of main theorem}
Finally, combining Proposition \ref{main estimate for scrM} with linear relations in Lemma \ref{lem:ZZ-linaer comb lemma} and Proposition \ref{prop: ZZ-linear comb prop}, we can finish the proof of Theorem \ref{prop slopes for A_1}. Recall $\a=p^{k-2}$.	
	\begin{proof}[Proof of Theorem \ref{prop slopes for A_1}]
		This follows from almost the same arguments in the proof of Proposition \ref{prop: estimates of X_j's in an example}. The non-strict inequalities ``$\geq$" follow from Lemma \ref{lem:ZZ-linaer comb lemma} and Proposition \ref{main estimate for scrM}.
		More precisely, each $\alpha^{-j}\tau_j$ can be expressed as a
		$\mathbb{Z}$-linear combination of
		$\{\alpha^{-t}c_t^{\new,(j)}:t=1,\dots,j\}$, and we have
		$v_p\left(c_i^{\new,(\ell)}\right)\geq i(k-2)-\boldsymbol{l}_\ell$ for all
		$i,\ell\in\{1,\dots,d\}$. Therefore we deduce that
		$
		v_p\left(\alpha^{-j}\tau_j\right)\geq -\boldsymbol{l}_j$.

		We prove the equalities as follows. For $t\in\left\{M,\dots,N\right\}$, let $j = 2(d_N + \cdots + d_t)$. Consider the first $\frac{j}{2}$ and the last $\frac{j}{2}$ equations of the $d$ linear equations associated with the index $j\in\left\{1,\dots,d\right\}$ in Proposition \ref{Cor for Simply F_i(u)}, and denote the corresponding $j \times j$ matrix by $D_d(j)$. Then we obtain linear equations in \eqref{eqeq3}, namely
		\[
			\left(	\begin{matrix}\a^{-1}c_{1}^{\new,(j)}\\
				\vdots\\
				\a^{-\frac{j}{2}}c^{\new,(j)}_{\frac{j}{2}}\\
				\a^{\frac{j}{2}-d-1}c^{\new,(j)}_{d-\frac{j}{2}+1}\\
				\vdots\\
				\a^{-d}c^{\new,(j)}_{d}\end{matrix}\right)=D_d(j)
			\cdot\left(\begin{matrix}
				\a^{-1}m_1^{+,(j)}\\
				\vdots\\
				\a^{-j}m_j^{+,(j)}
			\end{matrix}\right).
		\]
By Proposition \ref{main estimate for scrM}, 
	$v_p\left(\a^{-i}c_{i}^{\new,(j)}\right) > -\boldsymbol{l}_j$ for $i\in\left\{1,\dots,\frac{j}{2}-1\right\}\cup\left\{d-\frac{j}{2}+1,\dots,d\right\}$, and $v_p\left(\a^{-\frac{j}{2}}c_{\frac{j}{2}}^{\new,(j)}\right) = -\boldsymbol{l}_j$.
Therefore the equality $v_p(\tau_j)= j(k-2)-\boldsymbol{l}_j$ follows from Proposition \ref{prop: ZZ-linear comb prop}. This completes the proof of the theorem.
	\end{proof}

\end{document}